\documentclass[12pt]{article}

\usepackage{amssymb,amsmath}
 
\usepackage{stackengine}

\usepackage{accents} 

\let\Horig\H

\usepackage{dsfont}

\usepackage[usenames,dvipsnames,table]{xcolor}

\definecolor{ForestGreen}{RGB}{34,139,34}
\definecolor{mauve}{rgb}{0.7,0,0.43}
\definecolor{dkgreen}{rgb}{0,0.6,0}
\definecolor{darkgreen}{rgb}{0,0.6,0}
\definecolor{darkorange}{rgb}{1.0, 0.55, 0.0}
\definecolor{lightblue}{rgb}{0,0.2,0.5}
\definecolor{blue1}{rgb}{0,0.1,0.9}

\definecolor{lightblue}{rgb}{0,0.2,0.5}

\usepackage[colorlinks=true, urlcolor=blue,linkcolor=blue, citecolor=lightblue]{hyperref}

\usepackage{caption}
\usepackage{subcaption}

\usepackage{tikz}

\usetikzlibrary{automata,topaths}
\usetikzlibrary{shapes}
\usetikzlibrary{plotmarks}

\usetikzlibrary{calc,shapes}

\usetikzlibrary{positioning, fit, shapes.geometric}

\usetikzlibrary{hobby,backgrounds,calc,trees}

\pgfdeclarelayer{background}
\pgfsetlayers{background,main}

\newcommand{\hobbyconvexpath}[2]{
[   
    create hobbyhullnodes/.code={
        \global\edef\namelist{#1}
        \foreach [count=\counter] \nodename in \namelist {
            \global\edef\numberofnodes{\counter}
            \node at (\nodename)
[draw=none,name=hobbyhullnode\counter] {};
        }
        \node at (hobbyhullnode\numberofnodes)
[name=hobbyhullnode0,draw=none] {};
        \pgfmathtruncatemacro\lastnumber{\numberofnodes+1}
        \node at (hobbyhullnode1)
[name=hobbyhullnode\lastnumber,draw=none] {};
    },
    create hobbyhullnodes
]
($(hobbyhullnode1)!#2!-90:(hobbyhullnode0)$)
\pgfextra{
  \gdef\hullpath{}
\foreach [
    evaluate=\currentnode as \previousnode using int(\currentnode-1),
    evaluate=\currentnode as \nextnode using int(\currentnode+1)
    ] \currentnode in {1,...,\numberofnodes} {
    \ifnum\currentnode=1\relax
    \xdef\hullpath{([closed=true]$(hobbyhullnode\currentnode)!#2!180:(hobbyhullnode\previousnode)$)
  ..($(hobbyhullnode\nextnode)!0.5!(hobbyhullnode\currentnode)$)}
    \else
    \xdef\hullpath{\hullpath
  ..($(hobbyhullnode\currentnode)!#2!180:(hobbyhullnode\previousnode)$)
  ..($(hobbyhullnode\nextnode)!0.5!(hobbyhullnode\currentnode)$)}
    \fi
    \ifx\currentnode\numberofnodes
    \else
    \xdef\hullpath{\hullpath
  ..($(hobbyhullnode\nextnode)!#2!-90:(hobbyhullnode\currentnode)$)}
    \fi
}
}
\hullpath
}

\usepackage{pgfplots}
\usepackage{wasysym} 
\usepackage{float} 

\DeclareMathAlphabet{\eufrak}{U}{}{}{} 
\SetMathAlphabet\eufrak{normal}{U}{euf}{m}{n}
\SetMathAlphabet\eufrak{bold}{U}{euf}{b}{n}

\numberwithin{equation}{section}

\newenvironment{Proof}{\removelastskip\par\medskip
\noindent{\em Proof.} \rm}{\penalty-20\null\hfill\huge$\square$\par\medbreak}

\allowdisplaybreaks

 \def\real{{\mathord{\mathbb R}}}
 
 \def\inte{{\mathord{\mathbb N}}}
 
 \def\qu{{\mathord{\mathbb Z}}}

 \def\Var{{\mathrm{{\rm Var}}}}

 \def\real{{\mathord{{\rm I\kern-3pt R}}}}        
 
 \def\inte{{\mathord{{\rm I\kern-3pt N}}}}
 \def\sZZ{{\rm Z\kern-.45em{}Z}}

 \def\sQQ{{\kern 0.27em \vrule height1.45ex width0.03em depth0em
           \kern-0.30em \rm Q}}
 \def\qu{{\mathchoice
         {\sQQ}
         {\sQQ}
   {\kern 0.225em \vrule height1.05ex width0.025em depth0em \kern-0.25em \rm Q}
   {\kern 0.180em \vrule height0.78ex width0.020em depth0em \kern-0.20em \rm Q}
         }}
 \def\sGG{{\kern 0.27em \vrule height1.45ex width0.03em depth0em
           \kern-0.30em \rm G}}

 \newtheorem{prop}{Proposition}[section]

 \newtheorem{corollary}[prop]{Corollary}

 \def\P{{\mathord{\mathbb P}}}

\def\E{\mathop{\hbox{\rm I\kern-0.20em E}}\nolimits}

 \newcounter{hyp}
\usepackage{textcomp}

\usepackage{accsupp}    

 \newcommand{\noncopynumber}[1]{
    \BeginAccSupp{method=escape,ActualText={}}
    #1
    \EndAccSupp{}
}
 
\usepackage{listings}
\newcommand{\re}{\mathrm{e}}            

\title{\huge Asymptotic analysis of $k$-hop connectivity in the 1D unit disk random graph model}\date{\normalsize \today} 

\author{
\large 
 Nicolas Privault\thanks{nprivault@ntu.edu.sg} 
\\ 
\small
Division of Mathematical Sciences 
\\ 
\small 
School of Physical and Mathematical Sciences 
\\ 
\small
Nanyang Technological University 
\\ 
\small 
21 Nanyang Link 
\\ 
\small
Singapore 637371
}

\DeclareMathSymbol{\mlq}{\mathord}{operators}{``}
\DeclareMathSymbol{\mrq}{\mathord}{operators}{`'}

\let\savering\ring
\let\saveleftmoon\leftmoon
\let\saverightmoon\rightmoon
\let\savefullmoon\fullmoon
\let\savenewmoon\newmoon
\let\savediameter\diameter
\let\ring\relax
\let\leftmoon\relax
\let\rightmoon\relax
\let\fullmoon\relax
\let\newmoon\relax
\let\diameter\relax
\usepackage{mathabx}
\let\ring\savering
\let\leftmoon\saveleftmoon
\let\rightmoon\saverightmoon
\let\fullmoon\savefullmoon
\let\newmoon\savenewmoon
\let\diameter\savediameter

\usepackage{hhline} 

\begin{document}

\hyphenation{func-tio-nals} 
\hyphenation{Privault} 

\maketitle 

\vspace{-0.8cm}

\baselineskip0.6cm
 
\begin{abstract} 
We propose an algorithm for the closed-form recursive computation of joint moments and cumulants of all orders for $k$-hop counts in the 1D unit disk random graph model with Poisson distributed vertices. Our approach uses decompositions of $k$-hop counts into multiple Poisson stochastic integrals. As a consequence, using the Stein method we derive Berry-Esseen bounds for the asymptotic convergence of renormalized $k$-hop path counts to the normal distribution as the density of Poisson vertices tends to infinity. Computer codes for the recursive symbolic computation of moments and cumulants are provided in appendix. 
\end{abstract} 
 
 
\baselineskip0.7cm
 
\section{Introduction}
We consider the statistics and asymptotic behavior of
$k$-hop connectivity of the one-dimensional 
 unit disk random connection model with connection radius $r>0$
 on a finite interval, see \cite{drory1997}. 
 Random geometric graphs
 have the ability to model physical systems in e.g.
 wireless networks, 
 complex networks, 
 and statistical mechanics. 

 \medskip 
 
 Early results in the normal approximation of subgraph counts in random graphs
can be traced to the development of
the Erd{\H o}s-R\'enyi random graph $\mathbb{G}_n(p)$ 
 \cite{G}, \cite{ER}. 
Necessary and sufficient conditions for the asymptotic normality
 of the renormalized count of graphs in $\mathbb{G}_n(p_n)$
 that are isomorphic to a fixed graph $G$ 
 have been obtained in \cite{rucinski},
 and made more precise in \cite{BKR} by
 the derivation of explicit convergence rates in the Wasserstein distance,
 see also \cite{BarbourHolstJanson}
 for bounds on the total variation distance of subgraph counts
 to the Poisson distribution. 
 Such bounds have been strengthened in \cite{reichenbachsAoP} and \cite{roellin2}
 using the Kolmogorov distance in the case of triangle counts. 
 In \cite{PS2}, those results have been extended to any subgraph $G$ using the
 Kolmogorov distance. 

 \medskip

 In this paper, we focus on the counting of $k$-hops in the one-dimensional 
 unit disk random connection model with connection radius $r>0$
 on a finite interval, see \cite{drory1997}. 
 See also \cite{penrosebk} for the more general setting of random geometric graphs
 and \cite{wilsher} for the soft connection model.
 Here, the nodes are distributed on $[0,kr]$, $k\geq 1$, according to a
 Poisson point process $(N_t)_{t\in [0,kr]}$ with 
 intensity $\lambda (ds)$ of the form 
\begin{equation} 
\label{fjkldf1}
 \lambda (ds) = \sum_{l=1}^k {\bf 1}_{ ( (l-1)r,lr]}(s) \lambda_l((l-1)r + ds), 
\end{equation} 
 where $\lambda_1(ds)=\lambda_1(s) ds, \ldots , \lambda_k(ds)=\lambda_k(s)ds$ are
 absolutely continuous
 intensity measures on $[0,r]$, $l=1,\ldots , k$, as illustrated in the next graph. 

   \begin{center}
   \begin{tikzpicture}
   \begin{axis}[
       height=2cm,
       xmin=0, xmax=50,
    width=1\textwidth,
    axis x line=bottom,
    axis line style={-},
    hide y axis,    
    ymin=0,ymax=5,
    xticklabels={0,0,$r$,$2r$,$3r$,,,,,,$(k-1)r$,$kr$},
    scatter/classes={
        a={mark=o,draw=black}}
    ]

\addplot[scatter,only marks,
    mark size = 2pt,
    fill = black,
    scatter src=explicit symbolic]
table {
3 0 
6 0 
9 0 
12 0 
16 0 
18 0 
28 0 
35 0 
37 0 
42 0 
46 0 
58 0 
    };
\end{axis}

\end{tikzpicture}
\end{center}

   \vskip-0.3cm

   \noindent
 We are interested in the count $\sigma_k ( t )$ of 
 $k$-hops between additional nodes located respectively at $0$ and $t$ for some
 $t \in [0,kr]$, where two nodes $s,t\in [0,kr]$ are connected if and only
 if $|t-s|\leq r$. 
In this model,
the distribution of $k$-hop counts has been
expressed by a combinatorial approach in \cite{giles-privault}. 

\begin{figure}[H]
 \centering
  \includegraphics[width=0.94\linewidth]{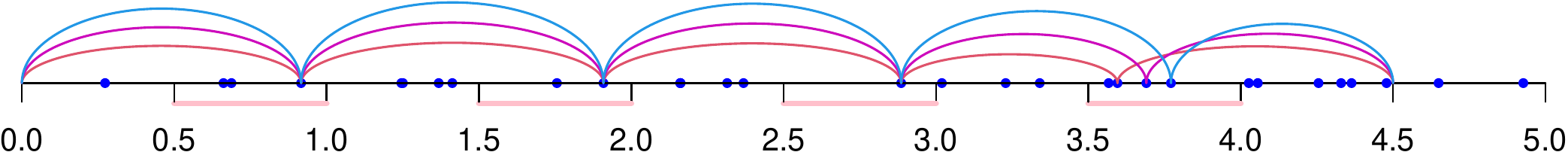}
\caption{Graphs of three $5$-hop paths linking $x=0$ to $y=4.5$ with $r=1$.} 
\label{fklds} 
\end{figure}

The moments of $k$-hop counts in the random-connection model
have been expressed in \cite{prkhp} as summations over
 non-flat partition diagrams,  
 however, those expressions are difficult to apply
 to the derivation of explicit bounds. 
 In this paper we use a different approach based on the representation of
 $k$-hop counts in terms of multiple Poisson
 stochastic integrals, which allows us to derive explicit expressions for
 moments and cumulants of all orders by recursive formulas. 

\medskip

 In Proposition~\ref{p1}
 we provide a combinatorial expression for
 the computation of the joint moments 
 of $k$-hop counts at different endpoint locations
 within $[(k-1)r,kr]$.
 This expression is then specialized to the computation
 of variance in Proposition~\ref{djklfs2}
 and Corollary~\ref{djklfs2-2}.

\medskip

 In Proposition~\ref{433},
 a recursive algorithm for the closed-form computation of
 joint moments is derived by representing
 $k$-hop counts as multiple Poisson stochastic integrals.
 A similar recursion formula is derived in Proposition~\ref{fjhkds}
 for the computation of joint cumulants,
 which yields a cumulant bound in Proposition~\ref{1fdjkl}. 

\medskip

 As a consequence, we obtain the bound   
$$ 
 \frac{ 
  c_{k,n}^{(\lambda)} (kr-t;\ldots ; kr-t)}{ (c_{k,2}^{(\lambda)}(kr-t;kr-t) )^{n/2}}
   \leq 
 (n!)^{k-2} O ( \lambda^{1-n/2} )
, \qquad n \geq 2,
$$
 on the cumulant $c_{k,n}^{(\lambda)} (kr-t;\ldots ; kr-t)$ of order $n\geq 2$ of
 $\sigma_k( t )$, for $t\in [(k-1)r, kr)$.
 Denoting by $\P_\lambda$ the distribution of the 1D unit disk graph
 with constant Poisson intensity $\lambda >0$, this implies the Berry-Esseen bound
$$ 
\sup_{x\in \real} | \P_\lambda  ( \widetilde{\sigma}_k (t) ) - \P({\cal N} \leq x ) |
\leq \frac{C(k,r)}{\sqrt{\lambda}}
$$ 
for the convergence of the renormalized $k$-hop count 
$$
\widetilde{\sigma}_k (t):= \frac{\sigma_k (t)- \E_\lambda [\sigma_k(t)]}{\sqrt{\Var_\lambda [ \sigma_k(t)]}}
$$
to the normal distribution ${\cal N}$ 
as $\lambda$ tends to infinity, see Proposition~\ref{fsklf34}.
A bound of same order is also obtained in Proposition~\ref{fsklf34} using the Wasserstein distance.

 \medskip
 
 The content of this paper can be summarized as follows.
  In Section~\ref{sec2} we show that $k$-hop counts can be represented in
terms of multiple Poisson stochastic integrals.
In Section~\ref{sec3} we specialize those expressions when the $k$-hops
are made of a single node per cell.
Section~\ref{sec4} presents moment expressions in terms of sums over non-flat
partitions based on results of \cite{prkhp}.
Sections~\ref{sec5} and \ref{sec6} develop recursive expressions for
the explicit calculation of joint moments and cumulants of any order.
In Sections~\ref{sec7} and \ref{sec8} we derive moment and cumulant bounds with
application to Berry-Esseen rates for the convergence normalized $k$-hop counts
to the normal distribution using the Stein method. 
The appendices contain specific moment and cumulant computations,
background on moment computations for Poisson point processes
 based on \cite{momentpoi, prob},  
 and Mathematica codes.

\subsubsection*{Set partitions, moments, cumulants, and M\"obius inversion}
This section gathers some preliminary facts
on the relationships between joint moments, cumulants, and sums over partitions
that will be useful in the sequel.
We let $\Pi[n]$ denote the set of partitions of $\{1,\ldots , n\}$, 
 and given a symmetric function $f(\rho ) = f(\pi_1,\ldots , \pi_l)$ where
 $\pi = \{\pi_1,\ldots , \pi_l \} \in \Pi [n]$ is a partition of $\{1,\ldots , n\}$
 of size $n\geq 1$ we will use the notation 
$$
 \sum_{\pi \in \Pi [n]} f(\pi) = \sum_{l=1}^n \sum_{\pi_1\cup \cdots \cup \pi_l} f(\pi_1,\ldots , \pi_l).
$$ 
 We will also use the M\"obius transform $\widehat{G}$ of a function $G$ on 
 partitions $\pi$ of $\{1,\ldots ,n\}$, defined as 
\begin{equation}
\label{g} 
 \widehat{G} (\sigma ) := \sum_{\pi \preceq \sigma } G(\pi ), 
 \qquad \sigma \in \Pi[n], 
\end{equation} 
where the sum \eqref{g} runs over all partitions $\pi$ of 
 $\{1,\ldots ,n \}$ that are {finer} than $\sigma$. 
 The M\"obius inversion formula,
 see e.g. \cite{rotabk} or \S~2.5 of \cite{peccatitaqqu},
 states that the function $G$ in \eqref{g}
 can be recovered from its M\"obius transform
 $\widehat{G}$ as 
\begin{equation}
\label{bpi} 
 G(\pi) = \sum_{ \sigma \preceq \pi } \mu ( \sigma , \pi ) \widehat{G} (\sigma),  
\end{equation} 
 where $\mu ( \sigma , \pi )$ 
 is the M\"obius function, with 
 $\mu ( \sigma , \widehat{\bf 1} ) = (|\sigma |-1)! (-1)^{|\sigma |}$,
 where $|\sigma |$ denotes the cardinality of the
 block $\sigma \in \Pi [n]$ and 
 $\widehat{\bf 1} : = \{\{1,\ldots , n\}\}$ is 
 the one-block partition of $\{1,\ldots , n\}$. 
 By \eqref{g} and \eqref{bpi} we also have the relation  
\begin{equation}
\label{fjlksa1} 
 G(\pi)
 = \sum_{ \sigma \preceq \pi } \mu ( \sigma , \pi )
 \sum_{\eta \preceq \sigma } G(\eta ) 
 = \sum_{ \eta \preceq \sigma \preceq \pi } 
 \mu ( \sigma , \pi ) G(\eta ), 
 \qquad \pi \in \Pi[n]. 
\end{equation} 
 Given $X=(X_1,\ldots , X_n)$ a random vector, 
 the joint {cumulants} of order $(l_1,\ldots , l_n)$
 are the coefficients 
 $\kappa \big( X_1^{l_1} ; \ldots ; X_n^{l_n} \big)$ 
 appearing in the log-moment generating (MGF) expansion 
\begin{equation} 
\nonumber 
\log \E\big[ \re^{t_1X_1+\cdots + t_n X_n} \big] 
 = 
 \sum_{l_1,\ldots , l_n\geq 1} \frac{t^{l_1}_1\cdots t^{l_n}_n}{l_1! \cdots l_n!}
 \kappa \big( X_1^{l_1} ; \ldots ; X_n^{l_n} \big)
 , 
\end{equation} 
for $(t_1,\ldots , t_n)$ in a neighborhood of zero in $\real^n$. 
 The joint moments of $(X_1,\ldots , X_n)$
 are given from its cumulants by the joint moment-cumulant relation
\begin{equation}
\nonumber 
   \E [ X_1 \cdots X_n ] 
 = 
 \sum_{\pi \in \Pi [n]} 
  \prod_{A \in \pi }
  \kappa \big( (X_i)_{i \in A} \big)
   = 
 \sum_{l=1}^n
   \sum_{\pi_1\cup \cdots \cup \pi_l = \{1,\ldots , n\}}
  \prod_{j=1}^l 
  \kappa \big( (X_i)_{i \in \pi_j} \big), 
\end{equation} 
 see Theorem~1 in \cite{elukacs} or Relation~(2.9) in \cite{mccullagh}.
 The M\"obius inversion relation \eqref{bpi} also allows us to
 recover joint cumulants from joint moments
 as 
\begin{eqnarray} 
\nonumber 
 \kappa (X_1 ; \ldots ; X_n) 
 & = & 
 \sum_{ \sigma \in \Pi [n] }
 \mu ( \sigma , \widehat{\bf 1} )
  \prod_{A \in \sigma}
  \E \Bigg[ \prod_{i\in A} X_i \Bigg]
  \\
  \label{jklsda1}
  & = &  
 \sum_{l=1}^n
  (l-1)!
  (-1)^{l-1}
  \sum_{\pi_1\cup \cdots \cup \pi_l = \{1,\ldots , n\}}
  \prod_{j=1}^l 
  \E \Bigg[ \prod_{i\in \pi_j} X_i \Bigg], 
\end{eqnarray} 
see Theorem~1 of \cite{elukacs} 
or Corollary~5.1.6 in \cite{stanley}. 
In particular,
 the cumulant of order $n\geq 1$ of a Poisson distributed random
 variable $X$ is its intensity is $\lambda >0$ for all $n\geq 1$, with 
$$ 
 \E [ X^n ] 
 = 
 \sum_{l=1}^n
   \sum_{\pi_1\cup \cdots \cup \pi_l = \{1,\ldots , n\}}
  \lambda^l 
 = 
 \sum_{l=1}^n
  S(n,l)
  \lambda^l 
$$
 where 
 $S(n,l)$ the Stirling number of the second kind, $1 \leq l \leq n$.
 For $\lambda = 1$ this yields the Bell number
$$ 
 B_n 
 = 
 \sum_{l=1}^n
  S(n,l)
$$
  which is the number of partitions of $\{1,\ldots , n \}$, i.e. the
  cardinality of $\Pi [n]$. 
\section{Multiple stochastic integral representation of $k$-hop counts} 
\label{sec2} 
   In the sequel we use the notations
   $u\wedge v := \min (u,v)$ and $u \vee v := \max (u,v)$, $u,v\geq 0$. 
   Our approach to the recursive computation of moments and cumulants
   relies on the following stochastic integral representation of 
   $k$-hop counts with respect to the Poisson process $(N_t)_{t\in \real_+}$
   with intensity \eqref{fjkldf1}. 
   \begin{prop}
   \label{fjkds}
   Let $k\geq 2$.
   The number of $k$-hops joining $0$ to $t \in [0,kr]$ can be written as
   the (non-compensated) multiple Poisson stochastic integral 
   \begin{equation}
     \label{fjkldsf} 
  \sigma_k (t)
 = \int_0^t 
  \cdots
  \int_0^t
  f_k(s_1,\ldots , s_{k-1}) dN_{s_1} \cdots dN_{s_{k-1}},
  \quad t\in [0,kr], 
\end{equation} 
  where $f_k$ is the function of $k-1$ variables defined as
  \begin{equation}
    \label{fjkl23} 
  f_k(s_1,\ldots , s_{k-1}) 
  :=
          \prod_{l=0}^{k-1}
               {\bf 1}_{ \{ s_{l+1} < s_l + r \} }, 
\end{equation}  
      $s_1,\ldots , s_{k-1}\in [0,t]$, with $s_0:=0$ and $s_k:=t$.
 \end{prop}
\begin{Proof}
   When $k=2$ we note that if a node is present at $s \in [0,r]$ it
  connects to every node inside $[0,s)$, and
    therefore it generates $\sigma_1 (s^-)$ new nodes,
    where $\sigma_1 (s^-)$ denotes the almost sure left limit
    $\sigma_1 (s^-):=\lim_{u\nearrow s} \sigma_1 (u)$. 
  When $s\in (r,2r]$, if a node is present at $s \in [0,r]$ 
  then 
  the count of $1$-hops up linking $0$ to $s-r$
  has to be deducted, which yields the evolution 
$$ 
 d\sigma_2(s) = 
 \left\{
 \begin{array}{ll}
   {\bf 1}_{[0,r]}(s) \sigma_1(s^-) dN_s, 
   \quad s\in [0,r],
   \\
   \\
   - {\bf 1}_{[r,2r]}(s) \sigma_1 (s^--r) dN_{s-r},
 \quad s\in (r,2r], 
 \end{array}
 \right.
 $$ 
 where $(T_n)_{n\geq 1}$
 denotes the jump times sequence of $(N_t)_{t\in [0,2r]}$. 
More generally, applying this argument by iterations to any $k\geq 3$ 
 leads to the system of jump stochastic differential equations 
$$ 
 d\sigma_k (s) = 
 \left\{
 \begin{array}{ll}
   {\bf 1}_{[0,(k-1)r]}(s) \sigma_{k-1} (s^-) dN_s, 
   \quad s\in [0,(k-1)r],
   \\
   \\
   - {\bf 1}_{[r,kr]}(s) \sigma_{k-1} (s^--r) dN_{s-r},
 \quad s\in ((k-1)r,kr], 
 \end{array}
 \right.
 $$ 
 or 
 $$ 
 d\sigma_k (s) = 
 {\bf 1}_{[0,(k-1)r]}(s) \sigma_{k-1} (s^-) dN_s
 - {\bf 1}_{[r,kr]}(s) \sigma^{(k-1)}_{s^--r} dN_{s-r},
 \quad s\in [0,kr], 
$$
 hence the recurrence relation
\begin{eqnarray*} 
  \sigma_k (t)
  & = & 
 \int_0^{((k-1)r)\wedge t} \sigma_{k-1} (s^-) dN_s
 - \int_r^{r\vee t} \sigma^{(k-1)}_{s^--r} dN_{s-r}
 \\
  & = & 
 \int_0^{((k-1)r)\wedge t} \sigma_{k-1} (s^-) dN_s 
 - \int_0^{0\vee (t-r)} \sigma_{k-1} (s^-) dN_s
 \\
  & = & 
 \int_{0\vee (t-r)}^{((k-1)r)\wedge t} \sigma_{k-1} (s^-) dN_s 
 \\
  & = & 
 \left\{
 \begin{array}{l}
   \displaystyle
   \int_{t-r}^{(k-1)r} \sigma_{k-1} (s^-) dN_s, \quad t\in [(k-1)r,kr], 
   \\
   \\
   \displaystyle
   \int_{t-r}^t \sigma_{k-1} (s^-) dN_s, \qquad t\in [r,(k-1)r], 
   \\
   \\
   \displaystyle
   \int_0^t \sigma_{k-1} (s^-) dN_s, \qquad t\in [0,r], 
 \end{array}
 \right.
\end{eqnarray*} 
 for $k \geq 2$. 
 Finally, by induction we obtain 
\begin{eqnarray} 
  \nonumber
  \sigma_k (t)
  & = & 
  \int_{0\vee (t-r)}^{((k-1)r)\wedge t}
  \int_{0\vee (s_{k-1}^--r)}^{((k-2)r)\wedge s_{k-1}^-} 
  \cdots
  \int_{0\vee (s_2^--r)}^{r\wedge s_2^-} dN_{s_1} \cdots dN_{s_{k-1}} 
\\
\label{fdksfa} 
  & = & 
  \int_{0\vee (t-r)}^{((k-1)r)\wedge t}
  \int_{0\vee (s_{k-1}^--r)}^{((k-2)r)\wedge t}
  \cdots
  \int_{0\vee (s_2^--r)}^{r\wedge t} dN_{s_1} \cdots dN_{s_{k-1}}, 
\end{eqnarray} 
and we conclude by letting
$$
f_k(s_1,\ldots , s_{k-1}) 
  : = 
  {\bf 1}_{ \{ t-r < s_{k-1} < (k-1)r\} } 
  {\bf 1}_{ \{ s_{k-1} -r < s_{k-2} < (k-2)r \} }
  \cdots
      {\bf 1}_{ \{ s_2 -r < s_1 < r \} } 
      ,
    $$ 
 $s_1,\ldots , s_{k-1}\in [0,t]$.
\end{Proof}
 In particular, the $2$-hop count is given by 
 \begin{equation}
   \label{fjlksd23} 
\sigma_2 (t) = \int_{0\vee (t-r)}^{r\wedge t} dN_s = N_t{\bf 1}_{[0,r]}(t) +
(N_r - N_{t-r} ){\bf 1}_{[r,2r]}(t). 
\end{equation}
 In the case of $3$-hops, we have 
\begin{align*} 
 & \sigma_3(t)
  = 
 \int_{0\vee (t-r)}^{(2r)\wedge t} \sigma_2 (s^-) dN_s 
 \\
  & =  
 \left\{
 \begin{array}{ll}
   \displaystyle
   \int_{t-r}^{2r} \sigma_2 (s^-) dN_s
   = \int_{t-r}^{2r} \int_{s^--r}^{s^-} dN_u dN_s
   = \int_{t-r}^{2r} (N_{s^-}-N_{s^--r}) dN_s,
   \quad t\in [2r,3r], 
   \\
   \\
   \displaystyle
   \int_{t-r}^t \sigma_2 (s^-) dN_s
   =
   \int_{t-r}^t N_{s^-}dN_s
   -
   \int_r^t N_{s^--r} dN_s
   \\
   \displaystyle
   \qquad
   \qquad
   \qquad
   \ \ 
   =
   \frac{1}{2} (N_t-1)N_t
   -
   \frac{1}{2} (N_{t-r}-1)N_{t-r}
   -
   \sum_{l=1+N_r}^{N_t} N_{T_l^--r} 
   , \quad t\in [r,2r], 
   \\
   \displaystyle
   \int_0^t N_{s^-} dN_s = \frac{1}{2} (N_t-1)N_t, 
   \qquad t\in [0,r]. 
 \end{array}
 \right.
\end{align*} 
 More generally, when $t\in [(l-1)r,lr]$ for some $l\in \{1,\ldots , k-1 \}$,
by \eqref{fdksfa} we have 
$$ 
  \sigma_k (t)
  = \int_{0\vee (t-r)}^t 
 \int_{0\vee (s_{k-1}^--r)}^t 
 \cdots \int_{0\vee (s_{l+1}^--r)}^t
 \int_{s_l^--r}^{(l-1)r}
  \cdots
  \int_{s_2^--r}^r dN_{s_1} \cdots dN_{s_{k-1}},
$$
 hence the identity in distribution
$$ 
 \sigma_k (t) \stackrel{d}{\simeq} \int_{r \vee t}^{t+r} 
 \int_{0\vee (s_{k-1}^--r)}^{t+r} 
 \cdots \int_{0\vee ( s_{l+1}^--r)}^{t+r}
 \int_{s_l^--r}^{lr}
 \int_{s_{l-1}^--r}^{(l-1)r}
  \cdots
  \int_{s_2^--r}^{2r} dN_{s_1} \cdots dN_{s_{k-1}}, 
$$ 
where we let $s_k:=t$. 
 The multiple compensated Poisson stochastic integral of order 
 $n\geq 1$ of a deterministic symmetric function
 $f_n \in L^2(\real_+^n, \lambda^{\otimes n})$ 
    is defined
	by 
	$$ 
	I_n (f_n) :=    
	n! 
	\int_0^\infty \int_0^{t_{n-1}} \cdots \int_0^{t_2} 
	{f}_n ( t_1,\ldots , t_n) 
	d(N_{t_1} - \lambda ( dt_1) )  
	\cdots 
	d(N_{t_n} - \lambda ( dt_n) ),  
	$$
	with the isometry property 
	\begin{equation}
	\label{isom1} 
	\mathbb{E} [ ( I_n (f_n) )^2 ] = (n!)^2
        \int_0^\infty \int_0^{t_{n-1}} \cdots \int_0^{t_2}
        f_n^2(t_1,\ldots , t_n) \lambda (dt_1)\cdots \lambda (dt_n),
        \quad n \geq 1, 
	\end{equation}
        see e.g. Propositions~2.7.1 and ~6.2.4 in \cite{privaultbk2}, 
        and references therein.
        The next corollary of Proposition~\ref{fjkds}
        gives the chaos decomposition of $k$-hop counts
        in terms of multiple Poisson stochastic integrals.
\begin{corollary}
\label{c1}
For any $t\in [0,kr]$, the number $\sigma_k (t)$ of $k$-hops linking $0$ to $t$ can be
represented as the sum of multiple
 compensated Poisson stochastic integrals 
 \begin{equation}
   \label{skt} 
 \sigma_k (t) 
 =  \frac{1}{(k-1)!} \sum_{l=0}^{k-1}
   {k-1 \choose l}
   I_l \left( {\bf 1}_{\{ * \in [0,t]^l\}}
   \int_0^t \cdots \int_0^t 
   \widetilde{f}_k (*,s_{l+1},\ldots , s_{k-1} )
   ds_{l+1}\cdots ds_{k-1}
   \right), 
\end{equation} 
 where $\widetilde{f}_k$ is the symmetrization in $k-1$ variables of the function
 $f_k$ defined in \eqref{fjkl23}. 
\end{corollary}
\begin{Proof}
  This is a direct consequence of
  Proposition~\ref{fjkds} and the binomial theorem applied
  to $(dN_{t_1} - \lambda (dt_1) )\cdots (dN_{t_{k-1}} - \lambda (dt_{k-1}) )$. 
\end{Proof}
\section{Single node per cell} 
\label{sec3}
 In the case where $k$-hop paths are constrained to have a single
node per cell $[(l-1)r,lr]$, $l=1,\ldots , k$,
we must have $t\in [(k-1)r,k]$, and $f_k(s_1,\ldots , s_{k-1}) =0$ unless
$$
 (s_1,\ldots , s_{k-1})
 \in [0,r] \times
     [r, 2r] \times
     \cdots
     \times
     [(k-2)r , (k-1)r]. 
$$ 
\subsubsection*{Multiple Poisson stochastic integral expression}
     In particular, when $t\in [(k-1)r,k]$ with $k\geq 2$,
     Relation~\eqref{fdksfa} yields the identity in distribution 
\begin{eqnarray}
  \nonumber 
      \sigma_k (t)
   & = & \int_{t-r}^{(k-1)r}
    \int_{s_{k-1}^--r}^{(k-2)r}
  \cdots
  \int_{s_2^--r}^r dN_{s_1} \cdots dN_{s_{k-1}}
  \\
  \label{fjdksdfa}
  & \stackrel{d}{\simeq} & 
 \int_t^{kr}
 \int_{s_{k-1}^--r}^{(k-1)r}
  \cdots
  \int_{s_2^--r}^{2r} dN_{s_1} \cdots dN_{s_{k-1}}
, 
\end{eqnarray} 
 where $s_k:=t$, 
 hence the following proposition. 

\begin{prop}
\label{fjkldsa} 
   Let $k\geq 2$.
 For $\tau \in [0,r]$ we have the identity in distribution 
 \begin{equation}
   \label{fjld3} 
 \sigma_k (kr-\tau) \stackrel{d}{\simeq} \int_0^\tau \int_0^{s_{k-1}^-} \cdots \int_0^{s_2^-} dN^{(1)}_{s_1} \cdots dN^{(k-1)}_{s_{k-1}},
\end{equation}
where $\big(N^{(l)}_s\big)_{s\in \real_+}$ is a family of
independent Poisson processes with respective intensities
 $\lambda_l(ds )=\lambda_l(s )ds$, $l=1,\ldots , k-1$. 
\end{prop}
\begin{Proof}
 When $t\in [(k-1)r,kr]$, by \eqref{fjdksdfa} we have the identity in distribution 
\begin{eqnarray*} 
  \sigma_k (t)
  & \stackrel{d}{\simeq} & 
  \int_t^{kr} 
  \int_{s_{k-1}^--r}^{(k-1)r}
  \cdots
  \int_{s_2^--r}^{2r} dN_{s_1} \cdots dN_{s_{k-1}}
  \\
  & = & 
 (-1)^{k-1} \int_0^{kr-t} 
  \int_r^{r+u_{k-1}^-}
  \cdots
  \int_{(k-2)r}^{r + u_2^-} dN_{kr-u_1} \cdots dN_{kr-u_{k-1}}
\end{eqnarray*} 
where we let $u_i:=kr-s_i$, $i=1,\ldots , k-1$, hence
 the identity in distribution 
$$
  \sigma_k (t)
  \stackrel{d}{\simeq} 
  \int_0^{kr-t} 
  \int_r^{r+v_{k-1}^-}
  \cdots
  \int_{(k-2)r}^{r+v_2^-} dN_{v_1} \cdots dN_{v_{k-1}}. 
$$
We note that
in the above integral we have $(l-1)r<v_l\leq l r$ for $l=1,\ldots ,k-1$,
hence the integration intervals are disjoint and the Poisson process
samples $(N_{v_i})$, $i=1,\ldots , k-1$, are independent on
their integration intervals, which yields the identity in distribution 
$$ 
\sigma_k(t) \stackrel{d}{\simeq} 
 \int_0^{kr-t} \int_0^{s_{k-1}^-} \cdots \int_0^{s_2^-} dN^{(1)}_{s_1} \cdots dN^{(k-1)}_{s_{k-1}},
   \qquad
 (k-1)r \leq t \leq kr. 
$$ 
\end{Proof} 
\subsubsection*{$U$-Statistics formulation} 
\noindent
 As noted in e.g. \cite{giles-privault}, when $\tau \in [0,r]$,
any node contributing to a $k$-hop path linking $x_0:=0$ to $x_{k+1}:=kr-\tau$
must belong to one of the {lenses} pictured in pink in Figure~\ref{fklds},
and defined as the intervals
$$
L_j := [jr-\tau,jr] = [0,\tau] + jr - \tau, \qquad j=1,\ldots , k-1, 
$$ 
 of identical length $\tau$.
 In particular, 
 any $k$-hop path linking $x_0:=0$ to $x_{k+1}:=kr-\tau$
 should have a single
 node per cell $[(j-1)r,jr]$, $j=1,\ldots , k$,
 hence it must be realized using a sequence 
 $(x_1, \ldots , x_{k-1})$ of nodes such that
 $$
 x_{i+1} < x_i + r, \qquad i = 0,1,\dots, k,
 $$
 with $x_0:=0$ and $x_k:=kr-\tau$.
 Therefore, any 
 $k$-hop path $(x_1,\ldots , x_{k-1})$ can be mapped to a
 sequence $(y_1,\ldots , y_{k-1}) \in [0,\tau ]^{k-1}$
 by the relation 
 $$
 y_j:=x_j - (jr-\tau), \qquad j =1,\ldots , k-1, 
 $$ 
 with $y_1 > \cdots > y_{k-1}$. 
Based on the above description, we can model the random graph using
a Poisson point process $\omega $ 
on $X:=[0,r]\times \{1,\ldots , k-1\}$, with intensity
$\mu$ of the form 
$$
\mu ( ds,\{i\}) := \lambda_l (ds )= \lambda_l (s )ds, \qquad l = 1,\ldots , k-1, 
$$ 
 and \eqref{fjld3} can be rewritten as in the next proposition.
\begin{prop}
 When $\tau \in [0,r]$, the count $\sigma_k (kr-\tau)$ of $k$-hop paths can be 
 represented as the $U$-statistics  
 \begin{equation}
   \label{fjda} 
\sigma_k (kr-\tau) 
= 
 \sum_{((x_1,l_1),\ldots , ((x_{k-1},l_{k-1}) ) \in \omega^{k-1}
   \atop (x_i,l_i)\not=(x_j,l_j), 1\leq i\not= j \leq d}
 f_\tau (x_1,l_1;\ldots ; x_{k-1},l_{k-1} )
\end{equation} 
 of order $k-1$, where $f_\tau :( [0,r] \times \{1,\ldots , k-1\})^{k-1} \to \{0,1\}$ is
 the function of $k-1$ variables in $[0,r]\times \{1,\ldots , k-1\}$
 given by 
$$
f_\tau (x_1,l_1;\ldots ; x_{k-1},l_{k-1})
=
  \prod_{i=0}^{k-1} {\bf 1}_{\{x_i<x_{i+1}, \ l_i<l_{i+1} \}}
=
  {\bf1}_{\{ l_1=1,\ldots , l_{k-1}=k-1\}} 
  {\bf1}_{\{ 0 < x_1 < \cdots < x_{k-1} < \tau \}}
$$ 
  with $(x_0,l_0):=(0,0)$
  and $(x_k,l_k):=(t,k)$.
\end{prop} 

\section{Joint moments of $k$-hop counts} 
\label{sec4}
 Proposition~\ref{p1} provides a combinatorial expression for the joint 
 moments of $(\sigma_k(kr-\tau_1),\ldots ,\sigma_k(kr-\tau_n))$ for any 
 $\tau_1 , \ldots , \tau_n \in [0,r]$, 
 using sums over partitions of $\{1,\ldots , n\}$. 
\begin{prop}
  \label{p1}
  Let $n \geq 1$. 
  For any $\tau_1,\ldots , \tau_n \in [0,r]$,
     letting 
   $\widehat{\tau}_\pi := \min_{i\in \pi} \tau_i$ for $\pi \subset \{1,\ldots , n\}$,
 we have 
\begin{align} 
\nonumber 
& \E  [ \sigma_k(kr-\tau_1) \cdots  \sigma_k(kr-\tau_n)  ] 
\\
\label{djkls12}
& =  
 \sum_{\pi^1, \ldots , \pi^{k-1} \in \Pi [n] }
 \int_{\prod_{l=1}^{k-1} \prod_{j=1}^{|\pi^l|} [0,\widehat{\tau}_{\pi^l_j}]} 
 \prod_{1 \leq l < k
   \atop
   { 
     1 \leq j \leq |\pi^l |
     \atop i\in \pi^l_j}
   }
        {\bf 1}_{\{
          z^1_{\zeta^1_i } < \cdots < z^{k-1}_{\zeta^{k-1}_i } \} } 
 \lambda_1(dz^1_{\pi^1} ) \cdots \lambda_{k-1} (dz^{k-1}_{\pi^{k-1}} )
 ,
\end{align} 
 where
$\zeta^j_i$ denotes the block of $\pi^j$ that contains the index $i\in \{1,\ldots , n\}$,  and $dz^j_{\pi^j} := (dz^j_i )_{ i \in \pi^j}$,
 $j=1,\ldots , k-1$. 
\end{prop}
\begin{Proof}
  For any $\tau_1,\ldots , \tau_n \in [0,r]$,
  by \eqref{fjda} and Corollary~\ref{c1-a} we have 
\begin{align*} 
\nonumber 
& \E  [ \sigma_k(kr-\tau_1) \cdots  \sigma_k(kr-\tau_n)  ] 
\\
& = 
 \sum_{
   \substack{
     \pi \in \Pi [n\times (k-1)] 
     \\
     \pi \wedge \rho = \hat{0}
   }}
   \int_{[0,r]^{|\pi|}}
   \sum_{1\leq l_q \leq k-1 \atop 1 \leq q \leq |\pi|}
\prod_{1 \leq j \leq |\pi | \atop 
   i\in \pi_j}
   f_{\tau_i} \big(z_{\zeta^\pi_{i,1}},l_{\zeta^\pi_{i,1}}; \ldots ; z_{\zeta^\pi_{i,k-1}},l_{\zeta^\pi_{i,k-1}}\big) 
   \lambda_{\bar{\pi}_1} ( d z_1 ) 
   \cdots
   \lambda_{\bar{\pi}_{|\pi|}} ( d z_{|\pi|} ) 
, 
\end{align*} 
 where
 $\widehat{\bf 0} : = \{\{1\},\ldots , \{n\}\}$ is 
the $n$-block partition of $\{1,\ldots , n\}$,
 $\bar{\pi}_i$ denotes the index $j\in \{1,\ldots , k-1\}$ of the unique block
$\eta_j = ((i,j))_{i=1,\ldots , n}$ containing $\pi_i$, $i=1,\ldots , |\pi|$,
and the sum is taken over non-flat 
 partitions $\pi$ in
${\rm NC} [n\times (k-1)]$ that are non-crossing in the sense that
if $(k,l)$ and $(k',l')$ belong to a same block of $\pi$ then we should have $l=l'$.
 This yields
\begin{align*} 
\nonumber 
& \E  [ \sigma_k(kr-\tau_1) \cdots  \sigma_k(kr-\tau_n)  ] 
\\
& = 
 \hskip-0.4cm
  \sum_{
   \substack{
     \pi \in \Pi [n\times (k-1)] 
     \\
     \pi \wedge \rho = \hat{0}
   }}
   \int_0^{\widehat{\tau}_{\pi_1}} 
   \cdots
   \int_0^{\widehat{\tau}_{\pi_{|\pi|}}}
   \hskip-0.4cm
   \sum_{1\leq l_q \leq k-1 \atop 1 \leq q \leq |\pi|}
   \prod_{1 \leq j \leq |\pi | \atop 
   i\in \pi_j}
   f_{\tau_i} \big(z_{\zeta^\pi_{i,1}},l_{\zeta^\pi_{i,1}}; \ldots ; z_{\zeta^\pi_{i,k-1}},l_{\zeta^\pi_{i,k-1}}\big) 
   \lambda_{\bar{\pi}_1} ( d z_1 ) 
   \cdots
   \lambda_{\bar{\pi}_{|\pi|}} ( d z_{|\pi|} ) 
\\
& = 
  \sum_{
   \substack{
     \pi \in \Pi [n\times (k-1)] 
     \\
     \pi \wedge \rho = \hat{0}
   }}
     {\bf 1}_{\{
       \bar{\pi}_1 \leq \cdots \leq \bar{\pi}_{|\pi|}
       \}
     }
   \sum_{l_1\leq \cdots \leq l_{|\pi|}} 
     \\
  & \qquad \quad 
   \int_0^{\widehat{\tau}_{\pi_1}} 
   \cdots
   \int_0^{\widehat{\tau}_{\pi_{|\pi|}}}
  \prod_{1 \leq j \leq |\pi | \atop 
   i\in \pi_j}
   f_{\tau_i} \big(z_{\zeta^\pi_{i,1}},l_{\zeta^\pi_{i,1}}; \ldots ; z_{\zeta^\pi_{i,k-1}},l_{\zeta^\pi_{i,k-1}}\big) 
   \lambda_{\bar{\pi}_1} ( d z_1 ) 
   \cdots
   \lambda_{\bar{\pi}_{|\pi|}} ( d z_{|\pi|} ) 
 \\
  & =  
 \sum_{
   \substack{
     \pi \in {\rm NC} [n\times (k-1)] 
     \\
     \pi \wedge \rho = \hat{0}
   }}
     {\bf 1}_{\{
       \bar{\pi}_1 \leq \cdots \leq \bar{\pi}_{|\pi|}
       \}
     }
    \int_0^{\widehat{\tau}_{\pi_1}} 
   \cdots
   \int_0^{\widehat{\tau}_{\pi_{|\pi|}}}
   \prod_{1 \leq j \leq |\pi |
     \atop i\in \pi_j}
      {\bf 1}_{\{ z_{\zeta^\pi_{i,1}} < \cdots < z_{\zeta^\pi_{i,k-1}} \} } 
 \lambda_{\bar{\pi}_1} (dz_1 ) \cdots \lambda_{\bar{\pi}_{|\pi|}} (dz_{|\pi |} ). 
\end{align*} 
 We conclude to \eqref{djkls12} by noting that
 any non-flat and non-crossing partition $\pi$ in
${\rm NC} [n\times (k-1)]$ can be written as
$$
\pi = \{\pi_1,\ldots , \pi_{|\pi|} \} = \bigcup_{l=1}^{k-1} \pi^l,
$$
where $\pi^l\in \Pi [n]$ is a partition of $\{1,\ldots , n\}$
for every $l=1,\ldots , k-1$.
\end{Proof}
Next, we present the application of Proposition~\ref{p1}
in the particular cases of first and second moments. 
\subsubsection*{First moment}
When $n=1$ there is only one non-flat and non-crossing partition
of $1 \times (k-1)$, which is given as 
$\rho = \{ \{(1,1)\}, \ldots , \{(1,k-1) \}\}$ and
can be represented as follows for $k=9$:
\\ 

\hskip-0.5cm
\begin{tikzpicture}[every node/.style={draw}] 
\centering
\node (l1) at (2,0) [circle] {};
\node (l2) at (4,0) [circle] {};
\node (l3) at (6,0) [circle] {};
\node (l4) at (8,0) [circle] {};
\node (l5) at (10,0) [circle] {};
\node (l6) at (12,0) [circle] {};
\node (l7) at (14,0) [circle] {};
\node (l8) at (16,0) [circle] {};

\draw let \p1=(l1), \p2=(l1), \n1={atan2(\y2-\y1,\x2-\x1)}, \n2={veclen(\y2-\y1,\x2-\x1)}
  in ($ (l1)!0.5!(l1) $) ellipse [x radius=\n2/2+12pt, y radius=0.4cm,rotate=0-\n1];
\draw let \p1=(l2), \p2=(l2), \n1={atan2(\y2-\y1,\x2-\x1)}, \n2={veclen(\y2-\y1,\x2-\x1)}
  in ($ (l2)!0.5!(l2) $) ellipse [x radius=\n2/2+12pt, y radius=0.4cm,rotate=0-\n1];
\draw let \p1=(l3), \p2=(l3), \n1={atan2(\y2-\y1,\x2-\x1)}, \n2={veclen(\y2-\y1,\x2-\x1)}
  in ($ (l3)!0.5!(l3) $) ellipse [x radius=\n2/2+12pt, y radius=0.5cm,rotate=0-\n1];
\draw let \p1=(l4), \p2=(l4), \n1={atan2(\y2-\y1,\x2-\x1)}, \n2={veclen(\y2-\y1,\x2-\x1)}
  in ($ (l4)!0.5!(l4) $) ellipse [x radius=\n2/2+12pt, y radius=0.5cm,rotate=0-\n1];
\draw let \p1=(l5), \p2=(l5), \n1={atan2(\y2-\y1,\x2-\x1)}, \n2={veclen(\y2-\y1,\x2-\x1)}
  in ($ (l5)!0.5!(l5) $) ellipse [x radius=\n2/2+12pt, y radius=0.5cm,rotate=0-\n1];
\draw let \p1=(l6), \p2=(l6), \n1={atan2(\y2-\y1,\x2-\x1)}, \n2={veclen(\y2-\y1,\x2-\x1)}
  in ($ (l6)!0.5!(l6) $) ellipse [x radius=\n2/2+12pt, y radius=0.5cm,rotate=0-\n1];
\draw let \p1=(l7), \p2=(l7), \n1={atan2(\y2-\y1,\x2-\x1)}, \n2={veclen(\y2-\y1,\x2-\x1)}
  in ($ (l7)!0.5!(l7) $) ellipse [x radius=\n2/2+12pt, y radius=0.5cm,rotate=0-\n1];
\draw let \p1=(l8), \p2=(l8), \n1={atan2(\y2-\y1,\x2-\x1)}, \n2={veclen(\y2-\y1,\x2-\x1)}
  in ($ (l8)!0.5!(l8) $) ellipse [x radius=\n2/2+12pt, y radius=0.5cm,rotate=0-\n1];
\end{tikzpicture}

\noindent
 This yields 
 \begin{eqnarray}
   \nonumber 
   \E  [ \sigma_k(kr-\tau) ] 
  & = & 
    \int_0^t \cdots \int_0^t  
 {\bf 1}_{\{ z_1 < \cdots < z_{k-1} \} } 
 \lambda_1 (dz_1 ) \cdots \lambda_{k-1} (dz_{k-1} )
 \\
 \label{fjkslff}
  & = & 
    \int_0^t \int_0^{z_{k-1}} \cdots \int_0^{z_2} 
 \lambda_1 (dz_1 ) \cdots \lambda_{k-1} (dz_{k-1} ),
\end{eqnarray} 
 and when $\lambda_i (s)$ is the constant density $\lambda_i >0$ on cell $i$,
 $i=1,\ldots , k-1$, we find
 \begin{equation}
\nonumber 
          \E [  \sigma_k (kr-\tau) ] 
 = \lambda_1 \cdots \lambda_{k-1} \frac{t^{k-1}}{(k-1)!}. 
\end{equation} 

\subsubsection*{Second moment}
\noindent
 When $n=2$ the count of blocks
of non-flat and non-crossing partitions of $[2\times (k-1)]$
ranges from $k-1$ to $2k-2$, each block has size either one or two,
as in the following example with $k=9$.
\\ 

\hskip-0.3cm
\begin{tikzpicture}[every node/.style={draw}] 
\centering
\node (l1) at (2,0) [circle] {};
\node (l2) at (4,0) [circle] {};
\node (l3) at (6,0) [circle] {};
\node (l4) at (8,0) [circle] {};
\node (l5) at (10,0) [circle] {};
\node (l6) at (12,0) [circle] {};
\node (l7) at (14,0) [circle] {};
\node (l8) at (16,0) [circle] {};

\node (m1) at (2,1) [circle] {};
\node (m2) at (4,1) [circle] {};
\node (m3) at (6,1) [circle] {};
\node (m4) at (8,1) [circle] {};
\node (m5) at (10,1) [circle] {};
\node (m6) at (12,1) [circle] {};
\node (m7) at (14,1) [circle] {};
\node (m8) at (16,1) [circle] {};

\draw let \p1=(l1), \p2=(l1), \n1={atan2(\y2-\y1,\x2-\x1)}, \n2={veclen(\y2-\y1,\x2-\x1)}
  in ($ (l1)!0.5!(l1) $) ellipse [x radius=\n2/2+12pt, y radius=0.4cm,rotate=0-\n1];
\draw let \p1=(m1), \p2=(m1), \n1={atan2(\y2-\y1,\x2-\x1)}, \n2={veclen(\y2-\y1,\x2-\x1)}
  in ($ (m1)!0.5!(m1) $) ellipse [x radius=\n2/2+12pt, y radius=0.4cm,rotate=0-\n1];
\draw let \p1=(l2), \p2=(m2), \n1={atan2(\y2-\y1,\x2-\x1)}, \n2={veclen(\y2-\y1,\x2-\x1)}
  in ($ (l2)!0.5!(m2) $) ellipse [x radius=\n2/2+12pt, y radius=0.5cm,rotate=0-\n1];
\draw let \p1=(l3), \p2=(m3), \n1={atan2(\y2-\y1,\x2-\x1)}, \n2={veclen(\y2-\y1,\x2-\x1)}
  in ($ (l3)!0.5!(m3) $) ellipse [x radius=\n2/2+12pt, y radius=0.5cm,rotate=0-\n1];
\draw let \p1=(l4), \p2=(l4), \n1={atan2(\y2-\y1,\x2-\x1)}, \n2={veclen(\y2-\y1,\x2-\x1)}
  in ($ (l4)!0.5!(l4) $) ellipse [x radius=\n2/2+12pt, y radius=0.4cm,rotate=0-\n1];
\draw let \p1=(m4), \p2=(m4), \n1={atan2(\y2-\y1,\x2-\x1)}, \n2={veclen(\y2-\y1,\x2-\x1)}
  in ($ (m4)!0.5!(m4) $) ellipse [x radius=\n2/2+12pt, y radius=0.4cm,rotate=0-\n1];
\draw let \p1=(l5), \p2=(l5), \n1={atan2(\y2-\y1,\x2-\x1)}, \n2={veclen(\y2-\y1,\x2-\x1)}
  in ($ (l5)!0.5!(l5) $) ellipse [x radius=\n2/2+12pt, y radius=0.4cm,rotate=0-\n1];
\draw let \p1=(m5), \p2=(m5), \n1={atan2(\y2-\y1,\x2-\x1)}, \n2={veclen(\y2-\y1,\x2-\x1)}
  in ($ (m5)!0.5!(m5) $) ellipse [x radius=\n2/2+12pt, y radius=0.4cm,rotate=0-\n1];
\draw let \p1=(l6), \p2=(m6), \n1={atan2(\y2-\y1,\x2-\x1)}, \n2={veclen(\y2-\y1,\x2-\x1)}
  in ($ (l6)!0.5!(m6) $) ellipse [x radius=\n2/2+12pt, y radius=0.4cm,rotate=0-\n1];
\draw let \p1=(l7), \p2=(m7), \n1={atan2(\y2-\y1,\x2-\x1)}, \n2={veclen(\y2-\y1,\x2-\x1)}
  in ($ (l7)!0.5!(m7) $) ellipse [x radius=\n2/2+12pt, y radius=0.4cm,rotate=0-\n1];
\draw let \p1=(l8), \p2=(m8), \n1={atan2(\y2-\y1,\x2-\x1)}, \n2={veclen(\y2-\y1,\x2-\x1)}
  in ($ (l8)!0.5!(m8) $) ellipse [x radius=\n2/2+12pt, y radius=0.4cm,rotate=0-\n1];
\end{tikzpicture}

\noindent 
 As a consequence, in the next proposition we obtain the 
 second moment of the count $\sigma_k(t)$ of $k$-hop paths.  
 Higher cumulants and moments of $\sigma_k(t)$ may also 
 be computed by this method 
 using Corollary~\ref{c1} above and Corollary~7.4.1 of
 \cite{peccatitaqqu}. 
\begin{prop} 
\label{djklfs2}
The variance 
of the $k$-hop path count $\sigma_k(kr-\tau)$, $\tau\in [0,r]$, is given by
 $$ 
 \Var_\lambda [ \sigma_k (kr-\tau)]
  =
  \sum_{l=1}^{k-1}
  \tau^{2k-2-l}
  \frac{  \lambda_1^2\cdots \lambda_{k-1}^2}{(2k-2-l)!} 
  \sum_{ j_0 + \cdots + j_l = k-1-l \atop j_0,\ldots , j_l\geq 0}
  \prod_{q=1}^l \frac{1}{\lambda_{j_0+\cdots + j_{q-1}+q}}
  \prod_{p=0}^l {2j_p \choose j_p} 
. 
$$
\end{prop}
\begin{Proof}
 We apply Proposition~\ref{p1}
 by noting that the blocks of size one are in even number, and 
 denoting by $i_1,\ldots , i_l$ their locations with  
 $i_1=2,i_2=3,i_3=6,i_4=7,i_5=8$ in the above example,
 when $t=t_1=t_2$,
 and letting $z_0:=0$ and $z_k:=\tau$, we obtain 
\begin{align*} 
 \E [\sigma_k^2 (kr-\tau) ] 
 & = 
  \lambda_1\cdots \lambda_{k-1}
  \sum_{l=0}^{k-1}
  \sum_{
  0=i_0<i_1< \cdots < i_l <i_{l+1}= k}
 \left( \prod_{1\leq q < k \atop
   q\notin \{i_1,\ldots ,i_l\}} \hskip-0.3cm \lambda_q
 \right)
 \\
  & \qquad \times \int_0^\tau
\int_0^{z_{i_l}}
\cdots
\int_0^{z_{i_2}}
\prod_{p=0}^l
\left( \frac{ (z_{i_{p+1}}-z_{i_p})^{(i_{p+1}-i_p-1)} }{(i_{p+1}-i_p-1)!}
\right)^2
dz_{i_1}\cdots dz_{i_l}
, 
\end{align*} 
      where
      we let $z_0:=0$ and $z_k:=\tau$,
 To conclude, we check that 
\begin{eqnarray} 
\nonumber 
\lefteqn{ 
 \! \! \! \! \! \! \! \! \! \! \! \! \! \! \! \! \int_0^{a_k}
\int_0^{z_{i_l}}
\cdots
\int_0^{z_{i_2}}
\prod_{p=0}^l
\frac{ (z_{i_{p+1}}-z_{i_p})^{2(i_{p+1}-i_p-1)} }{((i_{p+1}-i_p-1)!)^2}
dz_{i_1}\cdots dz_{i_l}
}
\\
\nonumber
& 
=
&  \frac{1}{((i_1-1)!)^2}
    \prod_{p=1}^l
\int_0^1
\frac{ (1-y)^{2(i_{p+1}-i_p-1)}
  y^{2i_p-p-1}
}{((i_{p+1}-i_p-1)!)^2}
dy
\\
\nonumber
& 
= & 
 \frac{1}{((i_1-1)!)^2}
\prod_{p=1}^l
\frac{B(2i_p-p,
2(i_{p+1}-i_p)-1)
  }{((i_{p+1}-i_p-1)!)^2}
\\
\nonumber
& 
= & 
  \frac{1}{(2k-2-l)!} 
 \prod_{p=0}^l
\frac{
(2(i_{p+1}-i_p-1))!}
     {
     ((i_{p+1}-i_p-1)!)^2}
     ,
\end{eqnarray} 
 where
 $$
 B(x,y) := \int_0^1 t^{x-1}(1-t)^{y-1} dt = \frac{(x-1)!(y-1)!}{(x+y-1)!}, 
 \qquad  x, y>0,
$$
 is the beta function, and
 $2^{-(i_{l+1}-i_l-1)}(2(i_{l+1}-i_l-1))! / (i_{l+1}-i_l-1)!$
 is the number of pair-partitions of $2(i_{l+1}-i_l-1)$. 
\end{Proof}
Alternatively, Proposition~\ref{djklfs2}
can be proved as a consequence of Corollary~\ref{c1},
and the It\^o isometry for multiple Poisson stochastic integrals. 
 For this, we can use the expression 
\begin{align*} 
& \sigma_k (kr-\tau ) = 
  \int_{(k-1)r - \tau}^{(k-1)r} 
  \int_{s_{k-1}^--r}^{(k-2)r}
  \cdots
  \int_{s_2^--r}^r dN_{s_1} \cdots dN_{s_{k-1}}
\\
 & = 
 \sum_{l=0}^{k-1}
 \sum_{
   0=i_0<i_1< \cdots < i_l <i_{l+1}= k}
 \int_0^{(k-1)r} 
  \cdots
 \int_0^{(k-1)r} 
  f_k(s_1,\ldots ,s_{k-1})
   \prod_{1\leq q \leq d \atop
     q\notin \{i_1,\ldots ,i_l\}}
   \hskip-0.3cm
   (\lambda_q ds_q) 
   \prod_{p=1}^l 
  (dN_{s_{i_p}} -\lambda_{i_p}ds_{i_p})  
\end{align*} 
 that follows from Corollary~\ref{c1} and
 the isometry and orthogonality property \eqref{isom1} of multiple Poisson
 stochastic integrals, to show that 
\begin{align*} 
& \E [\sigma_k^2 (kr-\tau) ] 
  \\
   & = 
 \sum_{l=0}^{k-1}
  \sum_{
  0=i_0<i_1< \cdots < i_l <i_{l+1}= k}
  \int_{[0,(k-1)r]^l}
  \bigg( 
  \int_{[0,(k-1)r]^{d-l}} 
  f_k(s_1,\ldots ,s_{k-1})
   \prod_{1\leq q < k \atop
     q\notin \{i_1,\ldots ,i_l\}}
   \hskip-0.3cm
  (\lambda_q ds_q)
  \bigg)^2
   \prod_{p=1}^l 
  (\lambda_{i_p}ds_{i_p}). 
\end{align*} 
 The following table provides variance formulas of $k$-hop counts
  computed from Proposition~\ref{djklfs2} by taking $\tau=1$ for simplicity.
\begin{table}[H]
\centering 
      \begin{tabular}{c|c|c|c|}
  \hhline{~-}
        \multicolumn{1}{c|}{} & 
           \addstackgap[3pt]{Variance} \cellcolor{gray!25}  
                                                  \\ 
                                                  \cline{1-2}
                                                  \multicolumn{1}{|c|}{$2$-hops} & \small $\displaystyle \lambda_1$
\\ 
\cline{1-2} 
\multicolumn{1}{|c|}{\addstackgap[10pt]{$3$-hops}} & \small $\displaystyle \frac{\lambda_1\lambda_2}{2} + 2\frac{\lambda_1^2\lambda_2+\lambda_1\lambda_2^2}{3!} 
$ 
\\ 
\cline{1-2} 
\multicolumn{1}{|c|}{\addstackgap[10pt]{$4$-hops}} & \small $\displaystyle \frac{\lambda_1 \lambda_2 \lambda_3}{3!} + 2 \frac{\lambda_1^2 \lambda_2 \lambda_3 + \lambda_1 \lambda_2^2 \lambda_3 + \lambda_1 \lambda_2 \lambda_3^2}{4!} + \frac{4 \lambda_1^2 \lambda_2 \lambda_3^2 + 6 \lambda_1^2 \lambda_2^2 \lambda_3 + 4 \lambda_1 \lambda_2^2 \lambda_3^2}{5!}$
\\ 
\cline{1-2}
\end{tabular} 
\caption{Variances of $k$-hop counts.} 
\end{table} 

\vspace{-0.3cm}
  
\noindent 
 In case the Poisson intensities are identical on all cells,  
 we obtain the following result.
\begin{corollary} 
\label{djklfs2-2}
Assume that $\lambda = \lambda_1 = \cdots = \lambda_{k-1}$
and let $\tau \in [0,r]$. 
 Then, the variance 
 of the $k$-hop path count $\sigma_k (kr-\tau)$ is given by
\begin{equation}
\label{dfjkvar} 
\Var_\lambda [ \sigma_k (kr-\tau)]
  =
\frac{1}{(k-1)!}
\sum_{l=0}^{k-2}
    {k-1 \choose l}
    (\lambda \tau )^{k-1+l}
    \frac{\Gamma ( (k-1-l)/2 +1)}{\Gamma ((k-1+l)/2+1)}, 
\end{equation} 
where $\Gamma$ denotes the gamma function defined as
$$
\Gamma (z) :=
\int_0^\infty x^{z-1} e^{-x} dx, \qquad z>0.
$$
\end{corollary} 
\begin{Proof} 
   Let $X_0,\ldots , X_k$ be independent standard normal random
 variables. From the moment relation
 $\E \big[ X_p^{2j_p} \big] = 2^{-j_p} (2j_p)!/j_p!$ and
 the fact that the sum $X_0^2 + \cdots + X_l^2$ has a Chi square
 distribution, we have
\begin{align} 
\nonumber 
  \sum_{j_0 + \cdots + j_l = k-1-l\atop j_0,\ldots , j_l\geq 0}
  \prod_{p=0}^l {2j_p\choose j_p} 
& =
\frac{2^{k-1-l}}{(k-1-l)!}
\E \left[
  \sum_{ j_0 + \cdots + j_l = k-1-l\atop j_0,\ldots , j_l\geq 0}
  \frac{(k-1-l)!}{j_0!\cdots j_l!}
    \prod_{p=0}^l X_p^{2j_p}
    \right]
\\
\nonumber
& 
=
\frac{2^{2(k-1-l)}}{(k-1-l)!}
\E \left[
 \left( \frac{X_0^2 + \cdots + X_l^2}{2} \right)^{k-1-l} 
\right]
\\
\nonumber
& 
= 2^{2(k-1-l)} \frac{\Gamma ( k-1-l + (l+1)/2)}{(k-1-l)!\Gamma ((l+1)/2)}. 
\end{align} 
 Hence from Proposition~\ref{djklfs2} we find 
\begin{eqnarray*} 
  \E  [ \sigma^2_k (kr-\tau)]
  & = & 
  \sum_{l=0}^{k-1}
  \frac{(\lambda \tau)^{2k-2-l}}{(2k-2-l)!} 
\frac{2^{2(k-1-l)}}{(k-1-l)!}
\frac{\Gamma ( k-1 + (1-l)/2)}{\Gamma ((l+1)/2)}
\\ & = & \sum_{l=0}^{k-1} \frac{(\lambda \tau)^{k-1+l}}{(k-1+l)!} \frac{2^{2l}}{l!} \frac{\Gamma ( (k-1+l+1)/2)}{\Gamma ((k-1-l+1)/2)}
\\
  & = & 
\frac{1}{(k-1)!}
\sum_{l=0}^{k-1}
{k-1 \choose l} (\lambda \tau)^{k-1+l} \frac{((k-1-l)/2)!}{((k-1+l)/2)!}. 
\end{eqnarray*} 
\end{Proof}
  The following table provides variance formulas of $k$-hop counts
  obtained from Corollary~\ref{djklfs2-2},
 by taking $\lambda =1$ for simplicity. 
\begin{table}[H]
\centering 
      \begin{tabular}{c|c|c|c|} 
          \hhline{~-}
        \multicolumn{1}{c|}{} & \cellcolor{gray!25}  
           \addstackgap[3pt]{Variance} 
\\ 
\cline{1-2} 
\multicolumn{1}{|c|}{\addstackgap[3pt]{$2$-hops}} & $\displaystyle \tau$
\\ 
\cline{1-2} 
\multicolumn{1}{|c|}{\addstackgap[10pt]{$3$-hops}} & \small $\displaystyle \frac{\tau^2}{2} + 2\frac{\tau^3}{3} 
$ 
\\ 
\cline{1-2} 
\multicolumn{1}{|c|}{\addstackgap[10pt]{$4$-hops}} & \small $\displaystyle \frac{\tau^3}{3!} + \frac{\tau^4}{4} + 2\frac{\tau^5}{15}$ 
\\ 
\cline{1-2} 
\multicolumn{1}{|c|}{\addstackgap[10pt]{$5$-hops}} & \small $\displaystyle \frac{\tau^4}{4!} + \frac{\tau^5}{15} + \frac{\tau^6}{24} + \frac{4\tau^7}{315}$
\\ 
\cline{1-2} 
\multicolumn{1}{|c|}{\addstackgap[10pt]{$6$-hops}} & \small $\displaystyle \frac{\tau^5}{5!} + \frac{\tau^6}{72} + \frac{\tau^7}{105} + \frac{\tau^8}{288} + \frac{2 \tau^9}{2835}$ 
\\ 
\cline{1-2}
\end{tabular} 
\caption{Variances of $k$-hop counts.} 
\end{table} 

 \vspace{-0.3cm}

\noindent
 Using the Legendre duplication formula
 $(2k-3)! \Gamma ( 3/2)
          = 
2^{2k-4}
(k-2)! \Gamma (k-1/2)$, 
 Corollary~\ref{djklfs2-2} also yields the
 following asymptotic variance. 
 In the sequel, for $f$ and $g$ two nonvanishing functions
 on $\real_+$ we write $f(\lambda ) \approx g(\lambda )$
 if $\lim_{\lambda \to \infty} f(\lambda ) / g(\lambda ) = 1$. 
 \begin{prop}
   As $\lambda$ tends to infinity we have the equivalence 
   \begin{equation}
     \label{fjkdls45}
    \Var_\lambda [ \sigma_k (kr-\tau)]
  \approx 
\frac{(2 \lambda \tau )^{2k-3}}{2(2k-3)!}. 
\end{equation} 
\end{prop} 
\section{Joint moments recursion} 
\label{sec5}
\noindent
  In the one-hop case 
we simply have
$\sigma_1 (t) = 1$ and $m^{(\lambda )}_{1,n} = 1$, $n \geq 0$.
 As for the two-hop count, \eqref{fjlksd23} 
   yields the joint Poisson moments formula 
\begin{equation} 
  \label{djfkls}
  m^{(\lambda )}_{2,n}(\tau_1,\ldots , \tau_n)
 =
 \E [ \sigma_2(2r-\tau_1)\cdots \sigma_2(2r-\tau_n) ]
 =
 \sum_{l=1}^n \lambda_1^l
     \sum_{\eta_1\cup \cdots \cup \eta_l = \{1,\ldots , n\}}
\prod_{j=1}^l \widehat{\tau}_{\eta_j}, 
\end{equation} 
 which follows from the expression 
 $\lambda_1 \min ( \tau_1, \ldots , \tau_n )$ 
 of the joint Poisson cumulants of
 $( \sigma_2(2r-\tau_1) , \ldots , \sigma_2(2r-\tau_n) )$.
 The direct application of Proposition~\ref{p1}
 to the evaluation of higher order joint moments 
 of $k$-hop counts 
 is not an easy task due to the complexity
 of the summations over partitions involved in \eqref{djkls12}.
 In Proposition~\ref{433} we propose to compute joint moments
 by a recursion argument 
 using the multiple Poisson stochastic integral
 representation \eqref{fjld3} %
 instead of the $U$-statistics expression \eqref{fjda}.
 Particular cases are considered with explicit computations for $n=1,2,3$ in
 Appendix~\ref{s9}.
\begin{prop}
  \label{433}
  For $k\geq 1$, the joint moments 
$$
 m^{(\lambda )}_{k,n} (\tau_1,\ldots , \tau_n)
 : =  \E_\lambda [ \sigma_k(kr-\tau_1)\cdots \sigma_k(kr-\tau_n) ],
 \qquad 0 \leq \tau_1,\ldots , \tau_n \leq r, 
$$
 satisfy the recursion 
\begin{eqnarray} 
  \label{fdfdf}
  \lefteqn{
 m^{(\lambda )}_{k+1,n} (\tau_1,\ldots , \tau_n)}
  \\
  \nonumber
  & = &  
  \sum_{l=1}^n
  \sum_{\pi_1\cup \cdots \cup \pi_l = \{1,\ldots , n\}}
  \int_0^{\widehat{\tau}_{\pi_l}}
  \cdots
  \int_0^{\widehat{\tau}_{\pi_1}}
  m^{(\lambda )}_{k,n} (\widebar{u}_{\pi_1} , \ldots , \widebar{u}_{\pi_l} ) 
  \lambda_k (du_1) \cdots \lambda_k (du_l), 
\end{eqnarray} 
 $0\leq \tau_1 , \ldots , \tau_n \leq r$,
where
$\widebar{u}_{\pi_i}:=( \underbrace{u_i,\ldots, u_i}_{|\pi_i|~ {\rm times }} )$ and
$\widehat{\tau}_{\pi}:= \min_{i\in \pi} \tau_i$ for $\pi \subset \{1,\ldots , n\}$. 
\end{prop}
\begin{Proof}
 By Proposition~\ref{fjkldsa}, when $\tau \in [0,r]$ we have 
 $\sigma^{(k+1)} (kr-\tau) \stackrel{d}{\simeq} Z^{(k+1)}_\tau$ in distribution,
 where $Z^{(k+1)}_\tau$ satisfies the recursion 
 \begin{equation}
   \label{kldsf} 
 Z^{(k+1)}_\tau  = \int_0^\tau Z^{(k)}_u dN^{(k)}_u, \qquad \tau \in [0,r], \quad k \geq 1, 
\end{equation} 
 and $(N^{(l)}_u)_{u\in \real_+}$ is a family of
 independent Poisson processes with respective intensities
 $\lambda_l(ds):=\lambda ( ds - (l-1) r)$,
 $l=1,\ldots , k$. 
 Hence, by Proposition~\ref{pr11-2}, 
 for $\tau_1 , \ldots , \tau_n \in [0,r]$ we have 
\begin{align} 
  \nonumber
  \lefteqn{
    m^{(\lambda )}_{k+1,n} (\tau_1,\ldots , \tau_n)  = 
  \E\big[ Z^{(k+1)}_{\tau_1} \cdots Z^{(k+1)}_{\tau_n} \big]
  }
  \\
\nonumber
    & = 
  \E\left[
    \int_0^{\tau_1} Z^{(k)}_u dN^{(k)}_u
    \cdots
    \int_0^{\tau_n} Z^{(k)}_u dN^{(k)}_u
\right] 
  \\
\nonumber
    & = 
  \sum_{l=1}^n
  \lambda_k^l
  \sum_{\pi_1\cup \cdots \cup \pi_l = \{1,\ldots , n\}}
  \int_{\real_+^l} 
  \E \Bigg[
    \prod_{j=1}^l
    \prod_{i\in \pi_j}
    \big( Z^{(k)}_{u_j} {\bf 1}_{[0,\tau_i]}(u_j) \big)
  \Bigg] 
  \lambda_k(du_1)\cdots \lambda_k(du_l)
  \\
\nonumber
      & = 
  \sum_{l=1}^n
  \lambda_k^l
  \sum_{\pi_1\cup \cdots \cup \pi_l = \{1,\ldots , n\}}
  \int_{\real_+^l} 
  \E \Bigg[ \prod_{j=1}^l
    \big(
    \big(Z^{(k)}_{u_j}\big)^{|\pi_j|}
    {\bf 1}_{[0,\widehat{\tau}_{\pi_j}]}(u_j)
    \big)
    \Bigg] 
  \lambda_k(du_1)\cdots \lambda_k(du_l)
  \\
\nonumber
      & = 
  \sum_{l=1}^n
  \lambda_k^l
  \sum_{\pi_1\cup \cdots \cup \pi_l = \{1,\ldots , n\}}
  \int_0^{\widehat{\tau}_{\pi_l}}
  \cdots
  \int_0^{\widehat{\tau}_{\pi_1}}
    \E \Bigg[ \prod_{j=1}^l
    \big(Z^{(k)}_{u_j}\big)^{|\pi_j|}
    \Bigg] 
  \lambda_k(du_1)\cdots \lambda_k(du_l)
  \\
  \label{kldsf-2}
  & = 
  \sum_{l=1}^n
  \lambda_k^l
  \sum_{\pi_1\cup \cdots \cup \pi_l = \{1,\ldots , n\}}
  \int_0^{\widehat{\tau}_{\pi_l}}
  \cdots
  \int_0^{\widehat{\tau}_{\pi_1}}
  m^{(\lambda )}_{k,n} (\widebar{u}_{\pi_1},\ldots , \widebar{u}_{\pi_l} ) 
  \lambda_k(du_1)\cdots \lambda_k(du_l). 
\end{align} 
\end{Proof}
 Table~\ref{table3} lists the first four joint moments
$m^{(\lambda )}_{2,n}(\tau_1,\ldots ,\tau_n )$ of
the two-hop counts, computed 
as an application of Proposition~\ref{433} 
from the command {\rm mk}[$\{ \tau_1 , \ldots , \tau_n \}, \{\lambda_1 \}$] in the Mathematica codes~\ref{code1}-\ref{code2} in appendix for $n=1,2,3,4$, 
when $\lambda_1(ds)=\lambda_1 ds =ds$. 
The case $\lambda_1 (ds)=\lambda_1 ds$, $\lambda_1>0$, is obtained by replacing
$\tau_i$ with $\lambda_1 \tau_i$, $i=1,2,3$. 
 
\begin{table}[H]
\centering 
      \begin{tabular}{c|l|} 
          \hhline{~-}
        \multicolumn{1}{c|}{} & \cellcolor{gray!25}  
 \addstackgap[3pt]{~~~~~~~~~~~~~~~~~~~~Joint moments of 2-hop counts}  
\\
\cline{1-2} 
\multicolumn{1}{|c|}{First} & $\tau_1$ 
\\ 
\cline{1-2} 
\multicolumn{1}{|c|}{Second} & $\displaystyle   \tau_1 + \tau_1 \tau_2
$ 
\\ 
\cline{1-2} 
\multicolumn{1}{|c|}{Third} & $\displaystyle
\tau_1 + \tau_1 \tau_3 + 2 \tau_1 \tau_2 + \tau_1\tau_2\tau_3 $  
\\ 
\cline{1-2}
\multicolumn{1}{|c|}{Fourth} & $\displaystyle
\tau_1 + \tau_1 \tau_4 + 2 \tau_1 \tau_3  + 4 \tau_1 \tau_2    
+ \tau_1\tau_3\tau_4 + 2 \tau_1\tau_2\tau_4   + 3 \tau_1\tau_2\tau_3
+ \tau_1\tau_2\tau_3\tau_4
$
  \\ 
\cline{1-2}
\end{tabular} 
\caption{Joint moments $m^{(\lambda )}_{2,n}(\tau_1,\ldots ,\tau_n )$ of $2$-hop counts of orders $n=1,2,3,4$.} 
\label{table3}
\end{table} 

\vskip-0.3cm

\noindent 
 Tables~\ref{table4} and \ref{table5} provide the first four moments of
 the three-hop and four-hop counts
 computed 
 from the commands {\rm mk}[$\{ \tau_1 , \ldots , \tau_n \}, \{\lambda_1,\lambda_2 \}$]
 and {\rm mk}[$\{ \tau_1 , \ldots , \tau_n \}, \{\lambda_1,\lambda_2 ,\lambda_3\}$]  in Mathematica for $n=1,2,3,4$, 
 where for simplicity we take $\lambda_1=\lambda_2 = \lambda_3 = 1$ and 
 $\tau_1=\tau_2=\tau_3=\tau$. 
 
\begin{table}[H]
\centering 
      \begin{tabular}{c|l|} 
          \hhline{~-}
        \multicolumn{1}{c|}{} & \cellcolor{gray!25}  
 \addstackgap[3pt]{~~~~~~~~~~Moments of 3-hop counts}  
\\ [0.3ex]  
\cline{1-2} 
\multicolumn{1}{|c|}{\addstackgap[10pt]{First}} & \small $\displaystyle \frac{\tau^2}{2} $ 
\\ [1ex]
\cline{1-2} 
\multicolumn{1}{|c|}{\addstackgap[10pt]{Second}} & \small $\displaystyle \frac{\tau^2}{2} + 2 \frac{\tau^3}{3} + \frac{\tau^4}{4} 
$ 
\\ [1ex]  
\cline{1-2} 
\multicolumn{1}{|c|}{\addstackgap[10pt]{Third}} & \small $\displaystyle \frac{\tau^2}{2} + 2\tau^3 + 5\frac{\tau^4}{2} + \tau^5 + \frac{\tau^6}{8} $  
\\ [1ex]   
\cline{1-2}
\multicolumn{1}{|c|}{\addstackgap[10pt]{Fourth}} & \small $\displaystyle
\frac{\tau^2}{2} + 14 \frac{\tau^3}{3} + 53 \frac{\tau^4}{4} + 66 \frac{\tau^5}{5} + 67 \frac{\tau^6}{12} + \tau^7 + \frac{\tau^8}{16}$ 
\\ [1ex]   
\cline{1-2}
\end{tabular} 
\caption{Moments $m^{(\lambda )}_{3,n}(\tau,\ldots ,\tau )$ of $3$-hop counts of orders $n=1,2,3,4$.} 
\label{table4}
\end{table} 

\vskip-0.3cm

\begin{table}[H]
\centering 
         \begin{tabular}{c|l|} 
          \hhline{~-}
        \multicolumn{1}{c|}{} & \cellcolor{gray!25}  
 \addstackgap[3pt]{~~~~~~~~~~Moments of 4-hop counts}  
\\ [0.3ex]  
\cline{1-2} 
\multicolumn{1}{|c|}{\addstackgap[10pt]{First}} & \small $\displaystyle \frac{\tau^3}{3!} $ 
\\ 
\cline{1-2} 
\multicolumn{1}{|c|}{\addstackgap[10pt]{Second}} & \small $\displaystyle \frac{\tau^3}{3!} + \frac{\tau^4}{4} + 2 \frac{\tau^5}{15} + \frac{\tau^6}{36}$
\\ 
\cline{1-2} 
\multicolumn{1}{|c|}{\addstackgap[10pt]{Third}} & \small $\displaystyle \frac{\tau^3}{3!} + 3 \frac{\tau^4}{4} + 5 \frac{\tau^5}{4} + \frac{59 \tau^6}{60} + \frac{13 \tau^7}{35} + \frac{\tau^8}{15} + \frac{\tau^9}{216}$
\\ 
\cline{1-2}
\end{tabular} 
\caption{Moments $m^{(\lambda )}_{4,n}(\tau,\ldots ,\tau )$ of $4$-hop counts of orders $n=1,2,3$.} 
\label{table5}
\end{table} 

\vspace{-0.3cm}

\noindent
The following figures plot higher joint moment formulas up to the order six, 
together with their confirmations by Monte Carlo simulations.
 
\begin{figure}[H]
  \centering
 \begin{subfigure}[b]{0.49\textwidth}
    \includegraphics[width=1\linewidth, height=5cm]{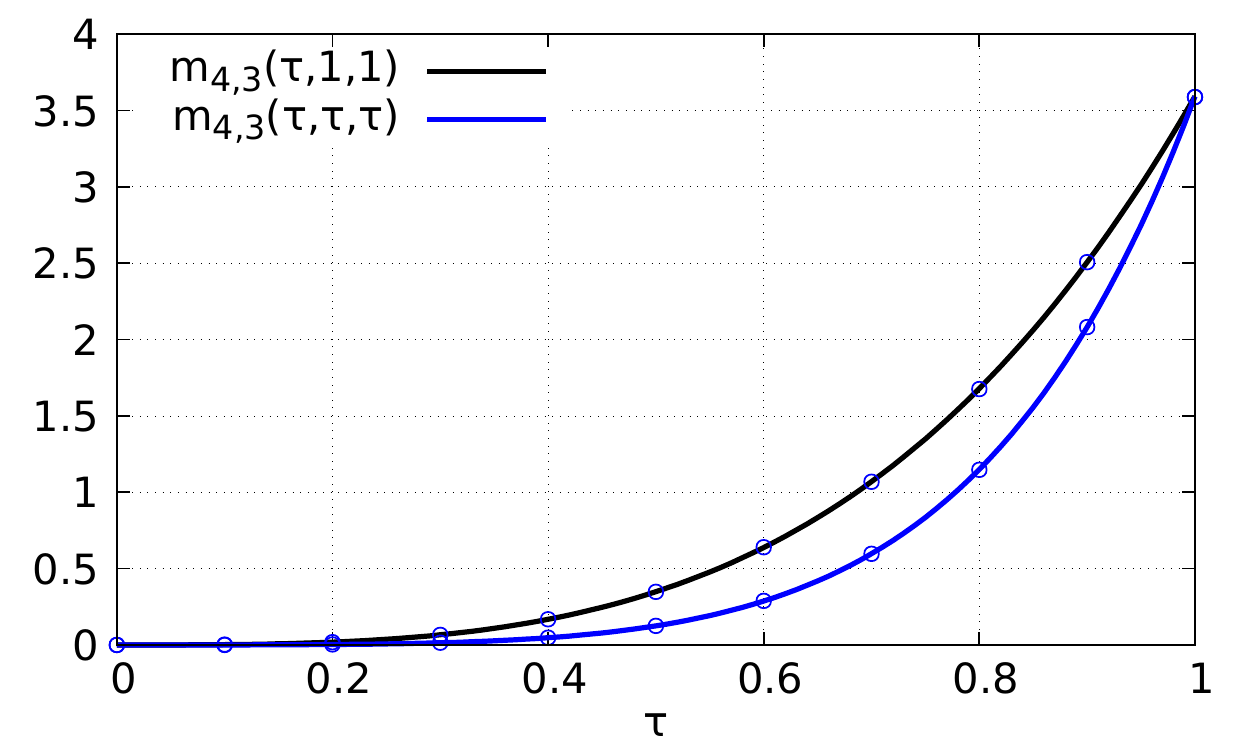}
    \caption{Third moments.} 
 \end{subfigure}
  \begin{subfigure}[b]{0.49\textwidth}
    \includegraphics[width=1\linewidth, height=5cm]{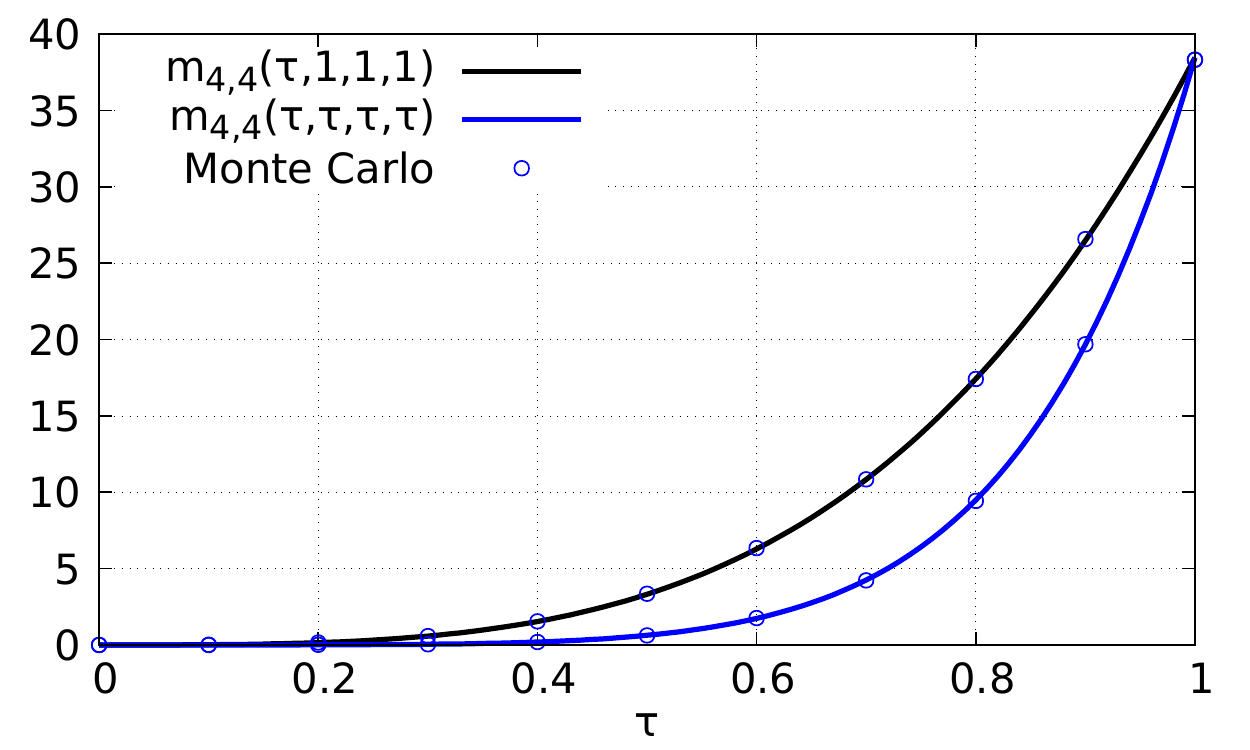} 
    \caption{Fourth moments.} 
  \end{subfigure}
    \caption{Third and fourth joint moments of $4$-hop counts.} 
\end{figure}

\begin{figure}[H]
  \centering
 \begin{subfigure}[b]{0.49\textwidth}
    \includegraphics[width=1\linewidth, height=5cm]{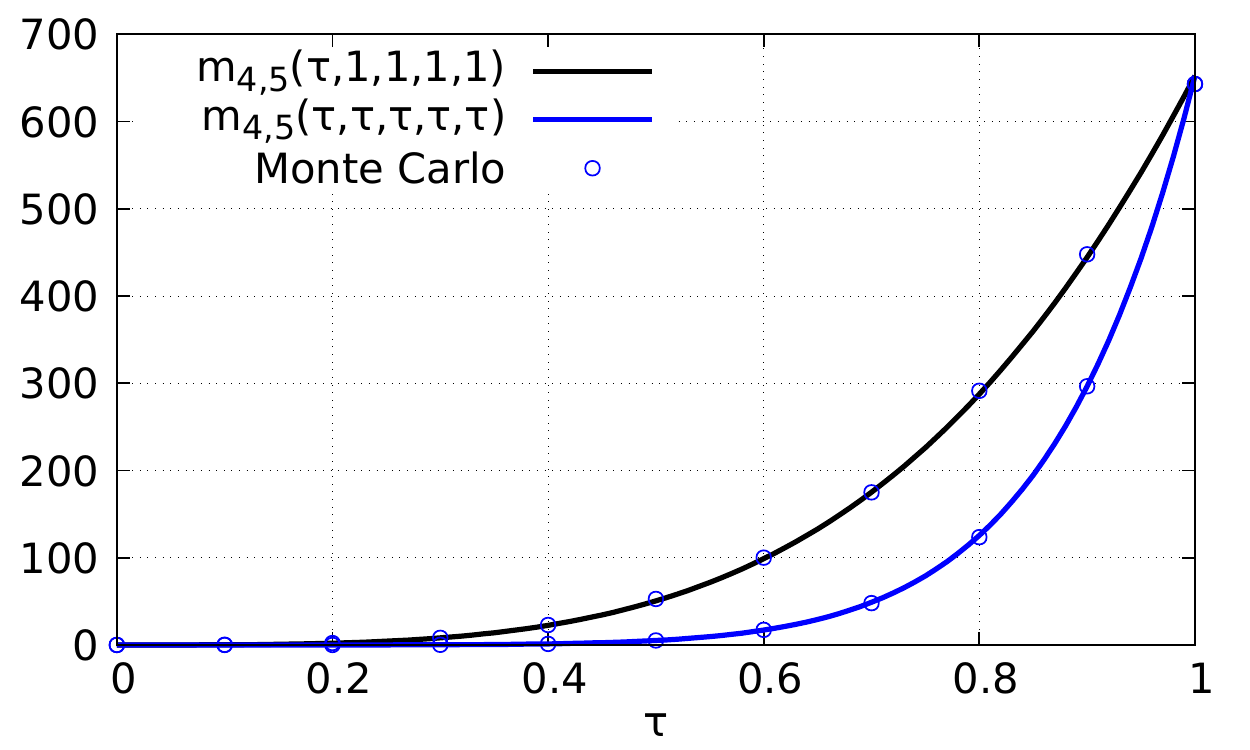}
    \caption{Fifth moments.} 
 \end{subfigure}
  \begin{subfigure}[b]{0.49\textwidth}
   \includegraphics[width=1\linewidth, height=5cm]{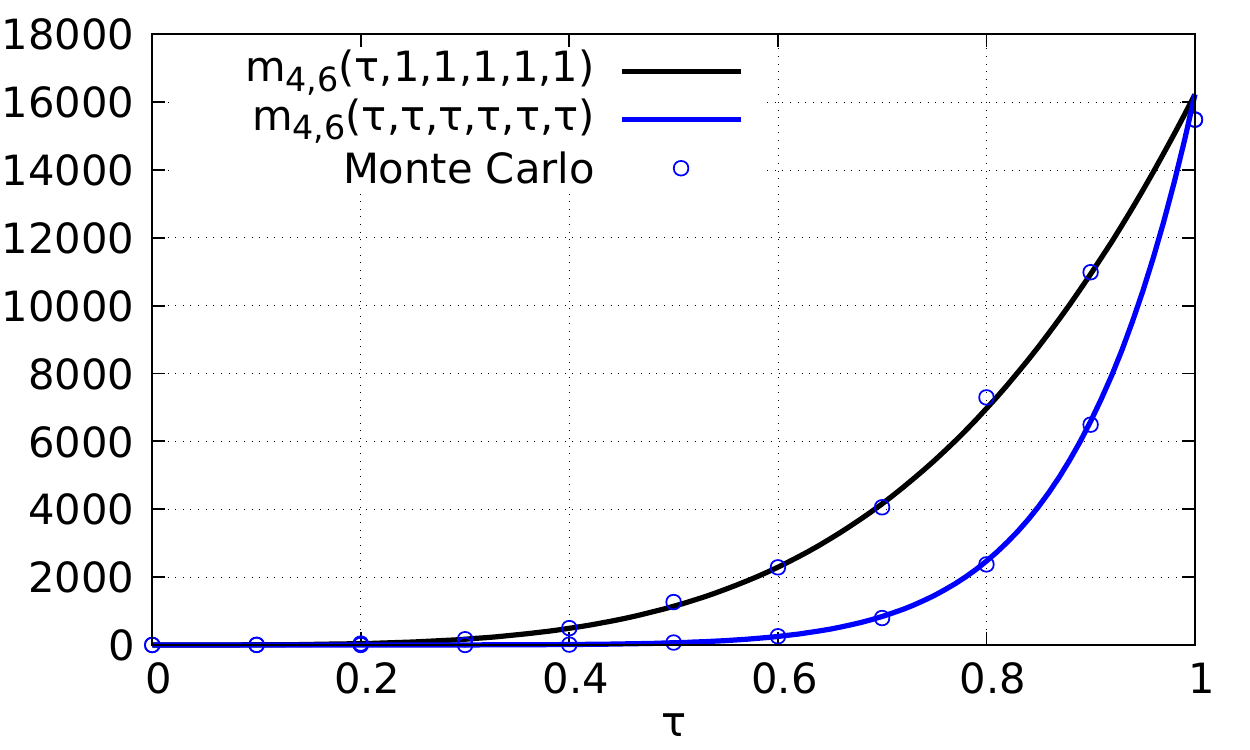} 
    \caption{Sixth moments.} 
  \end{subfigure}
    \caption{Fifth and sixth joint moments of $4$-hop counts.} 
\end{figure}
\section{Joint cumulants recursion} 
\label{sec6}
 In the sequel, given
 $\{ \pi_1,\ldots , \pi_l \}$ a partition of $\{1,\ldots , n\}$
 and $\tau_1,\ldots , \tau_l \in [0,r]$ we denote by
 $$
 c_{k,n}^{(\lambda)} ( \widebar{\tau}_{\pi_1} ; \ldots ; \widebar{\tau}_{\pi_l} ) : = 
 \kappa_\lambda \big( ( \sigma_k(kr - \tau_1) )^{|\pi_1|}, \ldots , ( \sigma_k(kr - \tau_l) )^{|\pi_l|}\big)
 $$ 
 the joint cumulant
 of $\big( ( \sigma_k(kr - \tau_1) )^{|\pi_1|}, \ldots , ( \sigma_k(kr - \tau_l) )^{|\pi_l|}\big)$, and we let 
$c_{k,n}^{(\lambda)} ( \tau_1 ; \ldots ; \tau_n )$ 
denote the joint cumulant  
$\kappa_\lambda ( \sigma_k(kr - \tau_1) , \ldots , \sigma_k(kr - \tau_n ) )$. 
 The following proposition is the counterpart of the moment recursion
 of Proposition~\ref{433} obtained by the M\"obius inversion relation 
 \eqref{bpi}. 
  Particular cases are considered with explicit computations for $n=2,3,4$ in
  Appendix~\ref{s10}.
\begin{prop}
\label{fjhkds}
 Let $\tau \in [0,r]$ and $k, n\geq 1$. 
 The cumulant of order $n$ of $\sigma_{k+1} (kr - \tau)$ 
 satisfies the recursion 
\begin{equation} 
  \label{fjkla}
  c_{k+1,n}^{(\lambda)} (\tau_1 ; \ldots ; \tau_n )
= 
  \sum_{l=1}^n
  \sum_{\pi_1\cup \cdots \cup \pi_l = \{1,\ldots , n\} }
  \int_0^{\widehat{\tau}_{\pi_1}} \cdots \int_0^{\widehat{\tau}_{\pi_l}} 
  c_{k,l}^{(\lambda)} (\widebar{s}_{\pi_1} ; \ldots ; \widebar{s}_{\pi_l} ) 
  \lambda_k ( ds_1 ) \cdots \lambda_k ( ds_l ) 
  , 
\end{equation} 
  where we let
  $\widehat{\tau}_{\pi}:= \min_{i\in \pi} \tau_i$ for $\pi \subset \{1,\ldots , n\}$. 
\end{prop}
\begin{Proof}
 When $k=1$ we have $\sigma_1(\tau ) =1$, $\tau \in [0,r]$,
hence $c_{1,n}^{(\lambda)} (\tau_1;\ldots ;  \tau_n ) = {\bf 1}_{\{ n =1\}}$, and 
 by \eqref{fjlksd23} the cumulants of the two-hop count 
 $\sigma_2 (2r-\tau) = N_r - N_{r-\tau} \stackrel{d}{\simeq} N^{(1)}_\tau$, 
 $\tau \in [0,r]$, 
 are the joint Poisson cumulants 
 \begin{equation}
   \label{fjkldsf-c} 
 c_{2,n}^{(\lambda)} (\tau_1 ; \ldots ; \tau_n ) = \min (\tau_1,\ldots ,  \tau_n ),
 \qquad
 n \geq 1, 
\end{equation} 
 which is consistent with the joint Poisson moments formula 
 \eqref{djfkls} and shows that \eqref{fjkla} holds at the rank $k=1$.
 Next,
 assuming that \eqref{fjkla} holds at the rank $k\geq 1$, 
 by the joint cumulant-moment inversion relation \eqref{jklsda1}
 we have 
\begin{align*} 
 & 
 c_{k+1,n}^{(\lambda)} (\tau_1 ; \ldots ; \tau_n ) 
 = 
 \sum_{l=1}^n
  (l-1)!
  (-1)^{l-1}
  \sum_{\pi_1\cup \cdots \cup \pi_l = \{1,\ldots , n\}}
  \prod_{q=1}^l 
  m^{(\lambda )}_{k+1,|\pi_q|} ( \tau_{\pi_q} ) 
  \\
 & = 
  \sum_{l=1}^n
  (l-1)!
  (-1)^{l-1}
  \\
   &  \quad \sum_{\pi_1\cup \cdots \cup \pi_l = \{1,\ldots , n\}}
  \prod_{q=1}^l 
  \sum_{\eta^q \preceq \pi_q } 
  \int_0^{\widehat{\tau}_{\eta^q_1}} 
  \cdots
  \int_0^{\widehat{\tau}_{\eta^q_{|\eta^q|}} }
  m^{(\lambda )}_{k,|\pi_q|} \big(\widebar{u}_{\eta^q_1} ;\ldots ; \widebar{u}_{\eta^q_{|\eta^q|}} \big) 
  \lambda_k (du_1) \cdots \lambda_k (du_{|\eta^q|}), 
\\
 & =  
  \sum_{l=1}^n
  (l-1)!
  (-1)^{l-1}
  \\
   & \quad \sum_{\pi_1\cup \cdots \cup \pi_l = \{1,\ldots , n\}}
  \sum_{\eta^i \preceq \pi_i 
    \atop { 1 \leq i \leq l }
    }
  \prod_{q=1}^l
  \int_0^{\widehat{\tau}_{\eta^q_1}} 
  \cdots
  \int_0^{\widehat{\tau}_{\eta^q_{|\eta^q|}} }
  m^{(\lambda )}_{k,|\eta^q|} \big(\widebar{u}_{\eta^q_1} ;\ldots ; \widebar{u}_{\eta^q_{|\eta^q|}} \big) 
  \lambda_k (du_1) \cdots \lambda_k (du_{|\eta^q|}), 
\\
 & =  
  \sum_{l=1}^n
  (l-1)!
  (-1)^{l-1}
  \\
  & \quad
  \sum_{\pi_1\cup \cdots \cup \pi_l = \{1,\ldots , n\}}
  \sum_{\eta^i \preceq \pi_i 
  \atop { 1 \leq i \leq l } }
  \prod_{q=1}^l
  \sum_{\psi^q \preceq \eta^q }
  \int_0^{\widehat{\tau}_{\eta^q_1}} 
  \cdots
  \int_0^{\widehat{\tau}_{\eta^q_{|\eta^q|}} }
  \prod_{m=1}^{|\psi^q|}
  c_{k,n}^{(\lambda)} \big(\big( \widebar{u}_{\eta^q_i}\big)_{i \in \psi^q_m} \big) 
  \lambda_k (du_1) \cdots \lambda_k (du_{|\eta^q|}) 
  \\
 & =  
  \sum_{l=1}^n
  (l-1)!
  (-1)^{l-1}
  \\
  & 
  \quad 
  \sum_{\pi_1\cup \cdots \cup \pi_l = \{1,\ldots , n\}}
  \sum_{\psi^i \preceq \eta^i \preceq \pi_i 
  \atop 
      { 1 \leq i \leq l }
      }
  \int_0^{\widehat{\tau}_{\eta^q_1}} 
  \cdots
  \int_0^{\widehat{\tau}_{\eta^q_{|\eta^q|}} }
  \prod_{m=1}^{|\psi^q|}
  \prod_{q=1}^l
  c_{k,n}^{(\lambda)} \big(\big( \widebar{u}_{\eta^q_i}\big)_{i \in \psi^q_m} \big) 
  \lambda_k (du_1) \cdots \lambda_k (du_{|\psi^q|})
  \\
 & =  
  \sum_{l=1}^n
  \mu ( \sigma , \widehat{\bf 1} ) 
  \sum_{\eta \preceq \sigma \preceq \widehat{\bf 1} } 
  \sum_{\psi \preceq \eta} 
  \int_0^{\widehat{\tau}_{\psi_1}} 
  \cdots
  \int_0^{\widehat{\tau}_{\psi_{|\psi |}} }
  \prod_{A\in \eta}
  c_{k,n}^{(\lambda)} \big(\big( \widebar{u}_{A \cap B} \big)_{B \in \psi} \big) 
  \lambda_k (du_1) \cdots \lambda_k (du_{|\psi |})
  \\
 & =  
  \sum_{l=1}^n
  \mu ( \sigma , \widehat{\bf 1} ) 
  \sum_{\eta \preceq \sigma \preceq \widehat{\bf 1} } 
  G(\eta)
  \\
    & =  G ( \widehat{\bf 1})
    \\
 & = 
  \sum_{\eta \preceq \widehat{\bf 1} } 
  \int_0^{\widehat{\tau}_{\eta_1}} 
  \cdots
  \int_0^{\widehat{\tau}_{\eta_{|\eta|}} }
    c_{k,n}^{(\lambda)} (\widebar{u}_{\eta_1};\ldots ; \widebar{u}_{\eta_{|\eta|}} ) 
    \lambda_k (du_1) \cdots \lambda_k (du_{|\eta|})
, 
\end{align*} 
 where we applied Relation~\eqref{fjlksa1} 
 with $\pi := \widehat{\bf 1}$ and
 $$
  G ( \eta ) := 
   \sum_{\psi \preceq \eta} 
  \int_0^{\widehat{\tau}_{\psi_1}} 
  \cdots
  \int_0^{\widehat{\tau}_{\psi_{|\psi |}} }
  \prod_{A\in \eta}
  c_{k,n}^{(\lambda)} \big(\big( \widebar{u}_{A \cap B} \big)_{B \in \psi} \big) 
  \lambda_k (du_1) \cdots \lambda_k (du_{|\psi |})
$$
 which shows \eqref{fjkla} by induction on $k\geq 1$. 
\end{Proof}
 In order to use
 Proposition~\ref{fjhkds} as an induction relation, 
 the higher order cumulants $c_{k,n}^{(\lambda)} ( \widehat{s}_{\pi_1} ;\ldots ; \widehat{s}_{\pi_l} )$
 appearing in \eqref{fjkla} can be computed by recurrence using the
following proposition. 
\begin{prop}
  \label{jkld}
 For any sequence $(X_1,\ldots , X_{n+1})$ of random variables we have
 the cumulant relation
$$
  \kappa (X_1,\ldots , X_n X_{n+1}) 
 = 
 \kappa (X_1,\ldots ,X_n,X_{n+1} ) + \sum_{\eta_1\cup \eta_2 = \{ 1, \ldots , n+1 \} \atop \eta_1 \ni n, \ \! \eta_2 \ni n+1, \ \! \eta_1 \cap \eta_2 = \emptyset 
 }
\kappa (X_{\eta_1} , X_n)\kappa (X_{\eta_2} , X_{n+1}). 
$$ 
\end{prop} 
\begin{Proof}
 This relation is a particular case of the cumulant-moment relationship 
$$ 
 \kappa \big(Z_1,\ldots , Z_n\big) 
 = 
 \sum_{l=1}^n
  (l-1)!
  (-1)^{l-1}
  \sum_{\pi_1\cup \cdots \cup \pi_l = \{1,\ldots , n\}}
  \prod_{j=1}^l 
  \E \Bigg[ \prod_{i\in \pi_j} Z_i \Bigg]
$$
 which yields, taking $Z_1:=X_1,Z_2:=X_2,\ldots , Z_n:=X_nX_{n+1}$,  
\begin{eqnarray*} 
  \lefteqn{
   \! \! \! \!   \kappa (X_1,\ldots , X_n X_{n+1}) 
 = 
 \sum_{l=1}^n
  (l-1)!
  (-1)^{l-1}
  \sum_{\pi_1\cup \cdots \cup \pi_l = \{1,\ldots , n\}}
  \prod_{j=1}^l 
  \E \Bigg[ \prod_{i\in \pi_j} X_i \prod_{ \pi_j \ni n } ( X_nX_{n+1}) \Bigg]
  }
 \\
   & = & 
 \sum_{l=1}^n
  (l-1)!
  (-1)^{l-1}
  \sum_{\pi_1\cup \cdots \cup \pi_l = \{1,\ldots , n\}}
  \prod_{j=1}^l 
  \sum_{p=1}^{|\pi_j|+1} 
  \sum_{\eta_1\cup \cdots \cup \eta_p =
    \pi_j \cup \{ n, n+1\} }
 \prod_{v=1}^p  \kappa \big( ( X_u)_{u\in \eta_p} \big) 
 \\
  & = &  
\kappa (X_1,\ldots ,X_n,X_{n+1} ) + \sum_{\eta_1\cup \eta_2 = \{ 1, \ldots , n+1 \} \atop \eta_1 \ni n, \ \! \eta_2 \ni n+1, \ \! \eta_1 \cap \eta_2 = \emptyset  }
\kappa (X_{\eta_1} , X_n) \kappa (X_{\eta_2} , X_{n+1}), 
\end{eqnarray*} 
\end{Proof}
\noindent
As an application of Proposition~\ref{jkld},
 Tables~\ref{table6} and \ref{table7} present the first five and first three cumulants
$c_{3,n}^{(\lambda)} (\tau ; \ldots ; \tau )$ of
the three-hop count
and  $c_{4,n}^{(\lambda)} (\tau ; \ldots ; \tau )$  of the four-hop count  
computed by the commands
{\rm ck}[$\{ \tau_1 , \ldots , \tau_n \}, \{ 1 , \ldots , 1 \} , \{ \lambda_1 ,\lambda_2 ,\lambda_2 \} ]$
and
{\rm ck}[$\{ \tau_1 , \ldots , \tau_n \}, \\ \{ 1 , \ldots , 1 \} , \{ \lambda_1 , \lambda_2 , \lambda_3 \} ]$
in the Mathematica codes~\ref{code3}-\ref{code4} in appendix for $n=1,2,3,4,5$,
 where for simplicity we take $\lambda_i = 1$, $i=1,2,3$ and 
 $\tau_j=\tau$, $j=1,\ldots , 5$. 

 \begin{table}[H]
\centering 
      \begin{tabular}{c|l|} 
          \hhline{~-}
        \multicolumn{1}{c|}{} & \cellcolor{gray!25}  
 \addstackgap[3pt]{~~~Cumulants of 3-hop counts}  
\\ [0.3ex]  
\cline{1-2} 
\multicolumn{1}{|c|}{\addstackgap[10pt]{First}} & \small $\displaystyle \frac{\tau^2}{2}$ 
\\ 
\cline{1-2} 
\multicolumn{1}{|c|}{\addstackgap[10pt]{Second}} & \small $\displaystyle \frac{\tau^2}{2}+\frac{2\tau^3}{3}$ 
\\ 
\cline{1-2} 
\multicolumn{1}{|c|}{\addstackgap[10pt]{Third}} & \small $\displaystyle \frac{\tau^2}{2} + 2\tau^3 + \frac{7\tau^4}{4}$  
\\ 
\cline{1-2}
\multicolumn{1}{|c|}{\addstackgap[10pt]{Fourth}} & \small $\displaystyle
\frac{\tau^2}{2}
+\frac{14 \tau^3}{3}
+\frac{23 \tau^4}{2}
+ \frac{36 \tau^5}{5}
$
\\ 
\cline{1-2}
\multicolumn{1}{|c|}{\addstackgap[10pt]{Fifth}} & \small $\displaystyle
\frac{\tau^2}{2} + 10 \tau^3 + \frac{215}{4} \tau^4 + 86 \tau^5 + 41 \tau^6
$
\\ 
\cline{1-2}
\end{tabular} 
\caption{Cumulants $c_{3,n}^{(\lambda)} (\tau ; \ldots ; \tau )$ of $3$-hop counts of orders $n=1,2,3$.} 
\label{table6}
\end{table} 
 
\vskip-0.3cm

\begin{table}[H]
\centering 
      \begin{tabular}{c|l|} 
          \hhline{~-}
        \multicolumn{1}{c|}{} & \cellcolor{gray!25}  
                           \addstackgap[3pt]{~Cumulants of $4$-hop counts}  
\\ [0.3ex]  
\cline{1-2} 
\multicolumn{1}{|c|}{\addstackgap[10pt]{First}} & \small $\displaystyle \frac{\tau^3}{3!}$ 
\\ 
\cline{1-2} 
\multicolumn{1}{|c|}{\addstackgap[10pt]{Second}} & \small $\displaystyle \frac{\tau^3}{3!} + \frac{\tau^4}{4} + \frac{2 \tau^5}{15}$ 
\\ 
\cline{1-2} 
\multicolumn{1}{|c|}{\addstackgap[10pt]{Third}} & \small $\displaystyle \frac{\tau^3}{3!} + \frac{3 \tau^4}{4} + \frac{5 \tau^5}{4} + \frac{9 \tau^6}{10} + \frac{69 \tau^7}{280}$ 
\\ 
\cline{1-2}
\end{tabular} 
\caption{Cumulants $c_{4,n}^{(\lambda)} (\tau ; \ldots ; \tau )$ of $4$-hop counts of orders $n=1,2,3$.} 
\label{table7}
\end{table}

\vskip-0.3cm 

\noindent
The following figures present third and fourth order cumulant plots 
for $4$-hop counts,
 together with their confirmations by Monte Carlo simulations.
 
\noindent
\begin{figure}[H]
  \centering
 \begin{subfigure}[b]{0.49\textwidth}
    \includegraphics[width=1\linewidth, height=5cm]{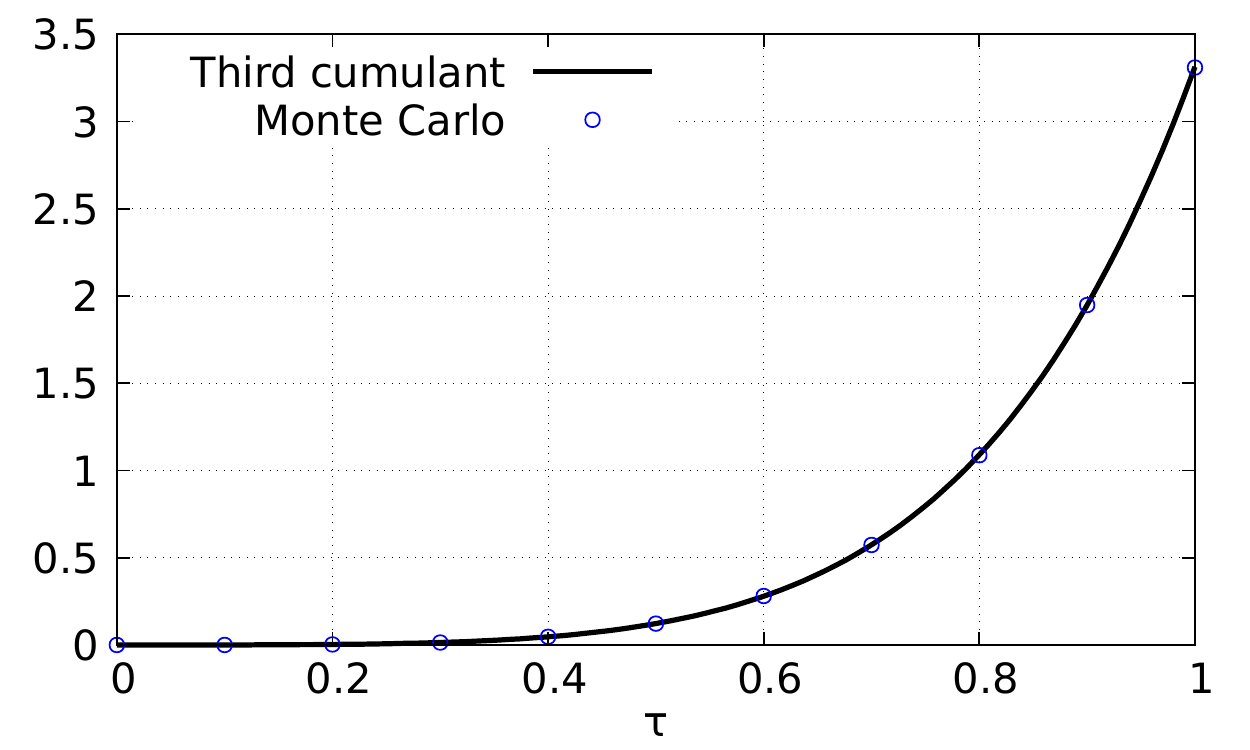}
    \caption{Third cumulant
      $c_{4,3}^{(\lambda)} (\tau ; \tau ; \tau )$.} 
 \end{subfigure}
  \begin{subfigure}[b]{0.49\textwidth}
   \includegraphics[width=1\linewidth, height=5cm]{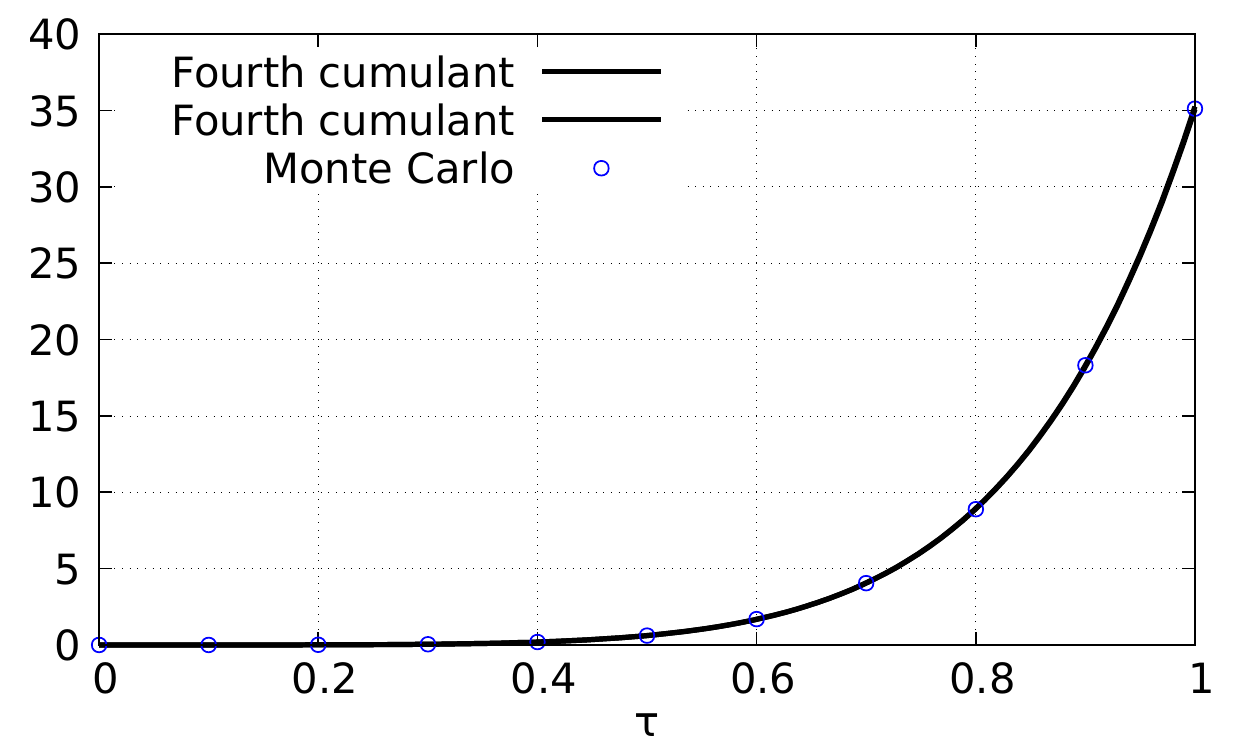} 
   \caption{Fourth cumulant
     $c_{4,4}^{(\lambda)} (\tau ;\tau ; \tau ; \tau )$.} 
  \end{subfigure}
    \caption{Third and fourth cumulants of the $4$-hop count 
 $\sigma_4(4-\tau )$.} 
\end{figure}
 
\begin{figure}[H]
  \centering
 \begin{subfigure}[b]{0.49\textwidth}
    \includegraphics[width=1\linewidth, height=5cm]{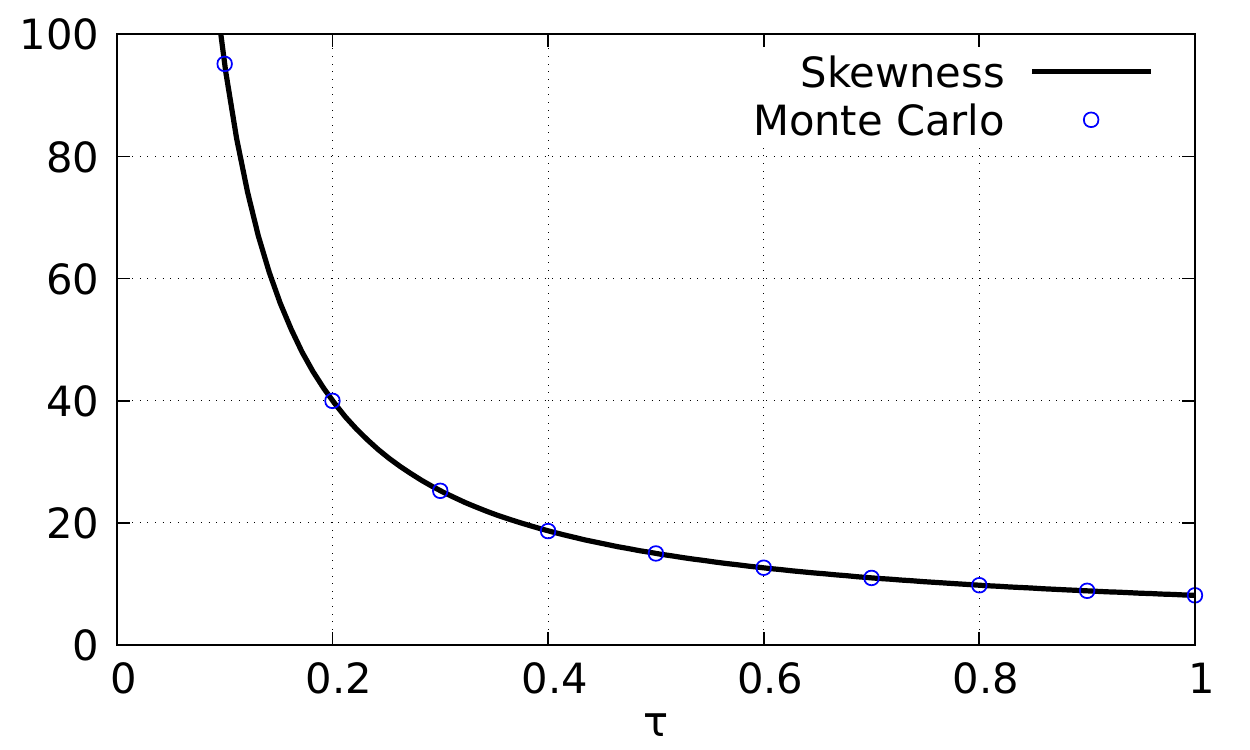}
    \caption{Skewness of $\sigma_4(4-\tau )$.} 
 \end{subfigure}
  \begin{subfigure}[b]{0.49\textwidth}
    \includegraphics[width=1\linewidth, height=5cm]{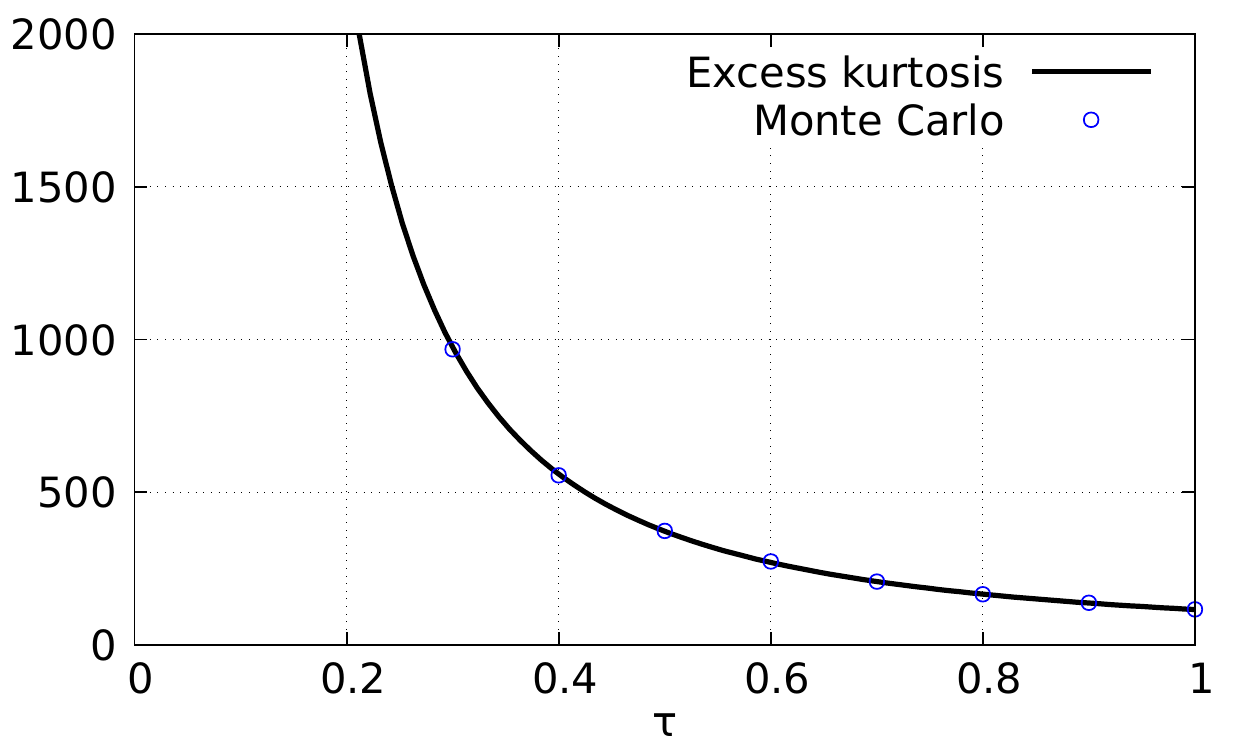} 
    \caption{Kurtosis of $\sigma_4(4-\tau )$.} 
  \end{subfigure}
   \caption{Skewness and kurtosis of
    the $4$-hop count $\sigma_4(4-\tau )$.} 
\end{figure}

\section{Moment and cumulant bounds}
\label{sec7}
In this section we take 
 $\lambda_l(ds):= \lambda ds$, $l=1,\ldots , k$,
 $\lambda >0$,
 and we write $f(\lambda ) = O(\lambda^n )$ if there exist $C_n>0$ and
$\lambda_n >0$ such that
$|f(\lambda )| \leq C_n \lambda^n$ for any $\lambda >\lambda_n$.
\begin{prop}
 Moment bound. For any $\tau \in [0,r]$ and $n\geq 0$, we have
$$ 
 \E_\lambda [ ( \sigma_k (kr - \tau ) )^n ]
  \leq 
  ( \E [ ( N_{\lambda \tau} )^n ] )^{k-1} 
 =
 O \big( ( \lambda \tau  )^{(k-1)n} \big)
 $$
 as $\lambda$ tends to infinity.
\end{prop}
\begin{Proof}
 We show by induction on $k\geq 1$ that 
$$ 
  m^{(\lambda )}_{k,n} (\tau_1 ,\ldots , \tau_n )
  \leq 
  ( \E [ ( N_{\lambda \tau } )^n ] )^{k-1}
  , 
  \qquad
  0 \leq \tau_1,\ldots , \tau_n \leq \tau. 
$$ 
  where,
  denoting by
  $S(n,l)$ the Stirling number of the second
 kind,  
  $$
  \E [ ( N_{\lambda \tau } )^n ] = 
\sum_{l=1}^n
S(n,l) ( \lambda \tau )^l
= O(\lambda^n).
$$
 The case $k=1$ is covered by the fact that 
  $\sigma_1 (r - \tau ) = 1$ and $m^{(\lambda )}_{1,n}(\tau ) = 1$, $n \geq 0$.
  Next, by the recurrence relation \eqref{fdfdf} we have
\begin{eqnarray*} 
  m^{(\lambda )}_{k+1,n} (\tau ,\ldots , \tau )
   & \leq &  
  \sum_{l=1}^n
  ( \lambda \tau )^l \sum_{\pi_1\cup \cdots \cup \pi_l = \{1,\ldots , n\}}
  \sup_{0 \leq \tau_1,\ldots , \tau_l \leq \tau }
  m^{(\lambda )}_{k,n} (\tau_1,\ldots , \tau_l ) 
  \\
   & = &  
  \sup_{0 \leq \tau_1,\ldots , \tau_l \leq \tau}
  m^{(\lambda )}_{k,n} (\tau_1,\ldots , \tau_l ) 
  \sum_{l=1}^n
  S(n,l) ( \lambda \tau )^l 
  \\
   & = &  
  \E [ ( N_{\lambda \tau } )^n ] 
  m^{(\lambda )}_{k,n} (\tau,\ldots , \tau ),
  \qquad k \geq 1. 
\end{eqnarray*} 
\end{Proof} 
The following bound on joint cumulants is also obtained by induction. 
 \begin{prop}
   \label{1fdjkl} 
   For any $\tau \in [0,r]$,
   $k\geq 2$ and $l_1,\ldots , l_p \geq 0$, $p\geq 1$,  
    we have the joint cumulant bound
 \begin{equation}
\label{c3789} 
   \kappa_\lambda \big( ( \sigma_k(kr - \tau ) )^{l_1}, \ldots , ( \sigma_k(kr - \tau ) )^{l_p}\big)
    \leq (2(\lambda + 1) \tau )^{(k-1)(l_1+\cdots + l_p)+1-p} (B_n)^{k-2}, 
\end{equation} 
 where $B_n$ denotes the Bell number of order $n\geq 1$. 
 In particular, we have
\begin{equation}
\label{fjkls343} 
 c_{k,n}^{(\lambda)} (\tau ; \ldots ; \tau ) \leq (2(\lambda + 1) \tau )^{ 1 + (k-2)n } (B_n)^{k-2},
 \quad \tau \in [0,r].
\end{equation} 
\end{prop}
\begin{Proof}
 We note that for $n=p$, as in \eqref{fjkldsf-c} we have 
 $c_{2,n}^{(\lambda)} (\tau_1 ; \ldots ; \tau_n) = \min (\tau_1,\ldots , \tau_n)$,
 hence by induction from Proposition~\ref{jkld} we obtain 
$$
 \kappa_\lambda \big( ( \sigma_2(2r - \tau_1) )^{l_1}, \ldots , ( \sigma_2(2r - \tau_p) )^{l_p}\big)
 \leq (2(\lambda + 1)\tau )^{l_1+\cdots + l_p +1-p}, 
$$ 
 $0 \leq \tau_1, \ldots , \tau_p \leq \tau$,
 which is \eqref{c3789} for $k=2$.
 Next, using Proposition~\ref{fjhkds},
 assuming that \eqref{c3789} holds at the rank $k \geq 2$,
 by induction for $\lambda \geq 1$
 we have 
\begin{eqnarray} 
   \lefteqn{
     \nonumber
     \! \! \! \! \! \! \! \! \! \! \! \! \! \! \! \! \!
     \big|
   \kappa_\lambda \big( \sigma_{k+1}((k+1)r - \tau ) , \ldots , \sigma_{k+1}((k+1)r - \tau ) \big)
    \big|  
   }
   \\
   & \leq & 
  \sum_{l=1}^n
  \lambda^l \sum_{\pi_1\cup \cdots \cup \pi_l = \{1,\ldots , n\} }
  \int_{[0,\tau ]^l} 
  |
  c_{k,l}^{(\lambda)} (\widebar{\tau}_{\pi_1} ;\ldots ; \widebar{\tau}_{\pi_l} )
  | 
  d\tau_1\cdots d\tau_l
    \\
  \nonumber 
  & = & 
  (B_n)^{k-2}
  \sum_{l=1}^n
  ( \lambda \tau ) ^l 
  \sum_{\pi_1\cup \cdots \cup \pi_l = \{1,\ldots , n\} }
  (2(\lambda +1) \tau )^{(k-1)n+1-l}
    \\
  \nonumber 
  & \leq & 
  (B_n)^{k-2}
  (2(\lambda +1) \tau )^{(k-1)n+1}
   \sum_{l=1}^n
  S(n,l)
    \\
  \nonumber 
  & = & 
  (B_n)^{k-1}
  (2(\lambda +1) \tau )^{(k-1)n+1}, 
\end{eqnarray} 
 which yields 
$$
 \kappa_\lambda \big( ( \sigma_{k+1}((k+1)r - \tau_1) )^{l_1}, \ldots , ( \sigma_{k+1}((k+1)r - \tau_p ) )^{l_p}\big)
 \leq (B_n)^{k-1}(2(\lambda + 1)\tau )^{k ( l_1 + \cdots + l_p ) +1-p}, 
$$ 
 $0 \leq \tau_1, \ldots , \tau_p \leq \tau$,
 by induction from Proposition~\ref{jkld}. 
\end{Proof}
\section{Berry-Esseen bounds} 
\label{sec8}
 In this section we will use the Wasserstein and Kolmogorov distances
 $d_W(X,Y)$ and $d_K(X,Y)$ 
 between the distributions of random variables $X, Y$, defined as
$$
 d_W (X,Y):
 =\sup_{h\in\mathrm{Lip}(1)} |\mathrm{E}[h(X)]-\mathrm{E}[h(Y)]|,
$$
 where $\mathrm{Lip}(1)$ denotes the class of real-valued
 Lipschitz functions with
 Lipschitz constant less than or equal to $1$,
 and
 $$
 d_K(X,Y) : = \sup_{x\in \real} | \P ( X \leq x ) - \P( Y \leq x ) |.
 $$ 
 For $\lambda>0$ and $t\in \real_+$ we let
 \begin{eqnarray}
   \label{fjkldsf1}
   \lefteqn{ 
 \sigma_k^{(\lambda )} (t)
  :=  \int_0^{\lambda t} 
  \cdots
  \int_0^{\lambda t}
  f_k(s_1/\lambda ,\ldots , s_{k-1} / \lambda ) dN_{s_1} \cdots dN_{s_{k-1}}
   }
   \\
   \nonumber
   & = & \frac{1}{(k-1)!} \sum_{l=0}^{k-1}
  \lambda^{k-1-l}
         {k-1 \choose l}
   I_l \left( \int_0^t \cdots \int_0^t 
   f_{k-1}(*,s_{l+1},\ldots , s_{k-1} ) ds_{l+1}\cdots ds_{k-1} \right) 
\end{eqnarray}  
 according to \eqref{fjkldsf} and \eqref{skt},
 so that the distribution of $\sigma_k^{(\lambda )} (t)$
 under $\P$ is the distribution of $\sigma_k (t)$ under
 the distribution $\P_\lambda$ of the 1D unit disk graph
 with constant Poisson intensity $\lambda >0$.

 \medskip
 
\begin{figure}[H]
  \centering
 \begin{subfigure}[b]{0.49\textwidth}
    \includegraphics[width=1\linewidth, height=5cm]{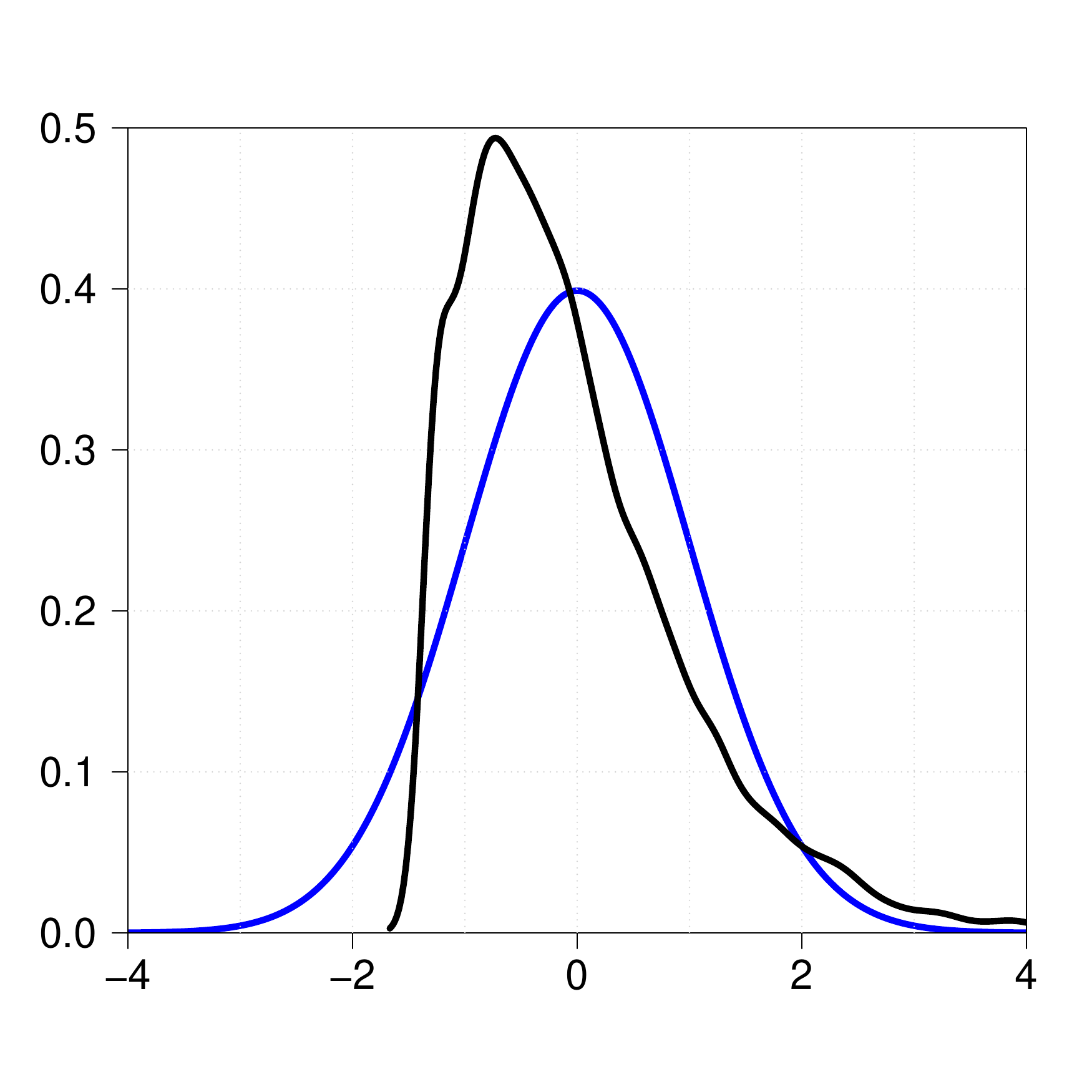}
       \vskip-.2cm
 \caption{$\lambda = 5$.} 
 \end{subfigure}
  \begin{subfigure}[b]{0.49\textwidth}
    \includegraphics[width=1\linewidth, height=5cm]{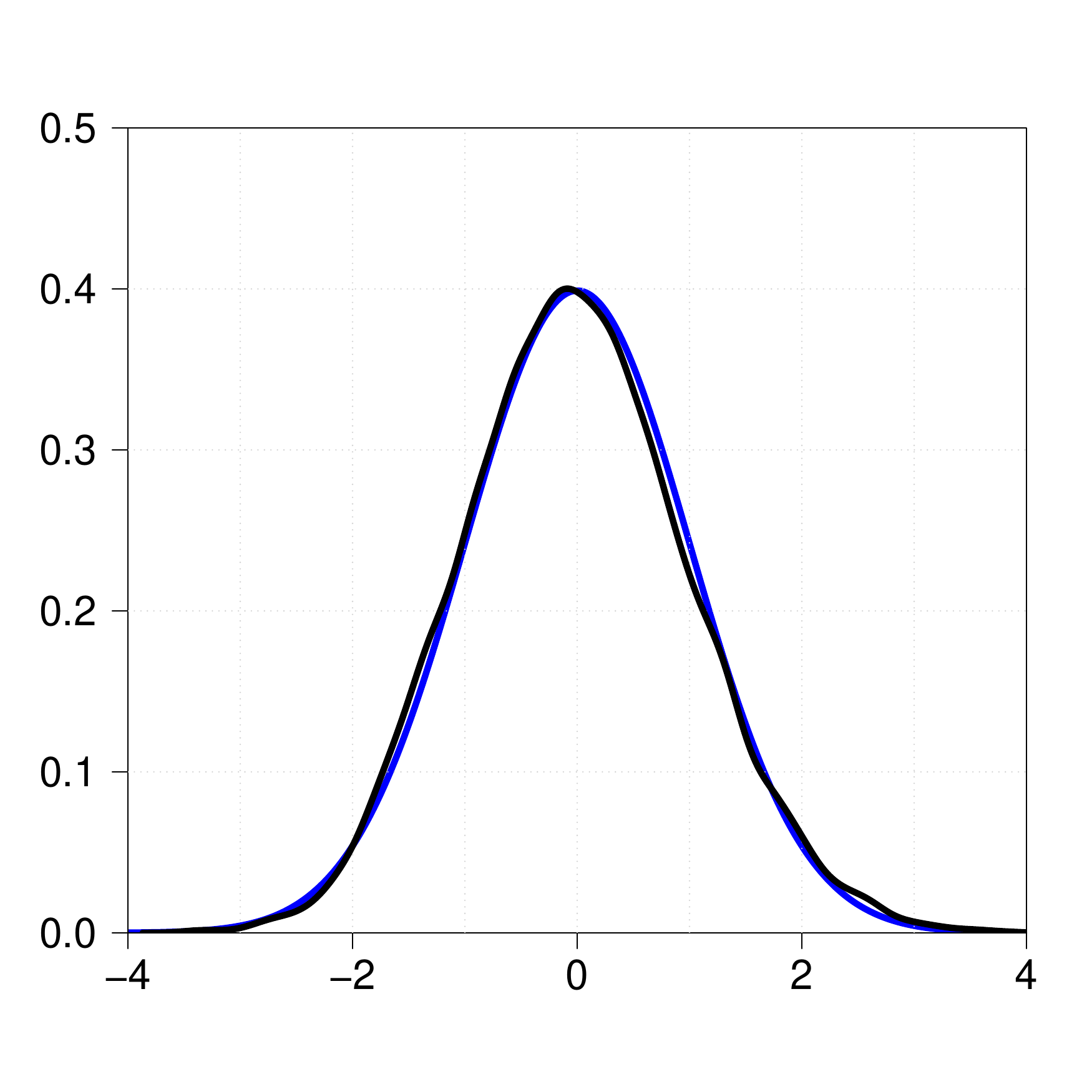}
       \vskip-.2cm
 \caption{$\lambda = 400$.} 
 \end{subfigure}
    \caption{Convergence of $3$-hop counts using probability density functions.} 
  \label{f1}
  \end{figure}

\begin{figure}[H]
  \centering
 \begin{subfigure}[b]{0.49\textwidth}
    \includegraphics[width=1\linewidth, height=5cm]{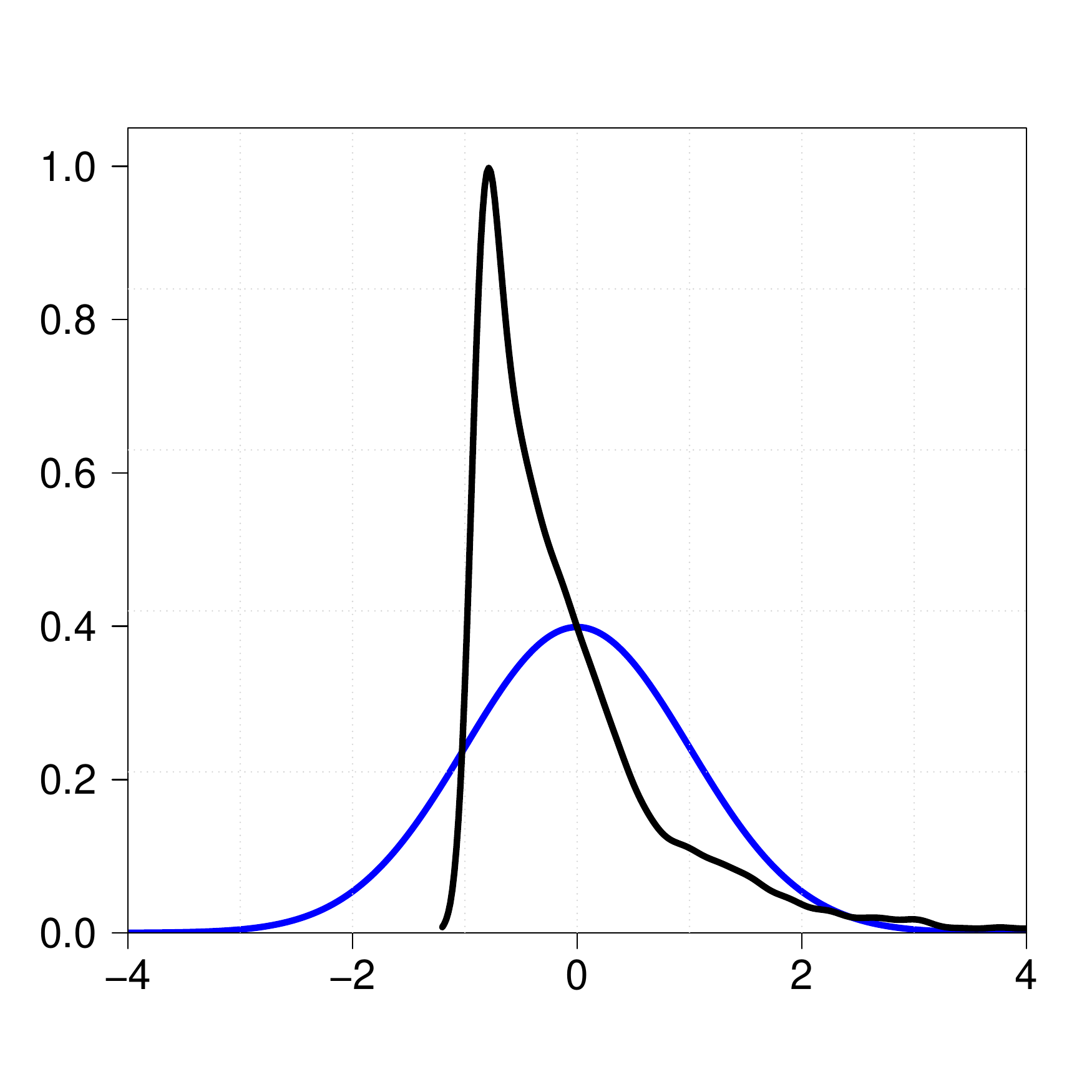}
        \vskip-.2cm
\caption{$\lambda = 5$.} 
 \end{subfigure}
 \begin{subfigure}[b]{0.49\textwidth}
    \includegraphics[width=1\linewidth, height=5cm]{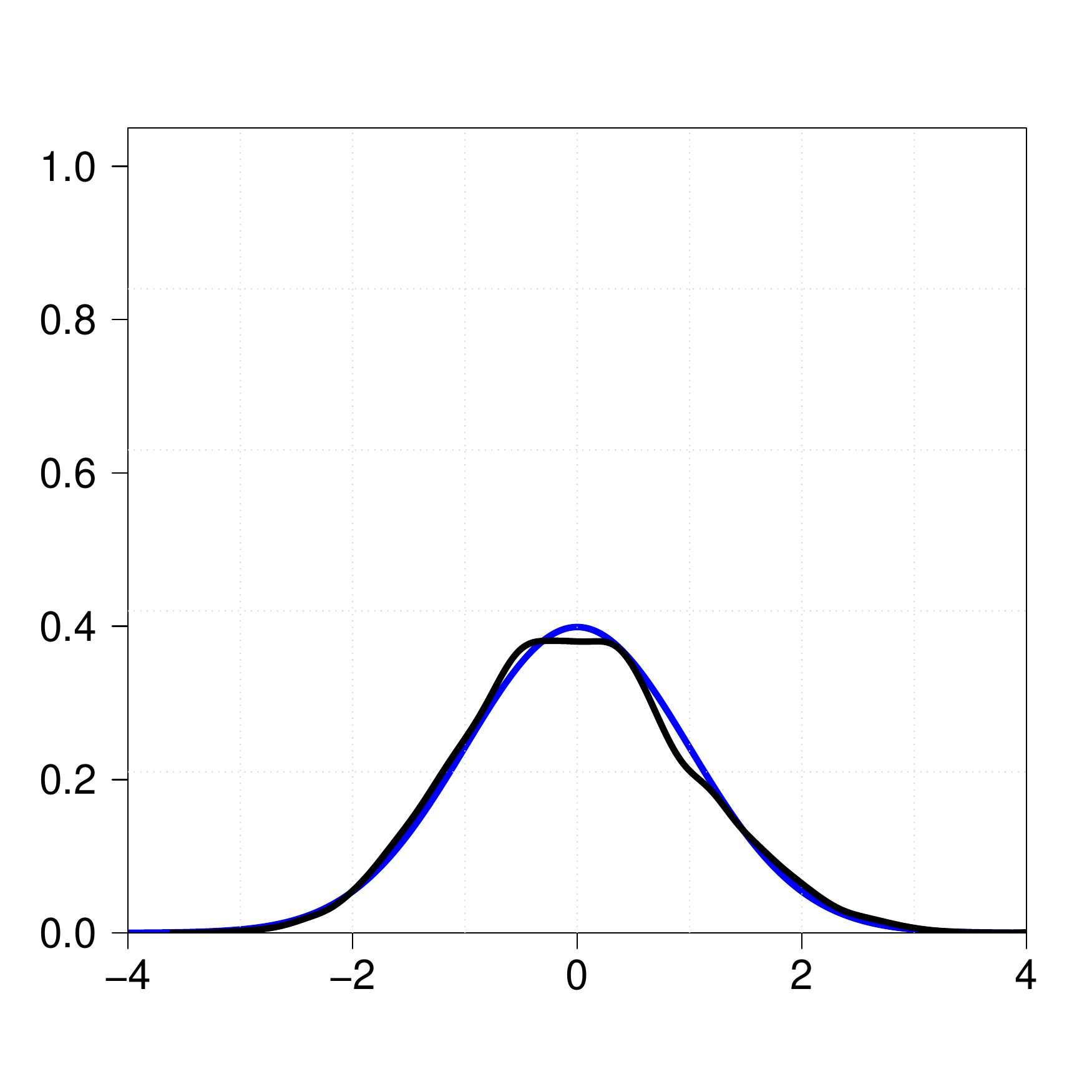}
    \vskip-.2cm
    \caption{$\lambda = 400$.} 
 \end{subfigure}
    \caption{Convergence of $4$-hop counts using probability density functions.} 
\label{f2}
  \end{figure}

 In addition, given $k\geq 2$ and $t\in [(k-1)r,kr)$,
 we consider the renormalized $k$-hop count
$$
\widetilde{\sigma}_k^{(\lambda )} (t):= \frac{\sigma_k^{(\lambda )} (t)- \E_\lambda [\sigma_k(t)]}{\sqrt{\Var_\lambda  [ \sigma_k(t)]}}. 
$$
From Proposition~\ref{1fdjkl}, for any $t \in [(k-1)r,kr)$, the skewness
 of $\sigma_k (t )$ satisfies  
$$
 \frac{\E_\lambda [ ( \sigma_k( t ) -  \E  [ \sigma_k(t ) ] )^3 ]}{( \Var_\lambda [ \sigma_k ( t )] )^{3/2}}
 =
 \E \big[ \big( \widetilde{\sigma}^{(\lambda )}_k( t ) \big)^3 \big] 
 =
 \frac{
   c_{k,3}^{(\lambda)} (kr - t ; kr - t ; kr - t )
 }{( \Var_\lambda [ \sigma_k (t)] )^{3/2}}
 \approx \frac{1}{\sqrt{\lambda}}. 
$$
 By Theorem~1 in \cite{Janson1988}, this shows the convergence in distribution
of $\widetilde{\sigma}_k^{(\lambda )} (t)$
to the standard normal distribution ${\cal N}(0,1)$ as $\lambda$ tends to infinity,
as illustrated in Figures~\ref{f1}-\ref{f2} using empirical probability
density plots.

\vspace{-0.6cm}

\vspace{-0.6cm}

  \vspace{-0.4cm}

\begin{figure}[H]
  \centering
 \begin{subfigure}[b]{0.49\textwidth}
    \includegraphics[width=1\linewidth, height=5cm]{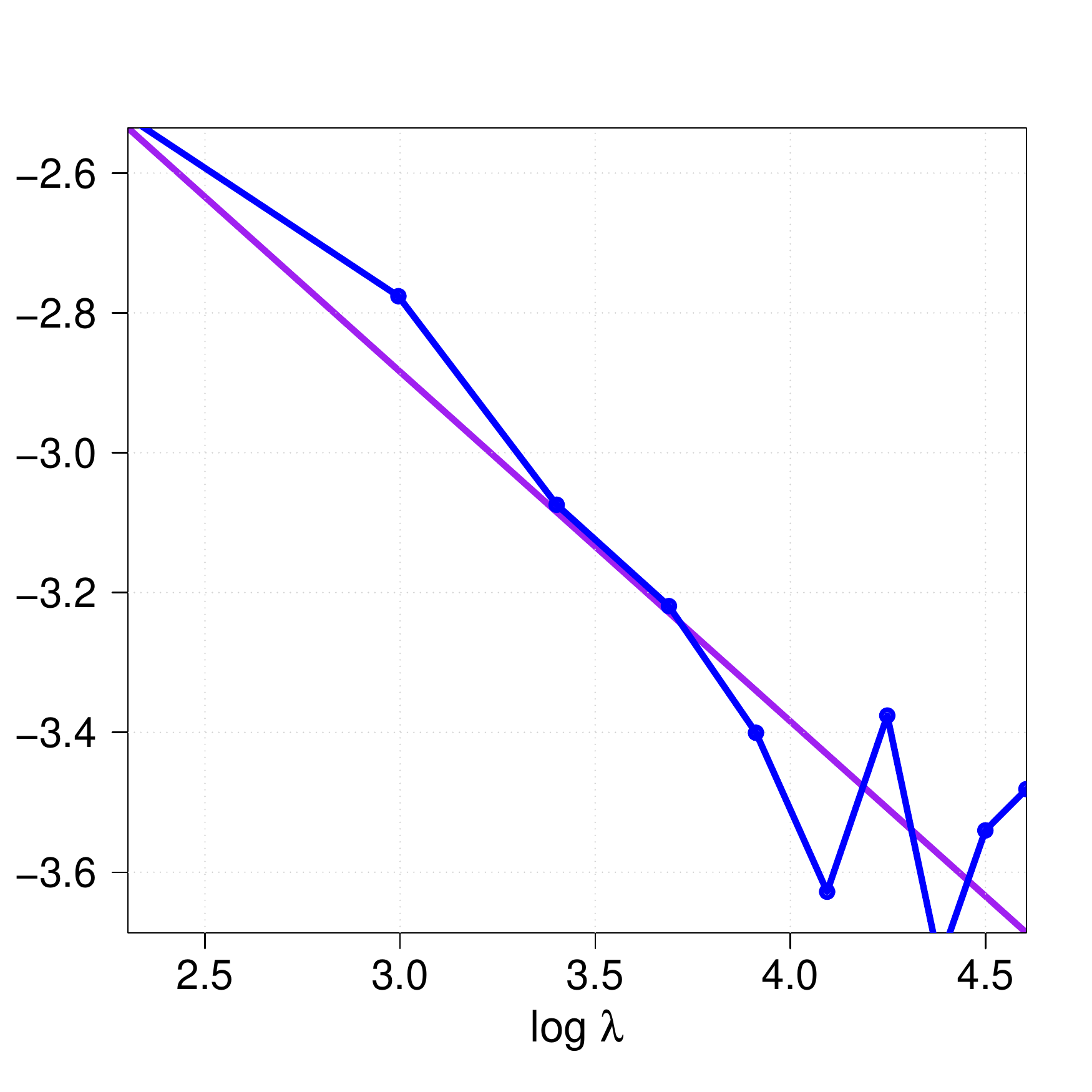}
    \vskip-.2cm
    \caption{Three-hop counts.} 
 \end{subfigure}
 \begin{subfigure}[b]{0.49\textwidth}
    \includegraphics[width=1\linewidth, height=5cm]{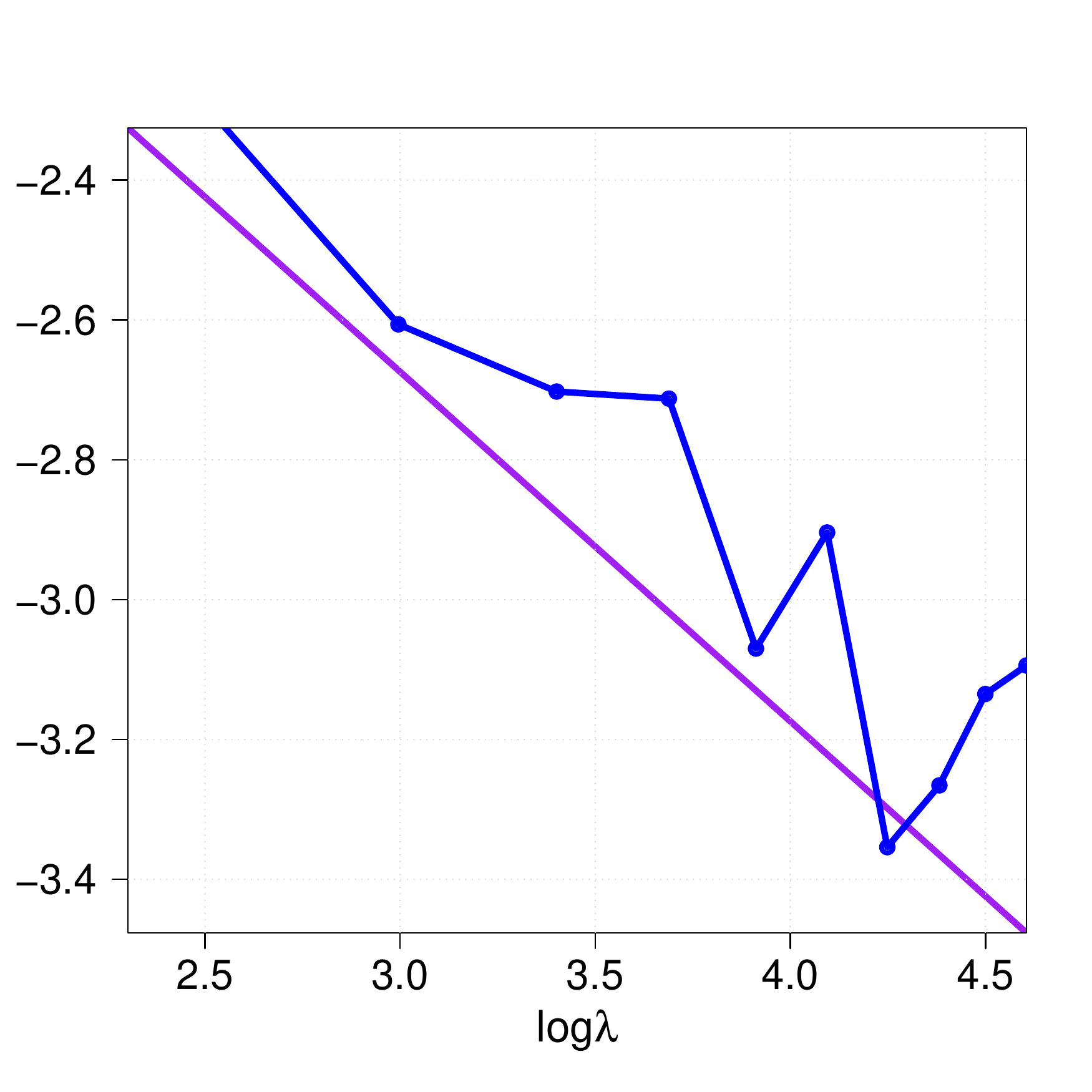}
    \vskip-.2cm
    \caption{Four-hop counts.} 
 \end{subfigure}
    \caption{Log-log plot of Kolmogorov distances.} 
  \label{f3}
  \end{figure}

\noindent
 In addition, by \eqref{fjkdls45} and \eqref{fjkls343} we have 
$$
 \frac{
 c_{k,n}^{(\lambda)} (kr - t ; \ldots ; kr - t )}{( \Var_\lambda  [ \sigma_k(t)] )^{n/2}} 
 = 
(B_n)^{k-2} O ( \lambda^{1-n/2} )
 \leq 
 (n!)^{k-2} O ( \lambda^{1-n/2} )
, \quad n \geq 2,
$$ 
 hence the Statulevi\v{c}ius condition 
 is satisfied with $\gamma := k-3$ and $\Delta := \sqrt{\lambda}$,
 which,
 by \cite{rudzkis},
   Corollary~2.1 in \S~1.3 in \cite{saulis},
   see also Theorem~2.4 in \cite{doering}, yields the Berry-Esseen bound
   \begin{equation}
   \label{fjkf} 
 d_K \big( \widetilde{\sigma}_k^{(\lambda )} (t) , {\cal N} \big) 
\leq \frac{C(k,r)}{\lambda^{1/(2 + 4(k-3))}}
\end{equation} 
   for $t\in [(k-1)r,kr)$ as $\lambda$ tends to infinity.
 More precisely, \eqref{fjkf} can be improved as in the following
 proposition.
\begin{prop}
  \label{fsklf34}
  Let $k\geq 2$ and $t\in [(k-1)r,kr)$.
  The renormalized $k$-hop count 
$\widetilde{\sigma}_k^{(\lambda )} (t)$ 
satisfies the Wasserstein and Kolmogorov bounds 
\begin{equation}
\label{fjkldsf-2} 
d_{K/W} \big( \widetilde{\sigma}_k^{(\lambda )} (t) , {\cal N} \big)
\leq \frac{C(k,r)}{\sqrt{\lambda}}
\end{equation} 
 for some constant $C(k,r)>0$ as $\lambda$ tends to infinity. 
\end{prop}
\begin{Proof}
 The kurtosis of $\sigma_k^{(\lambda )}(t )$ satisfies 
\begin{eqnarray*} 
 \frac{\E_\lambda [ ( \sigma_k(t) -  \E_\lambda [ \sigma_k(t) ] )^4 ]}{( \Var_\lambda [ \sigma_k (t)] )^2} - 3
 & = & 
 \E_\lambda \big[ \big( \widetilde{\sigma}_k^{(\lambda )}(t) \big)^4 \big] - 3 
 \\
 & = & 
 \frac{c_{k,4}^{(\lambda)} (kr-t;kr-t;kr-t;kr-t)}{( \Var_\lambda [ \sigma_k (t)] )^2} 
 \\
  & = & 
(B_4)^{k-2} O ( \lambda^{-1} ), 
\end{eqnarray*} 
 as $\lambda$ tends to infinity. 
 The Kolmogorov distance bound in \eqref{fjkldsf-2} 
 then follows from the fourth moment theorem for $U$-statistics and sums of
 multiple stochastic integrals Corollary~4.10 in \cite{eichelsbacher}
 applied to \eqref{fjkldsf1},
 see also Theorem~3 in \cite{lachieze-rey}. 

 \medskip

 Regarding the Wasserstein distance bound,
 according to \eqref{fjda} we can the represent $\sigma_k^{(\lambda )} (t)$ 
 as the $U$-statistics  
$$ 
 \sigma_k^{(\lambda )} (t) 
 = 
 \sum_{((x_1,l_1),\ldots , ((x_{k-1},l_{k-1}) ) \in \omega^{k-1}
   \atop (x_i,l_i)\not=(x_j,l_j), 1\leq i\not= j \leq d}
 \tilde{f}_{\lambda ( kr -t) } (x_1/\lambda ,l_1;\ldots ; x_{k-1}/\lambda ,l_{k-1} )
$$
 of order $k-1$, where $\tilde{f}_t:( [0,r] \times \{1,\ldots , k-1\})^{k-1} \to \{0,1\}$,
 given by 
$$
\tilde{f}_t (x_1,l_1;\ldots , x_{k-1};l_{k-1})
:= \frac{1}{(k-1)!} \prod_{i=0}^{k-1} {\bf 1}_{\{(l_{i+1}-l_i)x_i<(l_{i+1}-l_i)x_{i+1} \}}, 
$$ 
 $((x_1,l_1),\ldots , (x_{k-1},l_{k-1}))\in ([0,r]\times \{1,\ldots , k-1\})^{k-1}$
 is the symmetrization in
 $k-1$ variables in $[0,r]\times \{1,\ldots , k-1\}$ of $f_\tau$. 
  Theorem~4.7 in \cite{reitzner} yields the bound 
  $$
  d_W\big(\widetilde{\sigma}_k^{(\lambda )}, {\cal N}\big) \leq
  \sum_{1\leq i \leq j \leq r}
  \frac{\sqrt{M_{i,j}}}{\Var_\lambda [ \sigma_k(t)]}
 $$
 where   
 $M_{i,j}$ is defined in (14) therein satisfies
$$
 M_{1,1} \leq (k-1)^4 
 ( \lambda r)^{4(k-1)-3}, \qquad  i,j=1. 
$$ 
 and
  $$
 M_{i,j} \leq {k-1 \choose i}^2 {k-1 \choose j}^2
 ( \lambda r)^{4(k-1)-i-j}, \qquad  
 2\leq i \leq j \leq r. 
$$  
 Hence by \eqref{fjkdls45} and \eqref{dfjkvar} we have
\begin{eqnarray*} 
  \lefteqn{
    d_W\big(\widetilde{\sigma}_k^{(\lambda )}, {\cal N}\big)
  }
  \\
   & \leq & 
  \frac{1}{\Var_\lambda  [ \sigma_k(t)]}
  \left(
  (k-1)^2 ( \lambda r)^{2(k-1)-3/2}
  +
  \sum_{2\leq i \leq j \leq r}
  {k-1 \choose i} {k-1 \choose j}
  ( \lambda r)^{2(k-1)-i/2-j/2}
  \right)  
  \\
  & \leq & 
  \frac{C(k,r)}{\sqrt{\lambda r}} 
  +
  C(k,r) \sum_{2\leq i \leq j \leq r}
  ( \lambda r)^{1-i/2-j/2}. 
\end{eqnarray*}
 The above conclusions can also be reached by 
 noting that $\widetilde{\sigma}_k^{(\lambda )}$ admits a Hoeffding decomposition
 and by applying Theorem~1.3 in \cite{doblerpeccati} for the Wasserstein
 distance, or Theorem~6.3 in \cite{PS4} 
 for the Kolmogorov distance, which refine the central limit
 theorem of \cite{dejong1990}.
\end{Proof}
\noindent 
Figure~\ref{f3} presents numerical estimates that are consistent with the rate in 
\eqref{fjkldsf-2}, by plotting 
$\log d_K\big( \widetilde{\sigma}^{(\lambda )}_k (t) , {\cal N} \big)$ against $\log \lambda$
and their comparison with the line of slope $-1/2$.
 Kolmogorov distances $d_K$ have been estimated in R using the distrEx package.

\vspace{-0.4cm}

  \vspace{-0.3cm}

\appendix

\section{Computer codes}

\subsubsection*{Computation of joint moments}
\noindent
Explicit moment expressions are obtained using the following
Mathematica code which implements the recursion \eqref{fdfdf}
of Proposition~\ref{433} in two steps. 
 First, the computation of the multiple integral in
 $ \lambda_k (du_1) \cdots \lambda_k (du_l)$ in the following code.

 \medskip
  
\begin{lstlisting}[language=Mathematica, caption={1}, title={Mathematica Code~\nameref{code1}.}, label={code1}]
Needs["Combinatorica`"]
fseq[n_] := fseq[n] = (Module[{k, tmp}, tmp = {}; 
    If[n == 1, Return[{{1}}], For[k = 1, k <= n, k++, Do[tmp = Prepend[tmp, Sort[Append[a, k]]], {a, fseq[n - 1]}]]; tmp]])
r1[x__, t__, t2__, t3__, p_, f_] := (Module[{s}, 
   If [p == 1, Integrate[f[Prepend[x, s]], {s, 0, t[[2]]}], 
    If[t2[[p]] === t2[[p + 1]], 
     Integrate[r1[Prepend[x, s], ReplacePart[t, p -> s], t2, t3, p - 1, f], {s, t3[[p]], t[[p + 1]]}], 
     Integrate[r1[Prepend[x, s], t2, t2, t3, p - 1, f], {s, t3[[p]], t[[p + 1]]}]]]])
\end{lstlisting}

\vskip-0.2cm

\noindent
 This is followed by the summation over
 $\pi_1\cup \cdots \cup \pi_l = \{1,\ldots , n\}$ and
 $l=1,\ldots ,n$ in the next code.
 
 \medskip
 
\begin{lstlisting}[language=Mathematica, caption={2}, title={Mathematica Code~\nameref{code2}.}, label={code2}]
r2[x__, t1__, t2__, t3__, p_, q__, f_] := r2[x, t1, t2, t3, p, q, f] = (Module[{s}, If [p == 1, Integrate[f[Join @@ MapThread[Table, {Prepend[x, s], q}]], {s, 0, 
       t1[[2]]}], 
     If[t2[[p]] === t2[[p + 1]], 
      Integrate[r2[Prepend[x, s], ReplacePart[t1, p -> s], t2, t3, p - 1, q, f], {s, t3[[p]], t1[[p + 1]]}], 
      Integrate[r2[Prepend[x, s], t2, t2, t3, p - 1, q, f], {s, t3[[p]], t1[[p + 1]]}]]]])
mk[t__, lambda__] := mk[t, lambda] = (Module[{b, z1, q, n, m, z, zz, tmp},
    If[lambda == {}, Return[1]]; m = Length[t]; 
    If[Length[lambda] == 1, b[l__] = 1, b[l__] := (c[d__] := mk[d, Drop[lambda, -1]]; 
       With[{e = c[Array[s, Length[t]]]}, 
        h[z_] := Block[{s}, s[i_] := z[[i]]; e]; Return[h[l]]])]; 
    tmp = 0; Do[n = Length[pp]; q = Map[Length, pp]; z = t[[Map[Min, pp]]]; 
     zz = Prepend[Drop[z, -1], 0];
     Do[Do[If[r == r[[Ordering@p0]], z1 = Prepend[z[[r]], 0]; 
        tmp += Last[lambda]^n*r2[{}, z1, z1, zz[[r]], n, q[[Ordering@p0]], b]], {r, fseq[n]}], {p0, Permutations[Range[n]]}], {pp, SetPartitions[m]}]; Return[Expand[Flatten[{tmp}][[1]]]]])
\end{lstlisting}

\vspace{-0.3cm}

\noindent
 The joint moment $\E [ \sigma_k(kr-\tau_1)\cdots \sigma_k(kr-\tau_n) ]$
 is then computed from the command
 ${\rm mk}[\{ \tau_1, \ldots , \tau_n\}, \{\lambda_1,\ldots , \lambda_{k-1}]$.

\subsubsection*{Computation of joint cumulants}

Explicing cumulant expressions are obtained using the following
Mathematica code which implements the recursion \eqref{fjkla} 
 of Proposition~\ref{fjhkds} in two steps.
 First, the computation of the multiple integral in
 $ \lambda_k (du_1) \cdots \lambda_k (du_l)$ in the following code.

 \medskip
  
\begin{lstlisting}[language=Mathematica, caption={3}, title={Mathematica Code~\nameref{code3}.}, label={code3}]
Needs["Combinatorica`"]
fseq2[n_] := fseq2[n] = (Module[{k, tmp, tmp2}, tmp = {}; 
    If[n == 1, Return[{{1}}], 
     For[k = 1, k <= n, k++, Do[tmp2 = Join[{k}, a]; tmp = Append[tmp, ReverseSort[tmp2]], {a, fseq2[n - 1]}]]; tmp]])
r1c[x__, t__, t2__, t3__, p_, f_] := (Module[{s}, 
   If [p == 1, Integrate[f[Prepend[x, s]], {s, 0, t[[2]]}], 
    If[t2[[p]] === t2[[p + 1]], Integrate[r1c[Prepend[x, s], ReplacePart[t, p -> s], t2, t3, p - 1, f], {s, t3[[p]], t[[p + 1]]}], Integrate[r1c[Prepend[x, s], t2, t2, t3, p - 1, f], {s, t3[[p]], t[[p + 1]]}]]]])
\end{lstlisting}

\vskip-0.2cm

\noindent
This is followed by the recursive computation of
$c_{k,l}^{(\lambda)} (\widebar{s}_{\pi_1} ; \ldots ; \widebar{s}_{\pi_l} )$
by the induction relation of Proposition~\ref{jkld},
and the summation
 \eqref{fjkla} over
 $\pi_1\cup \cdots \cup \pi_l = \{1,\ldots , n\}$ and
 $l=1,\ldots ,n$ 
 in the following code

 \medskip
 
\begin{lstlisting}[language=Mathematica, caption={4}, title={Mathematica Code~\nameref{code4}.}, label={code4}]
ck[t__, l3__, lambda__] := ck[t, l3, lambda] = (Catch[
    Module[{b, l4, k, n, l5, l6, t5, tmp3, tmp4, oo, z, zz, z1, z2, u, q}, oo = 0; n = Length[t]; 
     Do[oo += 1; If[l0 > 1, l4 = l3; l4[[oo]] += -1; 
       l5 = Join[Range[oo], {oo}, Range[oo + 1, n + 1]]; 
       l6 = Join[l4[[Range[oo]]], {1}, l4[[Range[oo + 1, n]]]]; 
       t5 = Join[t[[Range[oo]]], {t[[oo]]}, t[[Range[oo + 1, n]]]]; 
       tmp4 = 0; Do[If[MemberQ[l5[[pi[[1]]]], oo] && MemberQ[l5[[pi[[2]]]], oo],
          tmp4 += ck[t5[[pi[[1]]]], l6[[pi[[1]]]], lambda]*ck[t5[[pi[[2]]]], l6[[pi[[2]]]], lambda]], {pi, KSetPartitions[Range[n + 1], 2]}]; tmp4 += ck[t5, l6, lambda]; Throw[tmp4]], {l0, l3}]; k = Length[lambda]; If[k == 1, Throw[lambda[[1]]*t[[1]]]]; 
     If [Length[lambda] == 1, b[l__] = 1, b[l__] := (c[d__] := ck[d, q0, Drop[lambda, -1]]; 
        With[{e = c[Array[s4s, Length[l]]]}, fl[z7_] := Block[{s4s}, s4s[i_] := z7[[i]]; e]; 
         Return[fl[l]]])]; m = Length[t]; q0 = {m}; 
     tmp3 = lambda[[1]]*r1c[{}, {0, t[[1]]}, {0, t[[1]]}, {0}, 1, b]; 
     For[n = 2, n <= m, n++, Do[zz = {0}; q = Map[Length, pp];
       z = t[[Map[Min, pp]]]; zz = Prepend[Drop[z, -1], 0]; 
       Do[q0 = q[[Ordering@p0]]; Do[If[bb == bb[[Ordering@p0]], z1 = Prepend[z[[Reverse[bb]]], 0]; 
          tmp3 += lambda[[1]]^n*r1c[{}, z1, z1, zz[[Reverse[bb]]], n, b]], {bb, 
          fseq2[n]}], {p0, Permutations[Range[n]]}], {pp, KSetPartitions[m, n]}]];
     Return[Expand[Flatten[{tmp3}][[1]]]]]])
\end{lstlisting}

\vspace{-0.3cm}

\noindent
 The joint cumulant
 $c_{k,n}^{(\lambda)} (\widebar{\tau}_{\pi_1} ; \ldots ; \widebar{\tau}_{\pi_l} )
 = \kappa_\lambda \big( ( \sigma_k(kr - \tau_1) )^{|\pi_1|}, \ldots , ( \sigma_k(kr - \tau_l) )^{|\pi_l|}\big)$
 is computed from the command
 {\rm ck}[$\{ \tau_1 , \ldots , \tau_l \}, \{ |\pi_1| , \ldots , |\pi_l| \} , \{ \lambda_1 , \ldots , \lambda_{k-1} \} ]$,
 and
 the joint cumulant 
 $c_{k,n}^{(\lambda)} (\tau_1 ; \ldots ; \tau_n )
 = \kappa_\lambda \big( \sigma_k(kr - \tau_1) , \ldots , \sigma_k(kr - \tau_n) \big)$
 is computed from the command
 {\rm ck}[$\{ \tau_1 , \ldots , \tau_n \}, \{ 1 , \ldots , 1 \} , \{ \lambda_1 , \ldots , \lambda_{k-1} \} ]$. 
\section{Moments of Poisson stochastic integrals}
In this section we review the background results on the
moments of Poisson stochastic integrals that are used in this paper.
Let $\omega (dx)$ denote a Poisson point process
with intensity measure $\mu (dx)$ on a measure space $X$.
The next proposition, see Proposition~3.1 in \cite{momentpoi}
 or Theorem~1 and Proposition~7 in \cite{prob},
 provides a moment identity for Poisson stochastic integrals
 with random integrands using sums over partitions. 
\begin{prop} 
\label{pr11-2} 
 Let $f_1,\ldots , f_n: X \times \Omega \longrightarrow \real$ 
 be deterministic functions, $n\geq 1$. 
 We have 
\begin{eqnarray} 
\nonumber 
\lefteqn{ 
 \E \left[ 
 \int_X f_1(x) \omega (dx) 
 \cdots 
 \int_X f_n(x) \omega (dx) 
 \right] 
} 
\\ 
\nonumber 
 & = & 
 \sum_{l=1}^n 
 ~ 
 \sum_{\pi_1 \cup \cdots \cup \pi_l = \{ 1, \ldots , n \} } 
  \int_{X^l} 
 \left( 
 \prod_{j=1}^l 
 \prod_{i\in \pi_j}
 f_i (x_j , \xi ) 
 \right)  
 \mu (dx_1) \cdots \mu (dx_l) 
 , 
\end{eqnarray} 
 where the sum runs over all partitions 
 $\pi_1,\ldots , \pi_l$ of $\{ 1 , \ldots , n \}$ of sizes 
 and $|\pi_l|$ denotes the cardinality of each block $\pi_l$,
 $l=1,\ldots , n$. 
\end{prop} 
Proposition~\ref{pr11-2} extends in the next Proposition~\ref{p2.2}
 as a moment identity for the stochastic integral of functions
 $f(z_1,\ldots z_p )$ of $p$ variables
 $z_1,\ldots z_p \in X^p$, see in Theorem~3.1 of \cite{bogdan}.
For this, let $\Pi [n\times p]$ denote the set of partitions
of 
 $$
 [n\times p] := \{1,\ldots , n\} \times \{1,\ldots , p\}
 =
 \big\{ (i,j) \ : \
 i=1,\ldots , n, \ j = 1,\ldots , p \big\}, 
 $$
 identified to $\{1,\ldots , np\}$, and let
$\rho : = (\rho_1,\ldots , \rho_n) \in \Pi [n\times p]$ denote the
 partition made of the $n$ blocks
$\rho_i := \{ (i,1), \ldots , (i,p )\}$ of size $p$, for 
$i=1,\ldots , n$.
\begin{prop} 
  \label{p2.2}
   For $f : X^p \times \Omega \longrightarrow \real$ be 
 a sufficiently integrable function of $p$ variables we have 
\begin{eqnarray} 
\nonumber 
\lefteqn{
  \E  \left[ 
    \left(
    \int_X
       \cdots \int_X
          f(x_1,\ldots , x_p; \omega ) \omega ( dx_1) \cdots
   \omega ( dx_p ) \right)^n 
 \right] 
}
\\
\nonumber
& = &  
  \sum_{ \pi \in \Pi [n\times r] } 
    \int_X
   \cdots \int_X
   \prod_{l=1}^n
   u \big(x_{\zeta^\pi_{l,1}}, \ldots , x_{\zeta^\pi_{l,p}}\big) 
   \mu ( d x_1 ) 
   \cdots
   \mu ( d x_{|\pi|}) 
   , 
\end{eqnarray} 
 where 
$
\zeta^\pi_{i,j}$ is the index of the block
 of $\pi$ that contains $(i,j)$. 
\end{prop} 
Proposition~\ref{p2.2} can also be extended 
as a joint moment identity for multiparameter processes
in the next proposition.
\begin{prop} 
  \label{p2.2-3}
  Let $f_i : X^p \longrightarrow \real$,
  $i=1,\ldots , n$ be sufficiently integrable functions of $p$ variables.
  We have 
\begin{eqnarray} 
\nonumber 
\lefteqn{
  \E  \left[ 
\prod_{i=1}^n 
    \int_X
      \cdots \int_X
         f_i (x_1,\ldots , x_p; \omega ) \omega ( dx_1) \cdots
   \omega ( dx_p ) 
 \right] 
}
\\
\nonumber
& = &  
  \sum_{ \pi \in \Pi [n\times r] } 
    \int_X
   \cdots \int_X
    \prod_{j=1}^{|\pi|}
   \prod_{i\in \pi_j}
   f_i \big(x_{\zeta^\pi_{l,1}}, \ldots , x_{\zeta^\pi_{l,p}}\big) 
   \mu ( d x_1 ) 
   \cdots
   \mu ( d x_{|\pi|}) 
   , 
\end{eqnarray} 
 where 
$
\zeta^\pi_{i,j}$ is the index of the block
 of $\pi$ that contains $(i,j)$. 
\end{prop} 
 In case the function $f(x_1,\ldots , x_p;\omega )$ vanishes on the diagonals in $X^p$,
 the integral of $f$ rewrites as the $U$-statistics 
 \begin{equation}
   \label{fjlf} 
 \int_{X^r} f(x_1,\ldots , x_p;\omega ) \omega ( dx_1 ) \cdots
 \omega ( dx_p ) 
 =
 \sum_{(x_1,\ldots , x_p ) \in \omega^p
   \atop x_i\not=x_j, 1\leq i\not= j \leq p}
 f(x_1,\ldots , x_p;\omega ), 
\end{equation}
 we have the following corollary of Proposition~\ref{p2.2-3}.
 For this, we let $\rho:=(\rho_1,\ldots , \rho_n)$ be the partition of $[n\times p]$ made
 of the blocks $\rho_i=((i,j))_{j=1,\ldots ,p}$, $i=1,\ldots , n$,
 and we say that a partition $\pi$ of $[n \times p]$
 is {\em non-flat}, i.e. $\pi \wedge \rho = \widehat{0}$,
 if every block of $\pi$ contains at most one element
 of $\rho_i$, $i=1,\ldots , n$, where 
 $\widehat{\bf 0} : = \{\{1\},\ldots , \{n\}\}$ is 
 the $n$-block partition of $\{1,\ldots , n\}$. 
\begin{corollary}
   \label{c1-a}
    Let $f_i : X^p \longrightarrow \real$,
    $i=1,\ldots , n$ be sufficiently integrable functions of $p$ variables
    such that the function $f_i(x_1,\ldots , x_p;\omega )$ vanishes on the diagonals in $X^p$, $i=1,\ldots , n$. 
 Then, the joint moments of the $U$-statistics \eqref{fjlf} can be computed as 
\begin{eqnarray} 
\nonumber 
\lefteqn{
  \E  \left[ 
\prod_{i=1}^n 
    \int_X
       \cdots \int_X
          f_i(x_1,\ldots , x_p; \omega ) \omega ( dx_1) \cdots
   \omega ( dx_p ) 
 \right] 
}
\\
\nonumber
& = &  
  \sum_{
   \substack{
     \pi \in \Pi [n\times p] 
     \\
     \pi \wedge \rho = \hat{0}
   }}
    \int_X
     \cdots \int_X
     \prod_{j=1}^{|\pi|} 
   \prod_{i \in \rho_j}
   u \big(x_{\zeta^\pi_{i,1}}, \ldots , x_{\zeta^\pi_{i,p}}\big) 
   \mu ( d x_1 ) 
   \cdots
   \mu ( d x_{|\pi|}).
\end{eqnarray} 
\end{corollary} 

\section{Explicit moment recursions}
\label{s9}
 In this section we confirm the joint moment
 induction of Proposition~\ref{433} via explicit calculations
 for $n=1,2,3$. 
\subsubsection*{First moment recursion}
 By \eqref{kldsf} we have 
\begin{eqnarray*} 
  m_{k+1,1}^{(\lambda)}(\tau )
 & = & 
  \E_\lambda [ \sigma_{k+1}((k+1)r- \tau ) ] 
  \\
   & = &  
  \E_\lambda \big[ Z^{(k+1)}_\tau \big] 
  \\
   & = & 
\int_0^\tau
\E_\lambda \big[ Z^{(k)}_s \big]
\lambda_k(ds)
\\
 & = & 
\int_0^\tau m_{k,1}^{(\lambda)} ( s ) \lambda_k(ds), 
\end{eqnarray*} 
 which recovers \eqref{fjkslff} as
$$
  m_{k,1}^{(\lambda)} (\tau )
=
\int_0^\tau \int_0^{s_{k-1}} \cdots
\int_0^{s_2}
\lambda_1(ds_1)\cdots \lambda_{k-1} (ds_{k-1}).  
$$

\subsubsection*{Second joint moment recursion}
For $n=2$, by \eqref{kldsf} the induction relation \eqref{kldsf-2} reads 
\begin{eqnarray*} 
  m_{k+1,2}^{(\lambda)} (\tau_1,\tau_2) & = & 
  \E_\lambda [ \sigma_{k+1}((k+1)r- \tau_1) \sigma_k((k+1)r-\tau_2) ]
  \\
  & = & 
  \E_\lambda \big[ Z^{(k+1)}_{\tau_1} Z^{(k+1)}_{\tau_2} \big] 
 \\
  & = & 
\E_\lambda \left[ 
\int_0^{\tau_1} Z^{(k)}_{s_1} dN^{(k)}_{s_1}
\int_0^{\tau_2} Z^{(k)}_{s_2} dN^{(k)}_{s_2}
\right]
\\ 
& = & 
\int_0^{\tau_1}
\E_\lambda \big[ \big( Z^{(k)}_{s_1} \big)^2 \big]
\lambda_k(ds_1) 
+ \int_0^{\tau_2} \int_0^{\tau_1}
\E_\lambda \big[ Z^{(k)}_{s_1} Z^{(k)}_{s_2} \big]
\lambda_k(ds_1) \lambda_k(ds_2)
\\
 & = & 
  \int_0^{\tau_1}
  m_{k,2}^{(\lambda)} (u_1,u_1)
  \lambda_k (du_1) 
 + 
  \int_0^{\tau_2} 
  \int_0^{\tau_1}
  m_{k,2}^{(\lambda)} (u_1, u_2 ) 
  \lambda_k (du_1) \lambda_k (du_2)
. 
\end{eqnarray*} 
\subsubsection*{Third joint moment recursion}
For $ \tau_1,\tau_2 , \tau_3 \in [0,r]$ the recursion \eqref{kldsf-2}
reads 
\begin{align*}
&
  m_{k+1,3}^{(\lambda)} (\tau_1,\tau_2, \tau_3)
  =
  \E_\lambda [ \sigma_{k+1}((k+1)r- \tau_1) \sigma_k((k+1)r-\tau_2) \sigma_k((k+1)r-\tau_3) ]
  \\
   & = 
  \E_\lambda \big[ Z^{(k+1)}_{\tau_1} Z^{(k+1)}_{\tau_2} Z^{(k+1)}_{\tau_3} ]
  \\
   & = 
  \E_\lambda \left[
    \int_0^{\tau_1} Z^{(k)}_u dN^{(k)}_u 
    \int_0^{\tau_2} Z^{(k)}_u dN^{(k)}_u 
    \int_0^{\tau_3} Z^{(k)}_u dN^{(k)}_u 
    \right]
  \\
   & = 
  \int_0^{\tau_1 \wedge \tau_2 \wedge \tau_3 }
  \E \big[ \big( Z^{(k)}_{s_1} \big)^3 \big] 
  \lambda_k(ds_1)
  \\
  & 
  \quad
  + 
  \int_0^{\tau_3}
  \int_0^{\tau_1 \wedge \tau_2}
  \E_\lambda \big[ \big( Z^{(k)}_{s_1} \big)^2 Z^{(k)}_{s_3} \big]
  \lambda_k(ds_1) \lambda_k(ds_3) 
  +
  \int_0^{\tau_2 \wedge \tau_3}
  \int_0^{\tau_1}
  \E_\lambda \big[ \big( Z^{(k)}_{s_1} \big)^2 Z^{(k)}_{s_2} \big]
  \lambda_k(ds_1) \lambda_k(ds_2) 
  \\
  & 
 \quad
   +
  \int_0^{\tau_2 \wedge \tau_3}
  \int_0^{\tau_1}
  \E_\lambda \big[ Z^{(k)}_{s_1} \big( Z^{(k)}_{s_2} \big)^2 \big]
  \lambda_k(ds_1) \lambda_k(ds_2) 
  + 
  \int_0^{\tau_3}
  \int_0^{\tau_2}
  \int_0^{\tau_1}
  \E_\lambda \big[ Z^{(k)}_{s_1} Z^{(k)}_{s_2} Z^{(k)}_{s_3} \big]
  \lambda_k(ds_1)\lambda_k(ds_2)\lambda_k(ds_3) 
\\
   & = 
 \int_0^{\tau_1 \wedge \tau_2 \wedge \tau_3 } 
  m_{k,3}^{(\lambda)} (u_1 , u_1 , u_1 ) 
  \lambda_k (du_1) 
     + 
  \int_0^{\tau_3} 
  \int_0^{\tau_1 \wedge \tau_2}
  m_{k,n}^{(\lambda)} (u_1,u_1,u_3) 
  \lambda_k (du_1) \lambda_k (du_3)
  \\
   &
\quad
    + 
  \int_0^{\tau_2} 
  \int_0^{\tau_1 \wedge \tau_3}
  m_{k,n}^{(\lambda)} (u_1,u_1,u_2) 
  \lambda_k (du_1) \lambda_k (du_2)
    + 
  \int_0^{\tau_2 \wedge \tau_3} 
  \int_0^{\tau_1}
  m_{k,n}^{(\lambda)} (u_1,u_2,u_2) 
  \lambda_k (du_1) (du_2)
  \\
  &
  \quad
  + 
  \int_0^{\tau_3} 
  \int_0^{\tau_2} 
  \int_0^{\tau_1} 
  m_{k,n}^{(\lambda)} (u_1,u_2,u_3) 
  \lambda_k (du_1) \lambda_k (du_2) \lambda_k (du_3)
  \\
 & =  
 \int_0^{\tau_1 \wedge \tau_2 \wedge \tau_3 }
  m_{k,3}^{(\lambda)} (u_1, u_1,u_1 ) 
  \lambda_k (du_1) 
    + 
  \int_0^{\tau_3} 
  \int_0^{\tau_1 \wedge \tau_2}
  m_{k,n}^{(\lambda)} (u_1,u_1,u_3) 
  \lambda_k (du_1) \lambda_k (du_3)
  \\
  & 
\quad
      + 
  \int_0^{\tau_2} 
  \int_0^{\tau_1 \wedge \tau_3}
  m_{k,n}^{(\lambda)} (u_1,u_2,u_1) 
  \lambda_k (du_1) \lambda_k (du_2)
      + 
  \int_0^{\tau_2 \wedge \tau_3} 
  \int_0^{\tau_1}
  m_{k,n}^{(\lambda)} (u_1,u_2,u_2) 
  \lambda_k (du_1) (du_2)
  \\
  & \quad
   + 
  \int_0^{\tau_3} 
  \int_0^{\tau_2} 
  \int_0^{\tau_1} 
  m_{k,n}^{(\lambda)} (u_1,u_2,u_3) 
  \lambda_k (du_1) \lambda_k (du_2) \lambda_k (du_3). 
\end{align*}
\section{Explicit cumulant recursions}
\label{s10}
 In this section we confirm the joint moment
 induction of Proposition~\ref{fjhkds} via explicit calculations
 for $n=2,3,4$. 
\subsubsection*{Second cumulant recursion}
By Proposition~\ref{433} and the joint cumulant inversion relation, we have
the second cumulant recursion 
\begin{align} 
  \nonumber
  & 
  c_{k+1,2}^{(\lambda)} (\tau_1 ; \tau_2 )
  = 
  m_{k+1,2}^{(\lambda)} (\tau_1,\tau_2) - 
  m_{k+1,1}^{(\lambda)} (\tau_1) m_{k+1,1}^{(\lambda)} (\tau_2) 
    \\
\nonumber
    & = 
      \int_0^{\tau_1}
  m_{k,2}^{(\lambda)} (u_1, u_1 ) 
  \lambda_k (du_1) 
 +  \int_0^{\tau_1} 
   \int_0^{\tau_2} 
  m_{k,2}^{(\lambda)} (u_1, u_2 ) 
  \lambda_k (du_1) \lambda_k (du_2)
  \\
\nonumber
     & \quad - 
 \int_0^{\tau_1} m_{k+1,1}^{(\lambda)} (s_1) ds_1
\int_0^{\tau_2} m_{k+1,1}^{(\lambda)} (s_2) ds_2
  \\
\nonumber
    & = 
      \int_0^{\tau_1}
  m_{k,2}^{(\lambda)} (u_1, u_1 ) 
  \lambda_k (du_1) 
 + 
   \int_0^{\tau_1} 
   \int_0^{\tau_2} 
(   m_{k,2}^{(\lambda)} (u_1, u_2 ) 
   -
    m_{k,1}^{(\lambda)} (u_1) m_{k,1}^{(\lambda)} (u_2) 
    )
    \lambda_k (du_1) \lambda_k (du_2)
\\ 
\nonumber 
& = 
  \int_0^{\tau_1}
 c_{k,2}^{(\lambda)} (u_1 , u_1) 
  \lambda_k (du_1) 
 + 
   \int_0^{\tau_1} 
   \int_0^{\tau_2} 
 c_{k,2}^{(\lambda)} ( u_1 ; u_2 ) 
 \lambda_k (du_1) \lambda_k (du_2). 
\end{align}
\noindent
\subsubsection*{Third cumulant recursion}
At the third order, we have the next cumulant expression. 
\begin{prop}
 For $k\geq 1$ we have the third cumulant recursion 
\begin{align} 
\nonumber
  & 
    c_{k+1,3}^{(\lambda)} (\tau_1 ; \tau_2 ; \tau_3) 
 =  \int_0^{\tau_1} 
  c_{k,3}^{(\lambda)} (u_1  , u_1, u_1 ) 
  \lambda_k (du_1) 
  + 
  \int_0^{\tau_1} 
  \int_0^{\tau_3} 
  c_{k,3}^{(\lambda)} (u_1,u_1 ;u_3) 
  \lambda_k (du_1) \lambda_k (du_3)
    \\
\nonumber
    &  
  + 
  \int_0^{\tau_1} 
  \int_0^{\tau_2} 
  c_{k,3}^{(\lambda)} (u_1,u_1;u_2) 
  \lambda_k (du_1) \lambda_k (du_2)
    + 
  \int_0^{\tau_1}
  \int_0^{\tau_2} 
  c_{k,3}^{(\lambda)} (u_1;u_2,u_2) 
  \lambda_k (du_1) (du_2)
  \\
\nonumber 
    &  
 + 
  \int_0^{\tau_1} 
  \int_0^{\tau_2} 
  \int_0^{\tau_3} 
  c_{k,3}^{(\lambda)} (u_1 ; u_2;u_3) 
  \lambda_k (du_1) \lambda_k (du_2) \lambda_k (du_3). 
\end{align} 
\end{prop}
\begin{Proof} 
  By the joint cumulant-moment relationship 
\begin{align} 
\label{fjkdsf12}
 & 
 c_{k+1,3}^{(\lambda)} (\tau_1 ; \tau_2 ; \tau_3) 
 = 
 \sum_{l=1}^n
  (l-1)!
  (-1)^{l-1}
  \sum_{\pi_1\cup \cdots \cup \pi_l = \{1,\ldots , n\}}
  \prod_{j=1}^l 
  m_{k+1,|\pi_j|}^{(\lambda)} (\widebar{\tau}_{\pi_j})
  \\
  \nonumber
  & = 
  m_{k+1,3}^{(\lambda)} (\tau_1,\tau_2,\tau_3) 
  - 
  \sum_{\pi_1\cup \pi_2 = \{ 1,2,3\} }
  m_{k+1,|\pi_1|}^{(\lambda)} (\widebar{\tau}_{\pi_1})
  m_{k+1,|\pi_2|}^{(\lambda)} (\widebar{\tau}_{\pi_2})
  +
 2 m_{k+1,1}^{(\lambda)} (\tau_1) 
 m_{k+1,1}^{(\lambda)} (\tau_2) 
 m_{k+1,1}^{(\lambda)} (\tau_3)
\end{align}
 and Proposition~\ref{433}, we have 
 \begin{align}
   \nonumber 
   & c_{k+1,3}^{(\lambda)} (\tau_1 ; \tau_2 ; \tau_3)
   \\
   \nonumber
   &
   = 
  \int_0^{ \tau_1 \wedge \tau_2 \wedge  \tau_3} 
  m_{k,3}^{(\lambda)} (u_1, u_1,u_1 ) 
  \lambda_k (du_1) 
  + 
  \int_0^{\tau_1 \wedge \tau_2} 
  \int_0^{\tau_3} 
  m_{k,3}^{(\lambda)} (u_1,u_1,u_3) 
  \lambda_k (du_3) \lambda_k (du_1)
  \\
   \nonumber 
  & 
  + 
  \int_0^{\tau_1 \wedge \tau_3} 
  \int_0^{\tau_2} 
  m_{k,3}^{(\lambda)} (u_1,u_1,u_2) 
  \lambda_k (du_2) \lambda_k (du_1)
  + 
  \int_0^{\tau_1}
  \int_0^{\tau_2 \wedge \tau_3} 
  m_{k,3}^{(\lambda)} (u_1,u_2,u2) 
  \lambda_k (du_2) \lambda_k (du_1) 
  \\
   \nonumber 
  &  + 
  \int_0^{\tau_3} 
  \int_0^{\tau_2} 
  \int_0^{\tau_1} 
  m_{k,3}^{(\lambda)} (u_1,u_2,u_3) 
  \lambda_k (du_1) \lambda_k (du_2) \lambda_k (du_3) 
  \\
   \nonumber 
  & 
  - 
  \int_0^{\tau_1}
  m_{k,1}^{(\lambda)} (u_1 ) 
  \lambda_k (du_1)
  \left(
      \int_0^{\tau_2 \wedge \tau_3}
  m_{k,2}^{(\lambda)} (u_2, u_2 ) 
  \lambda_k (du_2) 
 + 
   \int_0^{\tau_3} 
   \int_0^{\tau_2} 
  m_{k,2}^{(\lambda)} (u_2, u_3 ) 
  \lambda_k (du_2) \lambda_k (du_3)
  \right)
  \\
   \nonumber 
  & 
  - 
  \int_0^{\tau_2}
  m_{k,1}^{(\lambda)} (u_2 ) 
  \lambda_k (du_2)
  \left(
      \int_0^{\tau_1 \wedge \tau_3}
  m_{k,2}^{(\lambda)} (u_1, u_1 )  
  \lambda_k (du_1) 
 + 
   \int_0^{\tau_3} 
   \int_0^{\tau_1} 
  m_{k,2}^{(\lambda)} (u_1, u_3 ) 
  \lambda_k (du_1) \lambda_k (du_3)
  \right)
  \\
   \nonumber 
  & 
  - 
  \int_0^{\tau_3}
  m_{k,1}^{(\lambda)} (u_3 ) 
  \lambda_k (du_3)
  \left(
      \int_0^{\tau_1 \wedge \tau_2}
  m_{k,2}^{(\lambda)} (u_1, u_1 ) 
  \lambda_k (du_1) 
 + 
   \int_0^{\tau_2} 
   \int_0^{\tau_1} 
  m_{k,2}^{(\lambda)} (u_1, u_2 ) 
  \lambda_k (du_1) \lambda_k (du_2)
  \right)
  \\
  \label{fjklsdf}
  &  + 2 
  \int_0^{\tau_1}
  m_{k,1}^{(\lambda)} (u_1 ) 
  \lambda_k (du_1)
    \int_0^{\tau_2}
  m_{k,1}^{(\lambda)} (u_2 ) 
  \lambda_k (du_2)
    \int_0^{\tau_3}
  m_{k,1}^{(\lambda)} (u_3 ) 
  \lambda_k (du_3). 
\end{align}
 Applying \eqref{fjkdsf12} to $k$-hops, i.e. 
 \begin{eqnarray*}
  c_{k,3}^{(\lambda)} (u_1 ; u_2 ; u_3) 
   & = & 
  m_{k,3}^{(\lambda)} (u_1,u_2 ; u_3) 
  \\
   & &  - 
   m_{k,1}^{(\lambda)} (u_1 ) 
   m_{k,2}^{(\lambda)} (u_1 ;  u_2 ) 
    - 
   m_{k,1}^{(\lambda)} (u_1 ) 
   m_{k,2}^{(\lambda)} (u_1 ;  u_2 ) 
  - 
   m_{k,1}^{(\lambda)} (u_1 ) 
   m_{k,2}^{(\lambda)} (u_1 ;  u_2 ) 
   \\
   & &  + 2 
   m_{k,1}^{(\lambda)} (u_1 ) 
   m_{k,1}^{(\lambda)} (u_2 ) 
   m_{k,1}^{(\lambda)} (u_3 ),  
\end{eqnarray*}
 allows us to simplify \eqref{fjklsdf} to
\begin{align*}
   \nonumber 
 & 
     c_{k+1,3}^{(\lambda)} (\tau_1 ; \tau_2 ; \tau_3)
     \\
      & = 
  \int_0^{ \tau_1 \wedge \tau_2 \wedge  \tau_3} 
  m_{k,3}^{(\lambda)} (u_1, u_1,u_1 ) 
  \lambda_k (du_1) 
   + 
  \int_0^{\tau_1 \wedge \tau_2} 
  \int_0^{\tau_3} 
  m_{k,3}^{(\lambda)} (u_1,u_1,u_3) 
  \lambda_k (du_3) \lambda_k (du_1) 
  \\
  \nonumber 
  &   
   + 
  \int_0^{\tau_1 \wedge \tau_3} 
  \int_0^{\tau_2} 
  m_{k,3}^{(\lambda)} (u_1,u_1,u_2) 
  \lambda_k (du_2) \lambda_k (du_1) 
   + 
  \int_0^{\tau_1}
  \int_0^{\tau_2 \wedge \tau_3} 
  m_{k,3}^{(\lambda)} (u_1,u_2,u_2) 
  \lambda_k(du_2) \lambda_k (du_1) 
    \\
   \nonumber 
  &  
     + 
  \int_0^{\tau_1} 
  \int_0^{\tau_2} 
  \int_0^{\tau_3} 
  c_{k,3}^{(\lambda)} (u_1 ; u_2 ; u_3) 
  \lambda_k (du_1) \lambda_k (du_2) \lambda_k (du_3) 
  \\
   \nonumber 
  &  
     - 
  \int_0^{\tau_1}
  m_{k,1}^{(\lambda)} (u_1 ) 
  \lambda_k (du_1)
      \int_0^{\tau_2 \wedge \tau_3}
  m_{k,2}^{(\lambda)} (u_2, u_2 ) 
  \lambda_k (du_2) 
  \\
   \nonumber 
  &  
     - 
  \int_0^{\tau_2}
  m_{k,1}^{(\lambda)} (u_2 ) 
  \lambda_k (du_2)
      \int_0^{\tau_1 \wedge \tau_3}
  m_{k,2}^{(\lambda)} (u_1, u_1 ) 
  \lambda_k (du_1) 
  \\
   \nonumber 
  &  
     - 
  \int_0^{\tau_3}
  m_{k,1}^{(\lambda)} (u_3 ) 
  \lambda_k (du_3)
      \int_0^{\tau_1 \wedge \tau_2}
  m_{k,2}^{(\lambda)} (u_1, u_1 ) 
  \lambda_k (du_1)
  \\
  & =   
  \int_0^{ \tau_1 \wedge \tau_2 \wedge  \tau_3} 
  m_{k,3}^{(\lambda)} (u_1, u_1,u_1 ) 
  \lambda_k (du_1) 
     + 
  \int_0^{\tau_1 \wedge \tau_2} 
  \int_0^{\tau_3} 
  c_{k,3}^{(\lambda)} (u_1,u_1 ; u_3 ) 
  \lambda_k (du_3) \lambda_k (du_1) 
  \\
   \nonumber 
  &  
     + 
  \int_0^{\tau_1 \wedge \tau_3} 
  \int_0^{\tau_2} 
  c_{k,3}^{(\lambda)} (u_1,u_1 ; u_2) 
  \lambda_k (du_2) \lambda_k (du_1) 
     + 
  \int_0^{\tau_1}
  \int_0^{\tau_2 \wedge \tau_3} 
  c_{k,3}^{(\lambda)} (u_1 ; u_2,u_2) 
  \lambda_k(du_2) \lambda_k (du_1) 
    \\
   \nonumber 
  &  
    + 
  \int_0^{\tau_1} 
  \int_0^{\tau_2} 
  \int_0^{\tau_3} 
  c_{k,3}^{(\lambda)} (u_1 ; u_2 ; u_3) 
  \lambda_k (du_1) \lambda_k (du_2) \lambda_k (du_3).
\end{align*}
\end{Proof} 

\subsubsection*{Fourth cumulant recursion}
Taking $\tau_1=\tau_2=\tau_3=\tau_4=\tau$ for simplicity,
we have 
\begin{align*} 
 & 
 c_{k+1,4}^{(\lambda )} (\tau ;\tau ;\tau ;\tau ) 
 = 
 \sum_{l=1}^n
  (l-1)!
  (-1)^{l-1}
  \sum_{\pi_1\cup \cdots \cup \pi_l = \{1,\ldots , n\}}
  \prod_{j=1}^l 
  m_{k+1,|\pi_j|}^{(\lambda)} (\widehat{s}_{\pi_j})
    \\
& =   
  m_{k+1,4}^{(\lambda)} (\tau ,\tau ,\tau ,\tau ) 
  - 
  4 m_{k+1,1}^{(\lambda)} (\tau )m_{k+1,3}^{(\lambda)} (\tau ,\tau ,\tau ) 
  - 
  3 m_{k+1,2}^{(\lambda)} (\tau ,\tau )m_{k+1,2}^{(\lambda)} (\tau ,\tau ) 
  \\
  &  
  + 
  6 m_{k+1,1}^{(\lambda)} (\tau )m_{k+1,1}^{(\lambda)} (\tau )m_{k+1,2}^{(\lambda)} (\tau ,\tau ) 
  - m_{k+1,1}^{(\lambda)} (\tau )m_{k+1,1}^{(\lambda)} (\tau )m_{k+1,1}^{(\lambda)} (\tau ) m_{k+1,1}^{(\lambda)} (\tau )
\\
& =   
  \int_0^{\tau } 
    m_{k,4}^{(\lambda)} (u_1,u_1,u_1,u_1) 
  \lambda_k (du_1) 
  +
  4 \int_0^{\tau } 
  \int_0^{\tau }
  m_{k,4}^{(\lambda)} (u_1,u_2,u_2,u_2) 
  du_1du_2 
  \\
  & 
  +
  3 \int_0^{\tau } 
  \int_0^{\tau } 
  m_{k,4}^{(\lambda)} (u_1,u_1,u_2,u_2) 
  du_1du_2 
\\
  & 
  +
  6 \int_0^\tau  
  \int_0^\tau 
  \int_0^{\tau }
  m_{k,4}^{(\lambda)} (u_1,u_2,u_3,u_3) 
  du_1du_2du_3
    +
  \int_0^\tau  
  \int_0^\tau  
  \int_0^\tau  
  \int_0^\tau  
  m_{k,4}^{(\lambda)} (u_1,u_2,u_3,u_4) 
  du_1du_2du_3du_4
  \\
  &  
  - 
  4 m_{k+1,1}^{(\lambda)} (\tau )m_{k+1,3}^{(\lambda)} (\tau ,\tau ,\tau ) 
  - 
  3 m_{k+1,2}^{(\lambda)} (\tau ,\tau )m_{k+1,2}^{(\lambda)} (\tau ,\tau ) 
  \\
  &  
  + 
  6 m_{k+1,1}^{(\lambda)} (\tau )m_{k+1,1}^{(\lambda)} (\tau )m_{k+1,2}^{(\lambda)} (\tau ,\tau ) 
  - m_{k+1,1}^{(\lambda)} (\tau )m_{k+1,1}^{(\lambda)} (\tau )m_{k+1,1}^{(\lambda)} (\tau )m_{k+1,1}^{(\lambda)} (\tau )
  \\
& =   
  \int_0^{\tau } 
  m_{k,4}^{(\lambda)} (u_1,u_1,u_1,u_1) 
  \lambda_k (du_1) 
  +
  4 \int_0^{\tau } 
  \int_0^{\tau }
  m_{k,4}^{(\lambda)} (u_1,u_2,u_2,u_2) 
  du_1du_2 
  \\
  & 
  +
  3 \int_0^{\tau } 
  \int_0^{\tau } 
  m_{k,4}^{(\lambda)} (u_1,u_1,u_2,u_2) 
  du_1du_2 
\\
  & 
  +
  6 \int_0^\tau  
  \int_0^\tau 
  \int_0^{\tau }
  m_{k,4}^{(\lambda)} (u_1,u_2,u_3,u_3) 
  du_1du_2du_3
    +
  \int_0^\tau  
  \int_0^\tau  
  \int_0^\tau  
  \int_0^\tau  
  m_{k,4}^{(\lambda)} (u_1,u_2,u_3,u_4) 
  du_1du_2du_3du_4
  \\
  &  
  - 
  4 \int_0^\tau m_{k,1}^{(\lambda)} (u) du
  \int_0^\tau m_{k,3}^{(\lambda)} (u_1,u_1,u_1)du_1 
  - 
  12 \int_0^\tau m_{k,1}^{(\lambda)} (u) du
  \int_0^\tau \int_0^\tau m_{k,3}^{(\lambda)} (u_1,u_2,u_2)du_1 du_2 
\\
  & 
  - 
  4 \int_0^\tau m_{k,1}^{(\lambda)} (u) du
  \int_0^\tau \int_0^\tau \int_0^\tau m_{k,3}^{(\lambda)} (u_1,u_2,u_3)du_1 du_2 du_3 
  - 
  3
  \left(
  \int_0^\tau m_{k,2}^{(\lambda)} (u_1,u_1)du_1
  \right)^2
\\
  & 
  - 6 
  \int_0^\tau m_{k,2}^{(\lambda)} (u_1,u_1)du_1 
  \int_0^\tau \int_0^\tau m_{k,2}^{(\lambda)} (u_1,u_2)du_1 du_2 
  - 
  3
  \left(
  \int_0^\tau \int_0^\tau m_{k,2}^{(\lambda)} (u_1,u_2)du_1 du_2 
  \right)^2 
  \\
  &  
  + 
  6
  \left( \int_0^\tau m_{k,1}^{(\lambda)} (u) du \right)^2 
  \left(
  \int_0^\tau m_{k,2}^{(\lambda)} (u_1,u_1)du_1 
  + \int_0^\tau \int_0^\tau m_{k,2}^{(\lambda)} (u_1,u_2)du_1 du_2 
  \right)
\\
  & 
  - m_{k+1,1}^{(\lambda)} (\tau )m_{k+1,1}^{(\lambda)} (\tau )m_{k+1,1}^{(\lambda)} (\tau )m_{k+1,1}^{(\lambda)} (\tau )
 \\
& =   
  \int_0^{\tau } 
  m_{k,4}^{(\lambda)} (u_1,u_1,u_1,u_1) 
  \lambda_k (du_1) 
  +
  4 \int_0^{\tau } 
  \int_0^{\tau }
  m_{k,4}^{(\lambda)} (u_1,u_2,u_2,u_2) 
  du_1du_2 
  \\
  &
  +
  3 \int_0^{\tau } 
  \int_0^{\tau } 
  m_{k,4}^{(\lambda)} (u_1,u_1,u_2,u_2) 
  du_1du_2 
\\
  & 
  +
  6 \int_0^\tau  
  \int_0^\tau 
  \int_0^{\tau }
  m_{k,4}^{(\lambda)} (u_1,u_2,u_3,u_3) 
  du_1du_2du_3
    +
  \int_0^\tau  
  \int_0^\tau  
  \int_0^\tau  
  \int_0^\tau  
  c_{k,4}^{(\lambda)} (u_1 ; u_2 ; u_3 ; u_4) 
  du_1du_2du_3du_4
  \\
  &  
  - 
  4 \int_0^\tau m_{k,1}^{(\lambda)} (u) du
  \int_0^\tau m_{k,3}^{(\lambda)} (u_1,u_1,u_1)du_1 
  - 
  12 \int_0^\tau m_{k,1}^{(\lambda)} (u) du
  \int_0^\tau \int_0^\tau m_{k,3}^{(\lambda)} (u_1,u_2,u_2)du_1 du_2 
\\
  & 
  - 
  3
  \left(
  \int_0^\tau m_{k,2}^{(\lambda)} (u_1,u_1)du_1
  \right)^2
  - 6 
  \int_0^\tau m_{k,2}^{(\lambda)} (u_1,u_1)du_1 
  \int_0^\tau \int_0^\tau m_{k,2}^{(\lambda)} (u_1,u_2)du_1 du_2 
\\
  & 
  + 
  6
  \left( \int_0^\tau m_{k,1}^{(\lambda)} (u) du \right)^2 
  \int_0^\tau m_{k,2}^{(\lambda)} (u_1,u_1)du_1 
\\
& =   
  \int_0^{\tau } 
  m_{k,4}^{(\lambda)} (u_1,u_1,u_1,u_1) 
  \lambda_k (du_1) 
  +
  4 \int_0^{\tau } 
  \int_0^{\tau }
  c_{k,4}^{(\lambda)} (u_1 ; u_2,u_2,u_2) 
  du_1du_2 
  \\
  &
  +
  3 \int_0^{\tau } 
  \int_0^{\tau } 
  c_{k,4}^{(\lambda)} (u_1,u_1 ; u_2,u_2) 
  du_1du_2 
\\
  & 
  +
  6 \int_0^\tau  
  \int_0^\tau 
  \int_0^{\tau }
  m_{k,4}^{(\lambda)} (u_1,u_2,u_3,u_3) 
  du_1du_2du_3
    +
  \int_0^\tau  
  \int_0^\tau  
  \int_0^\tau  
  \int_0^\tau  
  c_{k,4}^{(\lambda)} (u_1 ; u_2 ; u_3 ; u_4) 
  du_1du_2du_3du_4
  \\
  &  
  - 
  12 \int_0^\tau m_{k,1}^{(\lambda)} (u) du
  \int_0^\tau \int_0^\tau m_{k,3}^{(\lambda)} (u_1,u_2,u_2)du_1 du_2 
\\
  & 
  - 6 
  \int_0^\tau m_{k,2}^{(\lambda)} (u_1,u_1)du_1 
  \int_0^\tau \int_0^\tau m_{k,2}^{(\lambda)} (u_1,u_2)du_1 du_2 
\\
  & 
  + 
  6
  \left( \int_0^\tau m_{k,1}^{(\lambda)} (u) du \right)^2 
  \int_0^\tau m_{k,2}^{(\lambda)} (u_1,u_1)du_1 
\\
& =   
  \int_0^{\tau } 
  c_{k,4}^{(\lambda)} (u_1,u_1,u_1,u_1) 
  \lambda_k (du_1) 
  +
  4 \int_0^{\tau } 
  \int_0^{\tau }
  c_{k,4}^{(\lambda)} (u_1 ; u_2,u_2,u_2) 
  du_1du_2 
  \\
  &
  +
  3 \int_0^{\tau } 
  \int_0^{\tau } 
  c_{k,4}^{(\lambda)} (u_1,u_1 ; u_2,u_2) 
  du_1du_2 
\\
  & 
  +
  6 \int_0^\tau  
  \int_0^\tau 
  \int_0^{\tau }
  c_{k,4}^{(\lambda)} (u_1 ; u_2 ; u_3,u_3) 
  du_1du_2du_3
    +
  \int_0^\tau  
  \int_0^\tau  
  \int_0^\tau  
  \int_0^\tau  
  c_{k,4}^{(\lambda)} (u_1 ; u_2 ; u_3 ; u_4) 
  du_1du_2du_3du_4
. 
\end{align*} 
 In order to reach we above conclusion we used the relations 
\begin{eqnarray*} 
  \lefteqn{
 4 
\int_0^\tau 
\int_0^\tau 
 c_{k,4}^{(\lambda)} (u_1; u_2,u_2,u_2) 
 du_1 du_2 
  }
  \\
  & = &
 4 
\int_0^\tau 
\int_0^\tau 
 m_{k,4}^{(\lambda)} (u_1, u_2,u_2,u_2 ) 
 du_1 du_2 
 -
 4 
\int_0^\tau 
\int_0^\tau 
 m_{k,1}^{(\lambda)} (u_1)
 m_{k,3}^{(\lambda)} ( u_2,u_2,u_2) 
 du_1 du_2,  
\end{eqnarray*} 
\begin{eqnarray*} 
  \lefteqn{
 3 
\int_0^\tau 
\int_0^\tau 
 c_{k,4}^{(\lambda)} (u_1,u_1 ; u_3,u_3 ) 
 du_1 du_2 
  }
  \\
  & = &
 3 
\int_0^\tau 
\int_0^\tau 
 m_{k,4}^{(\lambda)} (u_1,u_1 , u_3,u_3 ) 
 du_1 du_2 
 - 
 3 
\int_0^\tau 
\int_0^\tau 
 m_{k,2}^{(\lambda)} (u_1,u_2 )
 m_{k,2}^{(\lambda)} (u_3,u_3 ) 
 du_1 du_3,  
\end{eqnarray*} 
 and 
\begin{eqnarray*} 
  \lefteqn{
    6 
 \int_0^\tau 
\int_0^\tau 
\int_0^\tau 
 c_{k,4}^{(\lambda)} (u_1 ; u_2 ;u_3, u_3 ) 
 du_1 du_2 du_3
  }
  \\
  & = &
 6 
 \int_0^\tau 
\int_0^\tau 
\int_0^\tau 
 m_{k,4}^{(\lambda)} (u_1, u_2, u_3,u_3 ) 
 du_1 du_2 du_3
 - 
 12 
 \int_0^\tau 
 m_{k,1}^{(\lambda)} (u_1)
 du_1
 \int_0^\tau 
\int_0^\tau 
m_{k,3}^{(\lambda)} (u_2, u_3,u_3 )
du_2 du_3
\\
& &
- 
 6 
\int_0^\tau 
m_{k,2}^{(\lambda)} (u_3, u_3 )
du_3
\int_0^\tau 
\int_0^\tau 
m_{k,2}^{(\lambda)} (u_1, u_2) 
 du_1 du_2 
\\
 & &
 + 6 
 \int_0^\tau 
 \int_0^\tau 
 \int_0^\tau 
 m_{k,1}^{(\lambda)} (u_1) 
 m_{k,1}^{(\lambda)} (u_2) 
 m_{k,2}^{(\lambda)} (u_3,u_3) 
 du_1 du_2 du_3. 
\end{eqnarray*}

\small 

\subsubsection*{Acknowledgement}
I thank A.P. Kartun-Giles for useful discussions.

\footnotesize 

\def\cprime{$'$} \def\polhk#1{\setbox0=\hbox{#1}{\ooalign{\hidewidth
  \lower1.5ex\hbox{`}\hidewidth\crcr\unhbox0}}}
  \def\polhk#1{\setbox0=\hbox{#1}{\ooalign{\hidewidth
  \lower1.5ex\hbox{`}\hidewidth\crcr\unhbox0}}} \def\cprime{$'$}

\end{document}